\UseRawInputEncoding
\pdfoutput=1
\RequirePackage{fix-cm}
\documentclass[11pt]{article}

\usepackage{graphicx,amsmath,amsfonts,amscd,amssymb,bm,cite,epsfig,epsf,url,color,algorithm,algorithmic,enumerate}
\usepackage{fullpage}
\usepackage[small,bf]{caption}
\usepackage{subfigure}
\usepackage{lscape}
\usepackage[space]{grffile}
\usepackage{epsfig,psfrag,amstext,subfigure,subeqnarray,nccfoots,color,multirow,bm}
\usepackage{appendix}

\usepackage{amsfonts,amssymb,amsbsy,amsmath,paralist,theorem,bm,ifthen,color,xcolor}
\usepackage{graphicx}
\usepackage{psfrag,epsfig,subfigure,amsmath}
\usepackage[justification=centering]{caption}
\usepackage{multirow}

\usepackage{url}
\usepackage{subfigure,relsize}
\usepackage{cases,booktabs}
\usepackage{algorithmic,algorithm}
\usepackage{array,threeparttable}

\usepackage{pifont}       
\usepackage{bbding}       

\newcommand{\xmark}{\ding{55}}

\newcommand\Ec{\ensuremath{\mathcal{E}}}
\newcommand\Gc{\ensuremath{\mathcal{G}}}

\newcommand\Vc{\ensuremath{{\mathcal{V}}}}

\newcommand\xb{\ensuremath{{\bm x}}}

\newcommand\yb{\ensuremath{{\bm y}}}

\newcommand\ub{\ensuremath{{\bm u}}}

\newcommand\Ab{\ensuremath{{\bm A}}}

\newcommand\ab{\ensuremath{{\bm a}}}
\newcommand\Bb{\ensuremath{{\bm B}}}
\newcommand\bb{\ensuremath{{\bm b}}}
\newcommand\Cb{\ensuremath{{\bm C}}}

\newcommand\eb{\ensuremath{{\bm e}}}

\newcommand\Gb{\ensuremath{{\bm G}}}
\newcommand\gb{\ensuremath{{\bm g}}}
\newcommand\Ib{\ensuremath{{\bm I}}}

\newcommand\bsb{\ensuremath{{\bm s}}}

\newcommand\Xb{\ensuremath{{\bm X}}}
\newcommand\Yb{\ensuremath{{\bm Y}}}
\newcommand\Ub{\ensuremath{{\bm U}}}

\newcommand\Vb{\ensuremath{{\bm V}}}
\newcommand\vb{\ensuremath{{\bm v}}}
\newcommand\Wb{\ensuremath{{\bm W}}}
\newcommand\Zb{\ensuremath{{\bm Z}}}
\newcommand\zb{\ensuremath{{\bm z}}}

\newcommand\thetab{\ensuremath{{\bm \theta}}}
\newcommand\Thetab{\ensuremath{{\bm \Theta}}}

\newcommand\zerob{\ensuremath{{\bm 0}}}

\newcommand\phib{\ensuremath{{\bm \phi}}}

\newtheorem{Lemma}{Lemma}

\newtheorem{Theorem}{Theorem}

\newtheorem{Corollary}{Corollary}

\newtheorem{Remark}{Remark}
\newtheorem{Assumption}{Assumption}

\usepackage{geometry}
\allowdisplaybreaks[4]
\usepackage{epstopdf}
\graphicspath{{fig/}}

\begin{document}

\title{Gradient and Variable Tracking with Multiple Local SGD for Decentralized Non-Convex Learning}

\author{Songyang Ge and Tsung-Hui Chang
    \thanks{\smaller[1]	Songyang Ge and  Tsung-Hui Chang are with the Shenzhen Research Institute of Big Data, the School of Science and Engineering, The Chinese University of Hong Kong, Shenzhen, China 518172.
   E-mail: songyangge@link.cuhk.edu.cn; tsunghui.chang@ieee.org}
   \thanks{\smaller[1] Tsung-Hui Chang is the corresponding author.}}
\maketitle	

\begin{abstract}
    Stochastic distributed optimization methods that solve an optimization problem over a multi-agent network have played an important role in a variety of large-scale signal processing and machine leaning applications. Among the existing methods, the gradient tracking (GT) method is found robust against the variance between agents' local data distribution, in contrast to the distributed stochastic gradient descent (SGD) methods which have a slowed convergence speed when the agents have heterogeneous data distributions. However, the GT method can be communication expensive due to the need of a large number of iterations for convergence.  In this paper, we intend to reduce the communication cost of the GT method by integrating it with the local SGD technique. Specifically, we propose a new local stochastic GT (LSGT) algorithm where, within each communication round, the agents perform multiple SGD updates locally. Theoretically, we build the convergence conditions of the LSGT algorithm and show that it can have an improved convergence rate of $\mathcal{O}(1/\sqrt{ET})$, where $E$ is the number of local SGD updates and $T$ is the number of communication rounds. We further extend the LSGT algorithm to solve a more complex learning problem which has linearly coupled variables inside the objective function. Experiment results demonstrate that the proposed algorithms have significantly improved convergence speed even under heterogeneous data distribution.
\end{abstract}


\section{Introduction}
\label{sec:intro}

With the continuous acceleration of the digitization process of human society,
massive amounts of data are growing explosively. This phenomenon makes it difficult to store all the data in one device and process it by one processor.
In view of this,
distributed optimization methods arouse significant attention in various signal processing and machine learning fields, see, e.g.,  \cite{giannakis2016decentralized,BK:Bekkerman12,scutari2018parallel,chang2020distributed,zhang2021decentralized,yang2019survey,nedic2018distributed}. 
In particular, by considering a multi-agent network with $N$ agents linked via a connected graph, most of the existing distributed optimization methods focus on solving the following problem
\begin{align}\label{prob: HL 2}
    \min_{\yb \in \mathbb{R}^{p} } F(\yb) \triangleq
   \frac{1}{N} \sum_{n=1}^N f_n (\yb ),
\end{align}
where $\yb \in \mathbb{R}^p$ is the parameter vector to optimize and each $f_n: \mathbb{R}^p \!\rightarrow\! \mathbb{R}$ is a smooth and possibly non-convex local cost function of agent $n$, $\forall n\in[N]\triangleq \{1,\ldots,N\}$.
For a statistical learning problem, one may assume $f_n(\yb) = \mathbb{E}_{\xi\sim \mathcal{D}_n}[\ell_n(\yb, \xi)]$ where $\ell_n$ is a loss function of $\yb$, and the data sample $\xi$ is randomly drawn from a local dataset $\mathcal{D}_n$.

\vspace{-0.3cm}
\subsection{Literature Review}
There exist rich results in the literature \cite{chang2020distributed,yang2019survey,nedic2018distributed} for solving the nonconvex distributed problem \eqref{prob: HL 2}.
With the consideration of the full gradient $\nabla f_n$, some deterministic first-order methods have been developed.
The decentralized gradient descent (DGD) algorithm \cite{lian2017can} uses consensus gradient descent (GD) via a proper mixing matrix and a diminishing stepsize to achieve a stationary solution of problem \eqref{prob: HL 2}.
Primal-dual based methods, such as the the proximal gradient primal-dual algorithm (Prox-GPDA) \cite{hong2017prox} can converge to a stationary solution with a constant stepsize. However,  Prox-GPDA
can only achieve the communication and computation lower bounds in star and fully connected networks. With regard to the general network topology, a near-optimal scheme, xFilter, was proposed in \cite{sun2019distributed}. However, owing to the usage of full gradients, the aforementioned methods are suitable only for learning problems with a small-to-moderate dataset.

\begin{table*}
	  \centering
	  \caption{Comparison of different algorithms}
\tabcolsep=0.09cm
\renewcommand\arraystretch{2}
	  \label{table1}
	  \begin{tabular}{l c c c c c c c}
	  	\hline
	  	Algorithm &  function   & gradient & stepsize & comp. & comm. &  $\lambda_w$ &  $\Wb$   \\
	  	\hline
	  	GT \cite{qu2017harnessing} & ns-cvx. & full & cst. & $\textstyle \mathcal{O} (\frac{1}{\epsilon} )$ & $\textstyle \mathcal{O} (\frac{1}{\epsilon} )$ & $(-1, 1)$& ds. \\
	  	DSGT \cite{pu2021distributed} & cvx.  & stochastic & dimi. & $\textstyle \mathcal{O} (\frac{ 1}{\epsilon} )$ & $\textstyle \mathcal{O} (\frac{1}{\epsilon} )$ & $(-1, 1)$& ds. \\
	  	GNSD \cite{lu2019gnsd} & n-cvx. & stochastic & cst. & $\textstyle \mathcal{O} (\frac{ 1}{\epsilon^2} )$ & $\textstyle \mathcal{O} (\frac{1}{\epsilon^2} )$& $(-1, 1)$& ds. and sym.  \\
        GT-DSGD \cite{xin2021improved} & n-cvx. & stochastic, & cst. & $\textstyle \mathcal{O} (\frac{ 1}{\epsilon^2} )$
        & $\textstyle \mathcal{O} (\frac{1}{\epsilon^2} )$
        & $(-1, 1)$& ds.  \\
         $\mathrm{D^2}$ \cite{tang2018d} & n-cvx. & stochastic
& cst. & $\textstyle \mathcal{O} (\frac{ 1}{\epsilon^2} )$ & $\textstyle \mathcal{O}  (\frac{1}{\epsilon^2} )$
&$ (\textstyle - \frac{1}{3}, 1)$ &   rs.~and~sym.\\
    LU-GT \cite{nguyen2022performance}& n-cvx. & E-step full
& cst. & \xmark & \xmark
&$ (\textstyle -1, 1)$ &   ds. and sym. \\
	  	\hline
  LSGT (Proposed) & n-cvx. & E-step stochastic
& cst. & $\textstyle \mathcal{O} (\frac{ 1}{\epsilon^2} )$ & $\textstyle \mathcal{O} (\frac{1}{E\epsilon^2} )$
&$ (-1, 1)$ & ds.\\
MUST (Proposed) & n-cvx. & E-step stochastic
&cst. & $\textstyle \mathcal{O} (\frac{ 1}{\epsilon^2} )$ & $\textstyle \mathcal{O} (\frac{1}{E\epsilon^2} )$
&$ (-1, 1)$ & ds.\\
\hline
	  \end{tabular}\\
\vspace{0.2cm}
\noindent
``cvx.", ``n-cvx." and ``ns-cvx." manifests convex, nonconvex, and nonstrongly convex,\\
``comp." denotes  computation complexity,  ``comm." represents communication complexity,\\
 ``cst." and ``dimi." are the abbreviation of constant and diminishing,\\
 ``ds." indicates doubly stochastic, ``rs." shows right stochastic, ``sym." suggests symmetric,\\
 $\Wb$ is the mixing matrix, $\lambda_w$ is the second largest eigenvalue of $\Wb$, $\epsilon$ is the solution accuracy, and $E$ is the local update number.
  \end{table*}

Stochastic distributed methods based on mini-batch stochastic gradients have been considered for reducing the complexity of gradient computation. A stochastic variant of the DGD method, namely, the decentralized stochastic gradient descent (DSGD) algorithm, was proposed in \cite{nedic2017achieving,lian2017can}, by replacing the full gradients with stochastic gradients. Except for the requirement of a diminishing stepsize, it is found in \cite{lian2017can,tang2018d} that the DSGD-type methods are sensitive to the variance of data distribution of the agents. In particular, if the data distributions across the agents are heterogeneous, i.e., the so called non-IID data, the convergence performance of the DSGD algorithm degrades significantly. To alleviate the issue, the authors of \cite{tang2018d} proposed a new algorithm called $\textrm{D}^2$. While $\textrm{D}^2$ is more robust against the heterogeneous data, its convergence relies on a restrictive assumption on the mixing matrix (see \cite{chang2020distributed} and  Sec.~\ref{sec:models}). The primal-dual based stochastic distributed methods, such as the stochastic proximal primal dual algorithm with momentum (SPPDM) in \cite{wang2021distributed}, also exhibit better robustness against the data distribution, but involve more tunable parameters.


Recently, the gradient tracking (GT) methods \cite{qu2017harnessing,di2016next} have been proposed as a promising family of distributed methods for solving problem \eqref{prob: HL 2}. Specifically, the GT method can converge to the neighborhood of a stationary solution with a more relaxed condition on the mixing matrix than $\textrm{D}^2$ and involves one stepsize parameter only. The GT method is based upon a simple but effective idea which introduces auxiliary variables to help the distributed agents track the global gradient via consensus averaging. This enables the agents to imitate the centralized GD scheme and thereby has superior robustness than existing methods. The GT methods based on SGD have also been studied;
for instance,  \cite{pu2021distributed} for the strongly convex problems, and GNSD \cite{lu2019gnsd} and GT-DSGD \cite{xin2019distributed} for nonconvex problems.
However, these stochastic GT methods require frequent message exchanges between the agents. Since the communication is always expensive and constrained by limited bandwidth, it may not be easy to implement the above algorithms in applications with limited communication resources or with stringent delay constraints.

Thus, it is urgent to reduce the communication cost of the stochastic GT method. One possible way is to allow the agents to perform multiple local SGD updates within each communication round \cite{stich2018local,mcmahan2017communication}, in contrast to a single SGD update in existing GT methods. This local SGD technique has mature and successful applications in the federated learning scenarios \cite{li2019convergence} where the distributed agents are coordinated by a central server. Specifically, both theoretical and empirical studies \cite{stich2018local,mcmahan2017communication} have shown that, under proper conditions, the local SGD can reduce the communication cost for $E$ times, where $E$ is the number of the local SGD updates per communication round.
Nevertheless, the local SGD technique has not been thoroughly studied for the GT methods over the fully decentralized network. For example, while the recent work \cite{nguyen2022performance} proposed a GT method where the agents execute multiple local updates in each communication round, it considers the full gradient but not the stochastic gradient. Besides, the convergence analysis therein did not fully characterize the impact of local updates on the convergence speed, in addition that their algorithm requires an extra stepsize and a symmetric mixing matrix for proper convergence.

\vspace{-0.3cm}
\subsection{Contribution}
In this paper, we propose a new distributed algorithm, termed as the {local stochastic GT (LSGT)} algorithm, by incorporating the local SGD technique with the stochastic GT method. The LSGT algorithm neither introduces any extra stepsize nor requires a symmetric mixing matrix. As the major contribution,  we carry out the convergence analysis showing that the LSGT algorithm with a constant stepsize can converge sublinearly to the neighborhood of a stationary solution of problem \eqref{prob: HL 2}. Moreover, we provide the conditions for which the LSGT algorithm can benefit the convergence speedup brought by the local SGD technique, which has never been shown in previous works. Table \ref{table1} summarizes the comparison results between the proposed algorithms with the existing methods. One can see that the proposed LSGT algorithm has a $\textstyle \mathcal{O} (\frac{1}{E\epsilon^2} )$ communication complexity whereas that of the other methods is $\textstyle \mathcal{O} (\frac{1}{\epsilon^2} )$.

As the second contribution, we extend the idea of the LSGT algorithm to handle another type of optimization problems where the agents' local variables are linearly coupled inside the objective function (see \eqref{pro: HFL rewrite} and \eqref{eqn: Bni}). One motivating example is the distributed learning problem over hybrid data, where each distributed agent possesses only a subset of data samples and knows only part of the feature information.  There is scant attention on this challenging learning problem over hybrid data. The recent works \cite{zhang2020hybrid, FedHD2022} have studied such problem in the federated learning network but their algorithms are not applicable to the decentralized network.
Based on a similar idea as the GT methods, we introduce additional auxiliary variables to track the linearly coupled term and develop a new algorithm, called the \emph{Multiple-locally-Updated variable Sum Tracking (MUST)} algorithm, to solve the learning problem in \eqref{pro: HFL rewrite}-\eqref{eqn: Bni}.

Finally, the performance of the proposed LSGT and MUST algorithms are evaluated by extensive numerical experiments.

{\bf Synopsis:} In Sec.~\ref{sec: network problem}, the network model and assumptions are elaborated. For solving problem \eqref{prob: HL 2}, the LSGT algorithm is proposed and its convergence analysis are presented in Sec.~\ref{sec: LSGT alg}. The proof details of main theoretical results for the LSGT algorithm are shown in Sec.~\ref{sec: proof of thm 1}. Besides, extension of the LSGT algorithm for the learning problem over the hybrid data, i.e., the MUST algorithm, is investigated in Sec.~\ref{sec: hybrid}.
Numerical results are given in Sec.~\ref{section: simulation results} and
conclusions are drawn in Sec.~\ref{section: conclusions}.

{\bf Notation: } $\Ib_n$ is the $n$ by $n$ identity matrix, and $\mathbf{1}$ is the all-one vector. $A_{i,j}$ is the $(i,j)$-th element of matrix $\Ab$. $\otimes$ denotes the Kronecker product; $\ab^\top$ and $\Ab^\top$ perspectively represent the transpose operation of vector $\ab$ and matrix $\Ab$;
$\langle \ab, \bb \rangle$ represents the inner product of vectors $\ab$ and $\bb$,
$\|\ab\|$ is the Euclidean norm, and $\|\Ab\| $ represents the largest singular value of $ \Ab$; $\|\Ab\|_F$ denotes the matrix Frobenius norm.

\section{Network Model and Assumptions}\label{sec: network problem}



In this section, we present the network model and assumptions for problem \eqref{prob: HL 2}.

\vspace{-0.23cm}
\subsection{Network Model and Assumptions}
\label{sec:models}

We model the multi-agent network as an undirected graph $\Gc=(\Ec,\Vc)$, where $\Ec$ is the set of edges and $\Vc=[N]=\{1,\ldots,N\}$ is the set of agents.
Each agent owns the a local dataset $\mathcal{D}_n$ and can only communicate and exchange information with its neighbors, i.e., agents $n$ and $m$ can communicate with each other if and only if $(n,m)\in \Ec$.   Moreover, we have the following standard assumptions.
\begin{Assumption}\label{assum network}
The underlying graph is connected.
\end{Assumption}

To enable message exchange between agents,  we define a mixing matrix $ \Wb \in \mathbb{R}^{N\times N}$ with $W_{n,m} >0$, $\forall (n,m)\in \Ec$ and
$W_{n,m} =0$ otherwise. Moreover, we have the following assumption.
\begin{Assumption}\label{assum mixing matrix}
 The mixing matrix $\Wb$ is doubly stochastic satisfying
    \begin{align}
          \Wb \mathbf{1} = \mathbf{1},~~\mathbf{1}^\top \Wb = \mathbf{1}, ~~|\lambda_{w}| < 1 ,
    \end{align}
    where $\lambda_{w} =\|\Wb - \frac{1}{N} \mathbf{1}\mathbf{1}^\top\| $ is the second largest eigenvalue of $\Wb$.
\end{Assumption}
There are many ways to construct such mixing matrix $\Wb$ in Assumption 2  \cite[Remark 2]{lu2019gnsd}; for example the max-degree rule \cite{sayed2014adaptive}.
Note that our assumption on $\Wb$ is weaker than the $\mathrm{D}^2$ algorithm in \cite{tang2018d} which requires  $\lambda_{w}\in \textstyle \left(-\frac{1}{3}, 1\right)$. Besides, we don't assume a symmetric $\Wb$, unlike \cite{nguyen2022performance,tang2018d,lu2019gnsd}.

\subsection{Problem Assumptions}\label{sec: problem formulation}


For ease of algorithm development and convergence analysis, we have the following assumptions for problem \eqref{prob: HL 2}.

\begin{Assumption}
\label{assum: lower bound}
The loss function $F(\yb) \triangleq
  \textstyle \frac{1}{N} \sum_{n=1}^N f_n (\yb )$ is bounded below, denoted by $F(\yb)  \geq \underline{F},~\forall\yb$.
\end{Assumption}

\begin{Assumption} 
\label{assum: Lipschitz}
Each $f_n  $ is smooth and its gradient satisfies
\begin{align}
&\left \| \nabla f_n (\yb  ) - \nabla f_n (\yb' ) \right \|    \leq L   \left \|\yb  - \yb' \right \|,~\forall \yb,\yb',
\end{align}
where $L$ is the Lipschitz constant.
\end{Assumption}

In practice, calculating the full gradient can be computationally expensive. Alternatively, the stochastic gradient based on mini-batch data sampling is often used. Specifically, we define the following stochastic gradient for each agent $n$ as
\begin{align}\label{eq: h def}
\gb_n  \triangleq  \frac{1}{|\mathcal{I}_n|}\sum_{\xi  \in  \mathcal{I}_n }  \nabla \ell_n ( \yb_n; \xi  ) ,
\end{align}
where $\mathcal{I}_n  \subseteq \mathcal{D}_n$ is a randomly chosen mini-batch data set at agent $n$. Without loss of generality, we assume that all agents use the same mini-batch size, denoted as $|\mathcal{I}|$. The following assumption is standard for SGD methods.
\begin{Assumption} 
\label{assum: stochastic} For each agent $n\in[N]$,
 we have 

\noindent$\bullet$ Unbiased gradient:
$
\mathbb{E} [  \gb_n ] =   \nabla f_n(\yb_n  );
$%

\noindent$\bullet$ Uniform bounded variance:
$
\mathbb{E}  [  \|   \gb_n   -  \nabla  f_n(\yb_n  )  \|^2  ] \leq \frac{\sigma^2}{|\mathcal{I}|} .
$
\end{Assumption}
\section{Proposed LSGT Algorithm}\label{sec: LSGT alg}

In this section, we propose a new decentralized algorithm, called \emph{local stochastic gradient tracking (LSGT)}, for solving the non-convex problem \eqref{prob: HL 2}. The proposed LSGT algorithm is based on the celebrated stochastic GT method \cite{lu2019gnsd,xin2019distributed} and incorporates multiple local SGD updates in each communication round for improving the convergence speed. We first review the vanilla GT method, and then present the proposed algorithm and its convergence analysis.


\subsection{Review of Stochastic GT Method}
The GT method is a consensus based method with additional auxiliary variables to track the global gradient of the objective function for imitating the centralized gradient descent method \cite{lu2019gnsd,xin2019distributed}.
Specifically, except for the local variable $\yb_n$,
each agent $n$ possesses another auxiliary variable $\tilde{\vb}_n \in \mathbb{R}^{p}$ to estimate a stochastic approximation of the global gradient $\nabla F(\yb)=\frac{1}{N}\textstyle\sum_{n =1}^N   \nabla f_n (\yb)$. The algorithm starts with an initial $\yb_n^0$ and stochastic gradient $\tilde{\vb}_n^0 =  \gb_n^0$ for each agent $n$. Then, at the $r$-th iteration, all agents $n\in [N]$ perform the following two steps in parallel:
\vspace{-0.0cm}
\begin{subequations}\label{eq: GT}
\begin{align}
& \yb_{n}^{ r +1 }\!=\!\sum _{m =1}^N W_{n,  m}  \yb_m^{r}\!-\!
                 \gamma   \tilde{\vb}_n^r
                 , \label{eq: update y GT}\\
& \tilde{\vb}_{n}^{ r +1 }\!=\!\sum _{m =1}^N W_{n,  m} \tilde{\vb}_m^{r}+     \gb_n^{r+1} -   \gb_n^{r}
                 ,\label{eq: update v GT}
\end{align}
\end{subequations}
where $\gamma>0 $ denotes the stepsize. In \eqref{eq: update y GT}, the agents perform consensus averaging of the local variables from their neighbors followed by gradient descent along the direction of $\tilde{\vb}_n$; in \eqref{eq: update v GT}, the agents use consensus averaging to track the sum of the stochastic gradient $\sum_{n=1}^N \gb_n^{r+1}$. The idea behind the GT method is that, when the variable consensus  is approximately reached, i.e., $\yb_n\approx \frac{1}{N}\sum_{m=1}^N\yb_m$ and $\tilde{\vb}_n\approx \frac{1}{N}\sum_{m=1}^N\tilde{\vb}_m$ $\forall n\in [N]$, $\tilde{\vb}_n$ approximates the global gradient $\nabla F(\yb_n)$ and \eqref{eq: update y GT} is the same as the centralized gradient descent (GD) method.

Theoretically, it has been proved that, compared with DSGD algorithms, the convergence rate of the stochastic GT method does not depend on the variation of the local cost functions among the agents. In particular, the convergence analysis in \cite{lian2017can,tang2018d} shows that the convergence rate of DSGD depends on the bound
$\frac{1}{N}\sum_{n=1}^N\|\nabla f_n(\yb)-\nabla F(\yb)\| \leq \xi$ and can be slowed down if $\xi$ is large, whereas the stochastic GT method in \cite{lu2019gnsd,xin2021improved} does not. This implies that the stochastc GT method is more robust against so called ``heterogeneous data" (where the local datasets of agents have different statistical properties). 



However, like the DSGD methods \cite{nedic2017achieving,lian2017can}, the stochastic GT method in \eqref{eq: GT} requires many communication rounds to converge and thereby still has a large communication overhead. To reduce the communication cost, we propose to incorporate the local SGD technique, that is, allowing the agents to perform multiple steps of SGD in each communication round. The effectiveness of local SGD on speeding up the algorithm convergence has been studied both theoretically and empirically \cite{stich2018local,mcmahan2017communication}, but, except for the recent work \cite{nguyen2022performance}, has not been thoroughly considered for the GT-based method.

In the next subsection, we present a new communication-efficient LSGT algorithm by integrating the stochastic GT method and the local SGD technique, and then build its convergence properties.


\begin{algorithm}[!t]
	\caption{Proposed  LSGT algorithm for solving \eqref{prob: HL 2}}
	\label{alg: LSGT}
	\begin{algorithmic}[1]
	\STATE {\bfseries Initialize:} Let $\yb_1^0=\ldots=\yb_N^0$ and $\vb_n^{0}  =  \gb_n^0$, $\forall n\in[N]$.
	\FOR{communication round $r = 0$ {\bf to} $T$}
	\FOR{agent $n =1$ {\bf to} $N$ in parallel}
         \STATE
     Receive information from  neighbours and set
     \vspace{-0.1cm}
            \begin{align}\label{alg1:consensus aveg}
                \begin{bmatrix}
                \yb_n^{r,0} \\
                \vb_n^{r,0}
                \end{bmatrix}
                =\sum_{m=1}^N W_{n,m}
                 \begin{bmatrix}
                \yb_m^{r} \\
                \vb_m^{r}
                 \end{bmatrix},
                %
                 ~\gb_n^{r,0} \!= \!\gb_n^r.
            \end{align}
             \vspace{-0.0cm}
		\FOR{local update $q = 1,\ldots, E$,}
\STATE
\vspace{-0.2cm} \begin{subequations}\label{eq: LSGT q 2}
              \begin{align}
                & \yb_{n}^{r,q} =\textstyle \yb_{n}^{r,q-1} - \gamma   \vb_n^{r, q-1},\label{eq: y update LSGT FL 2}\\
                & \vb_{n}^{ r, q } =  \textstyle  \vb_n^{r, q-1} +   \gb_n^{r, q} -   \gb_n^{r, q-1} ,\label{eq: v update LSGT FL 2}
            \end{align}
        \end{subequations}
        where the  stochastic gradient
            $\gb_n^{r, q}$
        is computed like \eqref{eq: h def} using $\yb_n^{r, q}$ and a mini-batch  $\mathcal{I}_n^{r,q}\subseteq \mathcal{D}_n$.
    	\ENDFOR
    \STATE  \vspace{-0.0cm}
    Set
     $\yb_n^{r+1} \!= \!\yb_n^{r, E}, \vb_n^{r+1} \!= \!\vb_n^{r, E},  \gb_n^{r+1} \!=\! \gb_n^{r,E}, 
     $ and send $(\yb_n^{r+1},\vb_n^{r+1})$ to neighbors.
	\ENDFOR
	\ENDFOR
	\end{algorithmic}
\end{algorithm}

\vspace{-0.2cm}
\subsection{Proposed LSGT Algorithm}
The proposed LSGT algorithm is presented in Algorithm \ref{alg: LSGT}. Comparing to the vanilla stochastic GT method in \eqref{eq: GT}, in the proposed LSGT algorithm, each agent $n$ executes $E$ consecutive steps of SGD within each communication round. In particular, in each step $q \in [E]$, agent $n$ performs gradient descent along the direction $\vb_n^{r, q-1}$ as in \eqref{eq: y update LSGT FL 2}. Then, it randomly chooses a mini-batch dataset $\mathcal{I}_n^{r,q} $ (with size $|\mathcal{I}|$) and computes the associated stochastic gradient $\gb_n^{r, q}$
using $\yb_n^{r, q}$. The local auxiliary variable $\vb_n^{r, q-1}$ is locally updated following \eqref{eq: v update LSGT FL 2}. Finally, after $E$ local updates, each agent $n$ sends $(\yb_n^{r+1},\vb_n^{r+1})$ to its neighbors.

\vspace{-0.2cm}
\begin{Remark}{  (Comparison with  \cite{nguyen2022performance})
It is noticed that the locally updated GT (LU-GT) algorithm recently proposed in \cite{nguyen2022performance} also considers multiple local updates for the GT method. However, there are several distinctions between the LU-GT algorithm and the proposed LSGT algorithm in Algorithm \ref{alg: LSGT}. Firstly, the LU-GT algorithm does not adopt SGD but uses full gradient at agents. Secondly,  the LU-GT algorithm requires an additional stepsize parameter for the local auxiliary variable update, unlike ours in \eqref{eq: v update LSGT FL 2} which does not involve any additional parameter. Besides, LU-GT requires the mixing matrix $\Wb$ to be symmetric whereas our LSGT algorithm does not.

Thirdly and mostly importantly, the convergence analysis in \cite{nguyen2022performance} does not fully characterize the impact of local updates on the algorithm convergence. In fact, their results (e.g., \cite[Remark 2]{nguyen2022performance}) somehow implies that local updates slow down the algorithm even under the assumption that the network graph is well connected.

}
\end{Remark}

In the next subsection, we present a novel convergence analysis for the LSGT algorithm, which shows the conditions under which the LSGT algorithm indeed can benefit from the local updates and enjoys a linear speedup with the network size $N$ and the number of local updates $E$.

\vspace{-0.2cm}
\subsection{Convergence Rate Analysis}
Let us denote $\Yb^r \triangleq [\yb_1^r, \ldots, \yb_N^r]^\top$, $\Vb^r \triangleq [\vb_1^r, \ldots, \vb_N^r]^\top$, and define the average of local variables as $\bar \yb^r = \frac{1}{N} \sum_{n=1}^N \yb_n^r$, and $\bar \vb^r = \frac{1}{N} \sum_{n=1}^N \vb_n^r$. Then, we write the consensus and tracking error as the following compact form
\begin{align}
\phib^{r} =\begin{bmatrix}
             \phi_{y}^{r}\\
             \phi_{v}^{r}
        \end{bmatrix} = \begin{bmatrix}
             \mathbb{E}\left[ \|\Yb^{r} - \mathbf{1}(\bar{\yb}^{r})^\top\|_F^2\right]\\
               \mathbb{E}\left[ \|\Vb^{r} - \mathbf{1}(\bar{\vb}^{r})^\top\|_F^2\right]
        \end{bmatrix}\in \mathbb{R}^2.
 \end{align}

Our main theoretical result for the LSGT method is given by the following theorem.

\vspace{-0.2cm}
\begin{Theorem} \label{thm: FL} Suppose that Assumption \ref{assum network} to \ref{assum: stochastic} hold.
For a sufficiently small $\gamma<1$ (satisfying \eqref{condi: zeta bsb}, \eqref{condi: zeta det I -G}, \eqref{condi: zeta phi sum}, \eqref{condi: zeta phi diff}, and the conditions in Lemma \ref{lem: y r q} to Lemma \ref{lem: v css}),
we have
\begin{align}\label{eq: thm FL}
  &  \frac{1}{T}\sum_{r=0}^{T-1} \mathbb{E}\bigg[\bigg\|\frac{1}{N}\sum_{n=1}^N \nabla f_n(\yb_n^r)\bigg\|^2 \bigg] \leq
  \underbrace{\frac{4(  F(\bar \yb^0) \!-  \!\underline{F})}{\gamma   ET}
\!+\!\frac{40L\gamma   \sigma^2}{N|\mathcal{I}|}}_{\text{ \rm terms same as centralized~SGD}}
   \!     \notag \\
& +\underbrace{\frac{ 16 (1+7\lambda_w^2)^2 E^2 L^2 \gamma^2  }{ (1-\lambda_w^2)^4  } \bigg(\frac{ 2577  N   \sigma^2}{|\mathcal{I}|} +\frac{111\phi_v^0}{T}\bigg)}_{\text {\rm terms due to decentralized optimization}}.
\end{align}
\end{Theorem}
{\bf Proof:} The first key of the proof is to build the dynamics of the consensus and tracking error matrix as follows
$$\phib^{r+1} \leq \Ab\phib^r +\Cb \eb^r,$$
for some coefficeint matrices $\Ab$ and $\Cb$ 
and the perturbation vector
\begin{align}
\eb^r = \begin{bmatrix}
   \mathbb{E} \bigg[ \bigg\|\frac{1}{N}\sum\limits_{n=1}^N \nabla f_n(\yb_n^r)\bigg  \|^2\bigg]
   \\ \frac{\sigma^2}{|\mathcal{I}|}\end{bmatrix}.\label{eq: e def}
\end{align}
The second key is to analyze how $\phib^{r}$  affects the descent of the objective value in \eqref{prob: HL 2}. By combining the two, we can obtain the bound in \eqref{eq: thm FL}.   Details are relegated to Section \ref{sec: proof of thm 1}. \hfill $\blacksquare$

Theorem \ref{thm: FL} implies the convergence of the proposed LSGT algorithm to the neighborhood of a stationary solution to problem \eqref{prob: HL 2} under an appropriate  constant stepsize $\gamma$. Specifically, the first two terms in the right hand side (RHS) of \eqref{eq: thm FL} are independent of the network topologies and have the same order as those of the centralized SGD method \cite{bottou2018optimization}. Meanwhile, the last two terms are due to the decentralized optimization and depend on the network connectivity $\lambda_w$ and the initial tracking error $\phi_v^0$.

More importantly, the following corollary provides conditions under which the proposed LSGT algorithm enjoys linear speedup with the local SGD number $E$ and network size $N$.

\begin{Corollary}\label{coro: coro 2 Q larger than 1}
    Let $\gamma = \textstyle \sqrt{\frac{N}{ ET }}$ and $E \leq \textstyle (\frac{T}{N^5})^{\frac{1}{3}}$ where $T$ is sufficiently large so that $\gamma$ satisfies the conditions in Theorem \ref{thm: FL}.
    Then, for the proposed LSGT algorithm, we have
     \begin{align}\label{eq: coro 2 Q larger than 1}
     \frac{1}{T}\sum_{r=0}^{T \!- \!1} \mathbb{E}\bigg[\bigg\|\frac{1}{N}\sum_{n=1}^N \nabla f_n(\yb_n^r)\bigg\|^2 \bigg] &\leq
  \frac{4( F(\bar\yb^{0}) \!-  \!\underline{F})}{ \sqrt{NET}}
\!+\!\frac{40 L  \sigma^2}{\sqrt{NET}|\mathcal{I}|}\notag\\
&~~~
   \!   +  \frac{16  (1 \!+ \!7\lambda_w^2)^2   L^2  }{ (1 \!- \!\lambda_w^2)^4 \sqrt{NET} } \bigg(\frac{ 2577     \sigma^2}{  |\mathcal{I}|}  \!+ \!\frac{111\phi_v^0}{NT }\bigg).
\end{align}
\end{Corollary}%

Corollary \ref{coro: coro 2 Q larger than 1} shows that for a large $T$, the stationary gap of LSGT decays sublinearly at the rate of $\mathcal{O}(  1/\sqrt{NET})$.  This is faster than the convergence rate $\mathcal{O}(1/\sqrt{NT})$ of the vanilla stochastic GT method  \cite[Corollary 1]{xin2021improved}, and thus
well demonstrates the benefits of employing local SGD with $E>1$ for reducing the communication overhead.



\vspace{-0.3cm}
\begin{Remark} \label{rmk: tau FL}
 {  (Impact of network connectivity and stochastic gradient error)}
From both \eqref{eq: thm FL} and \eqref{eq: coro 2 Q larger than 1}, one can see that a smaller $\lambda_{w}$ can reduce the 3rd and 4th terms of the RHS bound. Since a smaller $\lambda_{w}$ implies a higher network connectivity, it shows that the LSGT algorithm can converge and reach variable consensus faster if the network is more connected. On the other hand, one can also see that a larger mini-batch size $ |\mathcal{I}|$ can improve the convergence performance.
\end{Remark}

\vspace{-0.2cm}
\section{Proof of Theorem 1}\label{sec: proof of thm 1}

In this section, we present the proof of Theorem \ref{thm: FL}. Readers who are not interested in the proof may skip this section and jump to Section \ref{sec: hybrid} for an extension of the LSGT algorithm.


\vspace{-0.2cm}
\subsection{Key Lemmas}
The proof of Theorem \ref{thm: FL} relies on four key lemmas given below. Their proofs are relegated to Appendix B to D and Section I in the supplementary material. For ease of presentation, a preliminary is provided in Appendix A.

The first lemma bounds the distance between local variables $(\yb_n^{r, q}, \vb_n^{r, q})$ and their corresponding average $(\bar \yb^r, \bar \vb^r)$.

\begin{Lemma} \label{lem: y r q}
\emph{(Local Update Gap)} Suppose that Assumptions \ref{assum network} to \ref{assum: stochastic} hold. For a sufficiently small $\gamma\!   \leq \!
   \textstyle \frac{  \lambda_{w}^2} {32 EN L  }$, we have
\begin{subequations}
\begin{align}
  \sum_{n=1}^N \sum_{q=1}^{E-1}  \mathbb{E}  [  \| \yb_n^{r,q} -\bar \yb^r  \|^2  ]
 & \leq  16(E-1)\lambda_{w}^2 \phi_{y}^r +8(E-1)\lambda_{w}^2[1+ E(E-1)  ] \gamma^2\phi_{v}^r\notag\\
    &~~~\!
    +  \!4E^2(E\!-\!1)N\gamma^2 \mathbb{E}  [  \|  \bar \vb^r \|^2  ] \!+ \! 64E(E\!-\!1)^2 N  \gamma^2 \frac{\sigma^2}{|\mathcal{I}|}, \!\label{eq: y r q}
    \end{align}
    and
    \begin{align}
 \sum_{n=1}^N \sum_{q=1}^{E-1}  \mathbb{E}   [  \| \vb_n^{r,q} -\bar \vb^r   \|^2  ]
 & \leq  8(E-1) L^2 (1+ 16 \lambda_{w}^2 )   \phi_{y}^r  +4(E-1)\lambda_{w}^2   \phi_{v}^r\notag\\
    &~~~ \!+ \!32   E^2(E\!-\!1) N  L^2\gamma^2  \mathbb{E}   [  \|  \bar \vb^r \|^2  ]
    \!+\!32(E\!-\!1)N   \frac{\sigma^2}{|\mathcal{I}|} .\label{eq: v r q}
\end{align}
\end{subequations}
\end{Lemma}

Based on the above Lemma \ref{lem: y r q}, in the next two lemmas we determine the contraction properties of the consensus error $\phi_{y}^{r}$ and tracking error $\phi_{v}^r$, respectively.
\begin{Lemma} \label{lem: y css}
\emph{(Consensus Error Bound)} Suppose that Assumptions \ref{assum network} to \ref{assum: stochastic} hold. For a sufficiently small $\gamma$  satisfying
\begin{align}
    \gamma   \leq \frac{(1-\lambda_{w}^2) ^2}{32 L^2(1+7\lambda_{w}^2)[8(1+16\lambda_{w}^2)(E-1)^2  ]},
\end{align}
 we have
\begin{align}
&  \phi_{y}^{r+1} = \mathbb{E}\left[\left\|\Yb^{r+1} - \mathbf{1}(\bar\yb^r)^\top\right\|^2\right]\notag\\
& \leq  \frac{(1+3\lambda_{w}^2) }{4}  \phi_{y}^r    +  \frac{2(1+7\lambda_{w}^2)  }{1-\lambda_{w}^2} \left[1
    + 8(E-1)^2   \lambda_{w}^2  \right] \gamma^2\phi_{v}^r\notag\\
    &~~~\! +  \!\frac{256(1\!+\!7\lambda_{w}^2) E^2(E\!-\!1)^2  N  L^2 \gamma^4 }{1\!-\!\lambda_{w}^2} \!\mathbb{E} \bigg[ \bigg\|\frac{1}{N} \!\sum_{n=1}^N\!\nabla f_n(\yb_n^r)\bigg  \|^2\bigg]
     +   \frac{256(1+7\lambda_{w}^2) (E-1)^2 N \gamma^2\sigma^2 }{(1-\lambda_{w}^2)|\mathcal{I}|} . \label{eq: y css}
\end{align}
\end{Lemma}

\begin{Lemma} \label{lem: v css}
\emph{(Tracking Error)} Suppose that Assumptions \ref{assum network} to \ref{assum: stochastic} hold. For a sufficiently small $\gamma\leq  \textstyle \frac{(1-\lambda_{w}^2)}{3600EN L }$,
we have
\begin{align}
&  \phi_{v}^{r+1} = \mathbb{E}\left[\left\|\Vb^{r+1} - \mathbf{1}(\bar\vb^r)^\top\right\|^2\right]\notag\\
& \leq \frac{30(1\!+\!7\lambda_{w}^2)(1\!+\!\lambda_{w}^2)L^2}{1\!-\!\lambda_{w}^2} \phi_{y}^r
 \!+ \!\frac{1\!+\!3\lambda_{w}^2}{4}  \phi_{v}^r
 \!+   \!\frac{12(1\!+\!7\lambda_{w}^2) N  \sigma ^2 }{  (1\!-\!\lambda_{w}^2)|\mathcal{I}| }
    \notag\\
    &~~~ +  \frac{60(1+7\lambda_{w}^2) }{1-\lambda_{w}^2}  E^2N L^2\gamma^2 \mathbb{E} \bigg[ \bigg\|\frac{1}{N} \sum_{n=1}^N\nabla f_n(\yb_n^r)\bigg  \|^2\bigg]
      . \label{eq: v css}
\end{align}
\end{Lemma}

Next, we characterize a key descent property of the global objective function.
\begin{Lemma} \label{lem: descent lemma for F HL}
    \emph{(Descent Lemma)} Let Assumptions \ref{assum network} to \ref{assum: stochastic} hold. Based on Lemma \ref{lem: y r q} to \ref{lem: v css}, we have $\forall r\geq0$,
\begin{align}\label{eq: f FL}
   &   F(\bar \yb^{r+1})  \leq   F(\bar \yb^{r}) \!  -\!\frac{ \gamma  E}{2}
                \mathbb{E}[  \| \nabla F(\bar \yb^{r })\|^2]
        \!-\!\frac{ \gamma }{2}\bigg[1\!-\! 2E^2  L\gamma
       \! - \!8 E^2(E\!-\!1) L^2 \gamma^2
        \bigg] \mathbb{E} \bigg[ \bigg\| \frac{1}{N}\sum_{n=1}^N \! \nabla f_n(\yb_n^{r }) \bigg\|^2\bigg] \notag\\
    &~~~~~~~~~~+   \frac{L^2 \gamma}{N}  [1+  16(E-1)\lambda_{w}^2] \phi_{y}^r
    +10 E L   \gamma^2    \frac{ \sigma^2}{N |\mathcal{I}|} + 8(E-1)\lambda_{w}^2[1+ E(E-1)  ]  \frac{ L^2 \gamma^3}{N}   \phi_{v}^r,
\end{align}
where $ \nabla F(\bar \yb^{r }) = \textstyle \mathbb{E}[\frac{1}{N}\sum_{n=1}^N \nabla f_n(\bar \yb^r) ]$.
\end{Lemma}

\subsection{Proof of Theorem \ref{thm: FL}}
By  Lemma \ref{lem: y css} and Lemma \ref{lem: v css}, we can establish a dynamics system of $\phib^{r}$ as follows
\begin{align}\label{eq: phi FL}
    \phib^{r+1} \leq \Ab\phib^{r} +\Cb \eb^r,
\end{align}
where the inequality is element-wise, and
the matrices $\Ab$ and $\Cb$ are given by
\begin{align}
    &\! \Ab\! =\!  \begin{bmatrix}
     \frac{(1+3\lambda_{w}^2) }{4}  \!&\!    \frac{2(1+7\lambda_{w}^2)\left[1
    + 8(E-1)^2   \lambda_{w}^2   \right]  \gamma^2  }{1-\lambda_{w}^2}\!\\
      \! \frac{30(1+7\lambda_{w}^2)(1+\lambda_{w}^2 ) L^2 }{1-\lambda_{w}^2}
      \!& \!\frac{1+3\lambda_{w}^2}{4}
     \end{bmatrix}\label{eq: G def},\\
    &\! \Cb\! = \! \begin{bmatrix}
    \frac{256(1 \! + \! 7\lambda_{w}^2) E^2(E \! - \! 1)^2  N  L^2 \gamma^4  }{1 \! - \! \lambda_{w}^2}
            \!&\!
            \frac{256(1 \! + \! 7\lambda_{w}^2) (E \! - \! 1)^2 N  \gamma^2  }{1 \! - \! \lambda_{w}^2}   \\
            \frac{60(1 \! + \! 7\lambda_{w}^2) E^2 N  L^2 \gamma^2  }{1 \! - \! \lambda_{w}^2}
           \! & \!  \frac{12(1 \! + \! 7\lambda_{w}^2) N  }{1 \! - \! \lambda_{w}^2}
     \end{bmatrix}. \!\!\label{eq: H def}
\end{align}

One can verify that there exists a positive vector $\bsb=[s_1, s_2]^\top \in \mathbb{R}^{2}$ satisfying
\begin{align}
    s_1& <\frac{ (1-\lambda_{w}^2)^2 }{40 (1+7\lambda_{w}^2)(1+\lambda_{w}^2) L^2}s_2,
\end{align}
such that  $ \Ab\bsb <\bsb$ as long as
\begin{align}
   \gamma & \leq       \frac{3(1-\lambda_{w}^2) s_1}{16(1+7\lambda_{w}^2) [1 + 8(E-1)   \lambda_{w}  ]  s_2   }  . \label{condi: zeta bsb}
\end{align}%
Thus, under \eqref{condi: zeta bsb}, the spectral radius of $\Ab$ satisfies $\rho(\Ab)<1$ according to \cite[Corollary 8.1.29]{horn2012matrix}.
Then, by \cite[Corollary 5.6.16]{horn2012matrix}, we can have
 \begin{align}\label{eq: I-G FL}
 & \det(\mathbf{I}  - \Ab) >0,~~
  \sum_{r=0}^{\infty} \Ab^r = (\mathbf{I} - \Ab)^{-1}.
 \end{align}
As a result \eqref{eq: phi FL} can be bounded as
\begin{align}\label{eqn: sum phi}
\phib^{r}  \leq \Ab  \phib^{r-1}  + \Cb\eb^{r-1}
\leq \Ab^r \phib_{0}  + \sum_{t=0}^{r-1}\Ab^t\Cb\eb^{r-1-t}.
\end{align}
Summing \eqref{eqn: sum phi} for $r$ from $0$ to $T$, we further have
\begin{align}\label{eq: sum phi FL 1}
\sum_{r=0}^{T}\phib^{r}
& \leq  \sum_{r=0}^{T}\Ab^r \phib^{0}  +
        \sum_{r=0}^{T}\sum_{t=1}^{r-1}\Ab^t\Cb\eb^{r-1-t}\notag\\
& \leq \bigg(\sum_{r=0}^{\infty}\Ab^r \bigg)\phib^{0}  +
       \bigg( \sum_{r=0}^{\infty}\Ab^r \bigg) \sum_{r=0}^{T-1}\Cb\eb^{r} \notag\\
& =  (\mathbf{I} - \Ab)^{-1}\phib^{0}
+ (\mathbf{I} - \Ab)^{-1}\sum_{r=0}^{T-1}\Cb\eb^{r},
\end{align}
where the last equality is owing to \eqref{eq: I-G FL}.

Since $\phib^r$ is always positive, we have $\sum_{r=0}^{T-1}\phib^{r} \leq \sum_{r=0}^{T}\phib^{r}$. Thus, \eqref{eq: sum phi FL 1} implies
\begin{align}\label{eq: sum phi FL}
\sum_{r=0}^{T-1}\phib^{r}
& \leq    (\mathbf{I} - \Ab)^{-1}\phib^{0}
+ (\mathbf{I} - \Ab)^{-1}\sum_{r=0}^{T-1}\Cb\eb^{r}.
\end{align}%

It is shown in Section III of the Supplementary materials that $(\Ib- \Ab)^{-1}$ can be bounded as follows.
\vspace{-0.2cm}
\begin{Lemma}\label{lem: det I -G}
 Suppose that $\gamma$ satisfies
\begin{align}\label{condi: zeta det I -G}
\gamma \leq \textstyle \frac{3(1-\lambda_{w}^2)^2}{320(1+7\lambda_{w}^2) (1+\lambda_{w}^2)[1
    + 8(E-1 )  \lambda_{w} ]N  L }.
    \end{align}
 Then, $(\Ib- \Ab)^{-1}$ has an element-wise upper bound given by
\begin{align}\label{eq: upper bound det I-G}
    &(\mathbf{I}-\Ab)^{-1}\notag\\
    & \!\leq\!
    \frac{8(1\!+\!7\lambda_{w}^2)  }{(1\!-\!\lambda_{w}^2)^3}
    \!\begin{bmatrix}
        1 \!& \! \!2    [1
    \!+\! 8(E\!-\!1)^2   \lambda_{w}^2   ]   \gamma^2\!\\
        30 (1\!+\!\lambda_{w}^2)   L^2
         \!&\! 1
    \end{bmatrix}\!.
\end{align}
 \end{Lemma}
%

Inserting \eqref{eq: upper bound det I-G} into \eqref{eq: sum phi FL} and by the fact of $\phi_{y}^0 = 0$, we obtain
\begin{align}
      \sum_{r=0}^{T-1} \phi_{y}^r
  &\leq \frac{64(1 \!+ \!7\lambda_w^2)^2}{(1 \!- \!\lambda_w^2)^4}\bigg[
    32(E \!- \!1)^2 \!+  \!15\lambda_w^2[1 \!+ \!8(E \!- \!1)^2 ]\bigg]
 E^2N L^2\gamma^4
    \sum_{r=0}^{T-1} \mathbb{E}\! \bigg[ \bigg\|  \frac{1}{N}\sum_{n=1}^N  \nabla f_n(\yb_n^{r }) \bigg\|^2\bigg] \notag\\
    & ~~~
        +\! \frac{64(1\!+\!7\lambda_w^2)^2}{(1\!-\!\lambda_w^2)^4}\bigg[
        32(E\!-\!1)^2 \!  + \!3\lambda_w^2[1\!+\!8(E\!-\!1)^2 ]\bigg]
      \frac{  N\gamma^2  T  \sigma^2}{|\mathcal{I}|} \notag\\
    &~~~
        + \frac{16(1+7\lambda_w^2)^2}{(1-\lambda_w^2)^4}\lambda_w^2
        [1+8(E-1)^2]\gamma^2 \phi_{v}^0 , \label{eq: sum phi y}
        \end{align}
        and
        \begin{align}
 &\sum_{r=0}^{T-1} \phi_{v}^r
    \leq\frac{8(1+7\lambda_w^2)^2}{(1-\lambda_w^2)^4}
          \phi_{v}^0
        \!  +
     \! \frac{32(1+7\lambda_w^2)^2}{(1-\lambda_w^2)^4}
[1920(E-1)^2   (1+\lambda_w^2)L^2\gamma^2 + 3]\frac{N   T \sigma^2}{|\mathcal{I}|}
      \! \notag\\
    & ~~~
+\!\frac{ 480(1\!+\!7\lambda_w^2)^2 E^2N L^2\gamma^2}{(1\!-\!\lambda_w^2)^4}   [128(E\!-\!1)^2 (1\!+\!\lambda^2)L^2\gamma^2 \!+ \!1]
 \sum_{r=0}^{T\!-\!1}\! \mathbb{E} \bigg[ \bigg\| \frac{1}{N}\sum_{n=1}^N \! \!\nabla f_n(\yb_n^{r }) \bigg\|^2\bigg]\!
   \label{eq: sum phi v proof pro}\\
  &\!\leq\frac{8(1+7\lambda_w^2)^2}{(1-\lambda_w^2)^4}
          \phi_{v}^0
     \! +\frac{128(1+7\lambda_w^2)^2 N   T \sigma ^2}{(1-\lambda_w^2)^4 |\mathcal{I}|}
        \!
+\!\frac{ 512(1\!+\!7\lambda_w^2)^2 E^2N L^2\gamma^2}{(1\!-\!\lambda_w^2)^4}
    \sum_{r=0}^{T\!-\!1}\!  \mathbb{E}\bigg[ \bigg\| \frac{1}{N}\!\sum_{n=1}^N \! \!\nabla f_n(\yb_n^{r }) \bigg\|^2\bigg],
        \label{eq: sum phi v}
 \end{align}
where, in order to obtain \eqref{eq: sum phi v}, we have assumed
\begin{align}\label{condi: zeta phi sum}
      \gamma &\leq
        \frac{1}{62E(1+\lambda_{w}^2)N L },
\end{align} and applied it to
the $2$nd and $3$rd terms in the RHS of \eqref{eq: sum phi v proof pro}.

On the other hand, by taking the telescope sum of \eqref{lem: descent lemma for F HL} in Lemma \ref{lem: descent lemma for F HL} from $r=0$ to $r=T-1$, we have
\begin{align} \label{eq: upper sum bar nabla}
  &\frac{ \gamma  E}{2}\sum_{r=0}^{T-1} \mathbb{E} [ \|\nabla F(\bar \yb^r)  \|^2  ] \notag \\
   & \leq  F(\bar\yb^{0}) \!-\! F(\bar\yb^{T})\!+ \!8(E\!-\!1)\lambda_{w}^2[1\!+\! E(E\!-\!1)  ] \frac{  L^2 \gamma^3 }{N}  \sum_{r=0}^{T\!-\!1} \! \phi_{v}^r
        \notag\\
     & ~~~ -\!\frac{ \gamma }{2}[1\!-\! 2E^2  L\gamma
            \!- \!8 E^2(E\!-\!1) L^2 \gamma^2
        ] \sum_{r=0}^{T-1}\!\mathbb{E}\! \bigg[ \bigg\| \frac{1}{N}\sum_{n=1}^N \! \nabla f_n(\yb_n^{r }) \bigg\|^2\bigg] \notag\\
    &~~~ +   \frac{L^2 \gamma  [1+  16(E-1)\lambda_{w}^2] }{N}\sum_{r=0}^{T-1} \phi_{y}^r
    +\frac{10E  L   \gamma^2 T  \sigma ^2 }{N|\mathcal{I}|} .
\end{align}
Besides, by Assumption \ref{assum: Lipschitz}, one have
\begin{align}\label{eq: sum real gr 1}
   & \sum_{r=0}^{T-1} \mathbb{E}\bigg[\bigg\| \frac{1}{N}\sum_{n=1}^N \nabla f_n(\yb_n^r)\bigg\|^2 \bigg] \notag\\
   & \leq  2 \sum_{r=0}^{T\!-\!1} \mathbb{E}\bigg[\bigg\| \frac{1}{N}\sum_{n=1}^N [\nabla f_n(\yb_n^r) \!- \!\nabla f_n(\bar \yb^r)]\bigg\|^2 \!+\!2  \| \nabla F(\bar \yb^r)  \|^2   \bigg] \notag\\
   &\leq  \frac{2 L^2 }{N}\sum_{r=0}^{T-1}\phi_y^r
   +  2 \sum_{r=0}^{T-1} \mathbb{E} [ \| \nabla F(\bar \yb^r)  \|^2  ].
\end{align}
By multiplying $\textstyle \frac{\gamma   E}{4}$ on both sides of \eqref{eq: sum real gr 1} followed by inserting it into \eqref{eq: upper sum bar nabla}, and by Assumption \ref{assum: lower bound}, we have
\begin{align} \label{eq: sum real gr 21}
  &  \frac{\gamma  E}{4} \sum_{r=0}^{T-1} \mathbb{E}\bigg[\bigg\|\frac{1}{N}\sum_{n=1}^N \nabla f_n(\yb_n^r)\bigg\|^2 \bigg]     \notag \\
   & \leq  F(\bar\yb^{0}) - \underline{F}  +    \frac{10E  L   \gamma^2 T \sigma^2}{N |\mathcal{I}| }
        \!-\!\frac{ \gamma }{2}[1\!-\! 2E^2  L\gamma
            \!- \!8 E^2(E\!-\!1) L^2 \gamma^2
        ] \mathbb{E}\! \bigg[ \bigg\| \frac{1}{N}\sum_{n=1}^N  \nabla f_n(\yb_n^{r }) \bigg\|^2\bigg] \notag\\
    & +\bigg\{2[1+  16(E-1)\lambda_{w}^2]+E\bigg\} \frac{ L^2 \gamma   }{2N}   \sum_{r=0}^{T-1} \phi_{y}^r
    + 8(Q-1)\lambda_{w}^2 [1+ E(E-1) ]  \frac{  L^2 \gamma^3  }{N} \sum_{r=0}^{T-1}  \phi_{v}^r .
\end{align}

Then, by inserting \eqref{eq: sum phi y} and \eqref{eq: sum phi v} into the above \eqref{eq: sum real gr 21}, we obtain
\begin{align}
  &  \frac{\gamma  E}{4} \sum_{r=0}^{T-1} \mathbb{E}\bigg[\bigg\|\frac{1}{N}\sum_{n=1}^N \nabla f_n(\yb_n^r)\bigg\|^2 \bigg]     \notag \\
& \leq  (F(\bar\yb^{0}) - \underline{F}) +    \frac{10E  L   \gamma^2 T \sigma^2}{N |\mathcal{I}| }
        \notag\\
     & + \bigg\{
          \frac{4(1\!+\!7\lambda_w^2)^2}{(1\!-\!\lambda_w^2)^4}\bigg[
        32(E\!-\!1)^2 \!
        + \!3\lambda_w^2[1\!+\!8(E\!-\!1)^2 ]\bigg]
         \bigg[ 2[1+  16(E-1)\lambda_{w}^2]+E \bigg] N  L \gamma
     \notag\\
     &~~~ 
     \!   +  \! \frac{128(1\!+\!7\lambda_w^2)^2}{(1\!-\!\lambda_w^2)^4}(E\!-\!1)\lambda_{w}^2 [1\!+\! E(E\!-\!1) ]
          N   L \gamma
     \bigg\}  \frac{8  L   \gamma^2    T {\sigma}^2 }{N |\mathcal{I}|} \notag\\
& + \frac{8(1+7\lambda_w^2)^2}{(1-\lambda_w^2)^4}\bigg\{
         8(E-1)\lambda_{w}^2[1+ E(E-1)  ]
    \!+\! \lambda_w^2[1\!+\!8(E\!-\!1)^2 ] \{2[1\!+\!  16(E\!-\!1)\lambda_{w}^2]+E\}
    \bigg\}\!  \frac{L^2  \gamma^3}{N}  \phi_{v}^0 \notag\\
& \!-\!\frac{ \gamma }{2}
    \bigg\{1\!-\! 2E^2  L\gamma
            \!- \!8 E^2(E\!-\!1)N^2L^2 \gamma^2
         -\frac{32(1 \!+ \!7\lambda_w^2)^2}{(1 \!- \!\lambda_w^2)^4}
        \{2[1+  16(E-1)\lambda_{w}^2]+E\}\notag\\
     &~~~~~~
        \{32(E \!- \!1)^2 \!+  \!15\lambda_w^2[1 \!+ \!8(E\!- \!1)^2 ]\}
        E^2 L^4\gamma^4
    \notag\\
     &~~~ -\frac{4096(1\!+\!7\lambda_w^2)^2}{(1\!-\!\lambda_w^2)^4}
         (E-1)\lambda_{w}^2[1+ E(E-1)  ]E^2  L^4\gamma^4
    \bigg\}
      \sum_{r=0}^{T-1}
    \mathbb{E}\! \bigg[ \bigg\|\frac{1}{N} \sum_{n=1}^N  \nabla f_n(\yb_n^{r }) \bigg\|^2\bigg] \label{eq: sum real gr 2 pro 1}\\
& \leq (F(\bar\yb^{0}) \!-  \!\underline{F})+   \!    \frac{10E  L   \gamma^2 T \sigma^2}{N |\mathcal{I}| }+
\frac{ 444(1+7\lambda_w^2)^2 E^3 L^2 \gamma^3 \phi_v^0 }{N (1 \!- \!\lambda_w^2)^4}
 \!
    +\frac{ 10308(1 \!+7\lambda_w^2)^2 E^3 L^2 \gamma^3  T \sigma^2}{ (1-\lambda_w^2)^4 |\mathcal{I}|  } ,\label{eq: sum real gr 2}
\end{align}
To obtain \eqref{eq: sum real gr 2}, we have assumed
\begin{align}\label{condi: zeta phi diff}
\gamma \leq \frac{1-\lambda_w^2}{320E^2 NL},
\end{align} so that the coefficient of the last term in the RHS of \eqref{eq: sum real gr 2 pro 1} is negative and the term can be removed. In addition, we have used the properties of $\lambda_w^2< 1$ and $E\geq 1$ to obtain bounds for coefficients of the 3rd and 4th terms in the RHS of \eqref{eq: sum real gr 2}.
Finally, after dividing $\textstyle \frac{T\gamma  E}{4}$ on both sides of \eqref{eq: sum real gr 2}, we obtain the results in Theorem \ref{thm: FL}.

\vspace{-0.1cm}
\section{Extension for Learning over Hybrid Data}\label{sec: hybrid}
In this section, we discuss an extension of the LSGT algorithm to solve the following optimization problem
\begin{align}
\min_{  (\xb,\thetab) \in \mathbb{R}^{J+K} }~~  & \frac{1}{NS}\sum_{i = 1}^S  f \left( \Bb_{ i}\xb, \thetab \right)\label{pro: HFL rewrite} 
\end{align}
where  \vspace{-0.2cm}\begin{align}\label{eqn: Bni}
\Bb_i\!=\!\sum_{n=1}^N \Bb_{n,i}  \in \mathbb{R}^{M\times J}, \forall i \in[S].
\end{align}
An instance of problem \eqref{pro: HFL rewrite} appears in a distributed learning problem over the hybrid data where the agents can access only a subset of data samples and knows only partial data features while they collaborate to train a deep neural network (DNN) model. Detailed descriptions about learning over hybrid data can be found in \cite{FedHD2022}.
Specifically, in
problem \eqref{pro: HFL rewrite},
$S$ denotes the number of data samples in the global dataset $\mathcal{D}=\cup_{n=1}^N\mathcal{D}_n$, $\xb$ and $\thetab$ are the coefficient parameters of the first layer and the remaining layers of the DNN model, respectively, and $\Bb_{i}$ is a matrix related to the $i$th data sample. By \eqref{eqn: Bni}, $\Bb_{n, i}$ is the matrix about the feature information of the $i$th data sample that agent $i$ knows exclusively. If $i \notin D_n$, i.e., agent $i$ does not own sample $i$, then  $\Bb_{n, i} =\zerob$.

Thus, comparing to problem \eqref{prob: HL 2}, the challenge of problem \eqref{pro: HFL rewrite} lies in that each agent $n$ does not know $\Bb_i$ but can access $\Bb_{n,i}$ only.
As a result, the LSGT algorithm in Algorithm \ref{alg: LSGT} cannot be directly applied to handle problem \eqref{pro: HFL rewrite}. 
To see this, observe that the local updates for problem \eqref{pro: HFL rewrite} would involve computing
\vspace{-0.2cm}
\begin{subequations}
\label{eq: consensus gradient descent HBL}
\begin{align}
 \thetab_{n}^{ r, q }& \!=  \! \thetab_n^{r,q\!-\!1}\!-\!
                 \alpha \frac{1}{NS} \sum_{i=1}^S  \nabla_{\thetab}f  \big( \underbrace{\Bb_i\xb_n^{r,q\!-\!1}}_{(a)}, \thetab_n^{r,q\!-\!1}   \big )
                 , \label{eq: update theta initial 1}\\
 \xb_{n}^{ r, q}&\!= \!   \xb_n^{r,q\!-\!1}
           \! -\! \beta \underbrace{ \frac{1}{NS} \sum_{i=1}^S\!\Bb_{i}^\top  \nabla_{\zb} f\big(\Bb_{i}\xb_n^{r,q-1}, \thetab_n^{r,q\!-\!1}   \big )}_{(b)}  ,\label{eq: update x initial 1}
\end{align}
\end{subequations} where $ \alpha,\beta >0$ are step size parameters, and $\nabla_{\zb} f(\zb,\theta)$ is the gradient with respect to $\zb$. One can see that the updates in
\eqref{eq: consensus gradient descent HBL} are unfortunately not realizable since each agent lacks the full information of $\Bb_i$ and cannot compute the terms (a) and (b) locally.

In view of this, we introduce a new auxiliary variable
$\zb_{n,i}^{r,q}$ to estimate the summation term (a), i.e., $\Bb_i\xb_n^{r,q}$, and another auxiliary variable $\ub_n^{r,q}$ to estimate term (b).
Given $\{\zb_{n,i}^{r,q}\}_{i\in D}$, agent $n$ can compute the mini-batch gradient locally by
\begin{subequations}\label{eq: stc gr HBL}
\begin{align}
    &
     \gb_{\theta,n}^{r,q} \triangleq     \frac{1}{N|\mathcal{I}|}\sum_{\xi   \in \mathcal{I}_n^{r,q}} \nabla_{\thetab} f( \zb_{n, \xi  }^{r,q}, \thetab_n^{r,q} ), ~\\
    &
    \gb_{x,n}^{r,q} \triangleq     \frac{1}{ |\mathcal{I}|}  \sum_{\xi  \in \mathcal{I}_n^{r, q}} \Bb_{n, \xi }^\top \nabla_{\zb}f( \zb_{n, \xi}^{r, q}, \thetab_n^{r, q}),
\end{align}
\end{subequations}
while, by consensus averaging, tracking term (a) and term (b) via exchanging the tracking variable $\zb_{n}^{r,q}=[(\zb_{n,1}^{r,q})^\top, \ldots, (\zb_{n,S}^{r,q})^\top]^\top$
and $\gb_{x,n}^{r,q}$
with neighbors, respectively.
Thus, \eqref{eq: consensus gradient descent HBL}
can be replaced by
  \begin{align}
                            & \thetab_n^{r, q} = \thetab_n^{r, q-1} - \alpha
                                 \gb_{\theta,n}^{r,q-1}, \\
                            & \xb_n^{r, q} = \xb_n^{r, q-1} -\beta
                                \ub_n^{r, q-1}.
\end{align}
Like the LSGT algorithm, the agents can still perform $E$ local SGD steps in each communication round for reducing the communication cost. We summarize the proposed algorithm for solving problem \eqref{pro: HFL rewrite} in
Algorithm \ref{alg: MUST}, which we refer to as the \emph{Multiple-locally-Updated variable Sum Tracking (MUST)} algorithm.

\begin{algorithm}[!t]
	\caption{Proposed MUST method for solving \eqref{pro: HFL rewrite}}
	\begin{algorithmic}[1]\label{alg: MUST}
		\STATE \textbf{Initialize} Let $ \xb_1^0 =\ldots = \xb_N^0$,
$ \thetab_1^0 =\ldots = \thetab_N^0$, $ \zb_{n,i}^{0} = N\Bb_{n,i}\xb_n^{0}$, and $  \ub_{n}^{0} = \gb_{x,n}^{0},\forall n\in [N]$.
        \FOR {communication round $r=0$ {\bf{to}} $T$}
		\FOR{agent $n=1$ {\bf{to}} $N$ in parallel}
        \STATE Receive information from neighbours and set
         \vspace{-0.0cm}
                \begin{align}
                    \begin{bmatrix}
                    \thetab_n^{r,0}\\
                    \xb_n^{r,0}\\
                    \zb_n^{r,0}\\
                    \ub_n^{r,0}
                    \end{bmatrix}
                    =  \sum_{m=1}^N W_{n,m}
                    \begin{bmatrix}
                    \thetab_m^{r,0}\\
                    \xb_m^{r,0}\\
                    \zb_m^{r,0}\\
                    \ub_m^{r,0}
                    \end{bmatrix},
                    \begin{bmatrix}
                    \gb_{\theta,n}^{r,0}\\
                     \gb_{x,n}^{r,0}
                     \end{bmatrix}
                     = \begin{bmatrix}
                     \gb_{\theta,n}^{r}\\
                     \gb_{x,n}^{r}
                     \end{bmatrix}
                \end{align}
                 \vspace{-0.6cm}
                \FOR{local update $q = 1, \ldots, E$}
                \STATE
          \vspace{-0.2cm} \begin{subequations}\label{eq: Alg2 q2}
                           \begin{align}
                            & \thetab_n^{r, q} = \thetab_n^{r, q-1} - \alpha
                                 \gb_{\theta,n}^{r,q-1}, \\
                            & \xb_n^{r, q} = \xb_n^{r, q-1} -\beta
                                \ub_n^{r, q-1},\\
                            & \zb_{n}^{r, q}\!= \!\zb_{n}^{r, q-1} \!+ \! N\Bb_{n}
                                (\xb_n^{r, q}\! - \!\xb_n^{r, q-1}),\label{eq: Alg2 q2 z}\\
                            & \ub_n^{r, q} = \ub_n^{r, q-1}
                                +    (\gb_{x,n}^{r,q} -\gb_{x,n}^{r,q-1}),\label{eq: Alg2 q2 u}
                         \end{align}
                         \end{subequations}
                    where the stochastic gradients $\gb_{\theta,n}^{r,q}$ and $\gb_{x,n}^{r,q}$ are computed based on \eqref{eq: stc gr HBL} using $\thetab_n^{r,q}$ and $\xb_n^{r,q}$ and a mini-batch $\mathcal{I}_n^{r, q}$, and $\Bb_n=[\Bb_{n, 1}^\top,\ldots, \Bb_{n, S}^\top]^\top$.
                \ENDFOR
                \STATE
                    Set $\thetab_n^{r+1} = \thetab_n^{r,E}$,
                       $ \xb_n^{r+1} = \xb_n^{r,E}$,
                        $\zb_n^{r+1} = \zb_n^{r,E}$,
                        $\ub_n^{r+1} = \ub_n^{r,E}$, $\gb_{\theta,n}^{r+1}= \gb_{\theta,n}^{r,E} $
                    $\gb_{x,n}^{r+1}= \gb_{x,n}^{r,E} $, and send $(\thetab_n^{r+1}, \xb_n^{r+1},\zb_n^{r+1},\ub_n^{r+1})$
                    to neighbours.
          \ENDFOR
          \ENDFOR
	\end{algorithmic}
\end{algorithm}

\vspace{-0.2cm}
\subsection{Convergence Rate Analysis}

In this subsection,  we
build the convergence conditions for the MUST algorithm.
Like Assumption \ref{assum: stochastic}, we make the following assumption for the mini-batch gradients in \eqref{eq: stc gr HBL}.
\begin{Assumption} 
\label{assum: stochastic HBL} For each agent $n\in[N]$, it holds that

\noindent$\bullet$ Unbiased gradient: $\mathbb{E} [  \gb_{\theta,n} ] = \frac{1}{NS} \sum_{i=1}^S \nabla_{\thetab} f(\zb_{n,i} , \thetab_n ) $,

$\mathbb{E} [  \gb_{x,n} ] =  \frac{1}{ S} \sum_{i=1}^S \Bb_{n,i}^\top\nabla_{\zb} f(\zb_{n,i} , \thetab_n ) $;

\noindent$\bullet$ Uniform bounded variance:

$\mathbb{E}  \textstyle [  \|   \gb_{\theta,n}   -   \frac{1}{NS} \sum_{i=1}^S \nabla_{\thetab} f(\zb_{n,i} , \thetab_n ) \|^2  ] \leq   \frac{\sigma^2}{|\mathcal{I}|},$

$\mathbb{E}   \textstyle [  \|   \gb_{x,n}   -  \frac{1}{ S} \sum_{i=1}^S \Bb_{n,i}^\top\nabla_{\zb} f(\zb_{n,i} , \thetab_n )  \|^2   ] \leq   \frac{\sigma^2}{|\mathcal{I}|}.$
\end{Assumption}%
We also assume that the cost function $f$ in \eqref{pro: HFL rewrite} is lower bounded and smooth, like Assumption \ref{assum: lower bound} and \ref{assum: Lipschitz}.
The main theoretical result for the MUST algorithm is given by the following theorem.
\begin{Theorem}\label{thm: HBL}
Let Assumption \ref{assum network} to  Assumption \ref{assum: Lipschitz}, and  Assumption \ref{assum: stochastic HBL} hold. Define the average $\bar \xb^r = \textstyle \frac{1}{N} \sum_{n=1}^N \xb_n^r$, and denote $\bar \thetab^r , \bar \zb^r,\bar \ub^r$ in the same fashion.
Choosing $\alpha = \beta = \textstyle \sqrt{\frac{N}{ET}}$, for $E \leq \textstyle \frac{T^{\frac{1}{3}}}{N^{3}}$, $T\geq N^{\frac{5}{2}}E^3$, we have
\begin{align}\label{eq: gr sum theta x HBL coro}
&\frac{1}{T}\sum_{r=0}^{T-1} \frac{1}{N} \sum_{n=1}^N \bigg\{ \mathbb{E}\bigg[\bigg\| \frac{1}{NS}\sum_{i=1}^S\nabla_{\thetab} f\bigg(\sum_{t=1}^N \Bb_{t,i}\xb_n^r, \thetab_n^r\bigg)\bigg\|^2\bigg] \notag\\
&~~~ +  \mathbb{E}\bigg[\bigg\|   \frac{1}{NS}\sum_{i=1}^S\sum_{t=1}^N \Bb_{t,i}^\top  \nabla_{\Bb_i\xb} f\bigg(\sum_{t=1}^N \Bb_{t,i}\xb_n^r, \thetab_n^r)\bigg\|^2\bigg]\bigg\} \notag\\
& \leq \frac{1000 [\tilde{F}(\bar \thetab^0, \!\bar\xb^0) \! -\!\underline{F}]}{\sqrt{N ET}}
        \!+\!\frac{14L   \sigma^2}{  \sqrt{N ET} |\mathcal{I}| }
       \!  + \!
        \frac{   222 L^2(1\!+\!7\lambda_w^2)^2 B_{\max}^2 \sigma^2}
            {  \sqrt{N ET}  N (1\!-\!\lambda_{w}^2)^6 |\mathcal{I}|}  \notag\\
    &~~~ +  \frac{2L^2 B_{\max}^2  }
            {\sqrt{NET} N  (1-\lambda_w^2)^4}
            \bigg[\frac{990}{SE N }\phi_z^0
                + \frac{59 (1+7\lambda_w^2)}{(1-\lambda_w^2) T }
            \phi_u^0
            \bigg]   ,
\end{align}
where $\tilde{F}(\bar \xb^0, \bar\thetab^0) \triangleq \frac{1}{NS} \sum_{i=1}^S f(\Bb_i \bar\xb^0, \bar\thetab^0)$, $ \phi_z^0  \triangleq \mathbb{E}\left[ \|\Zb^{r} - \mathbf{1}(\bar{\zb}^{r})^\top\|_F^2\right]$ with $\Zb^r  = [\zb_1^r, \ldots, \zb_N^r]^\top$, and $\phi_u^0$ is defined in the same manner.
\end{Theorem}

Theorem \ref{thm: HBL} can be proved in a similar fashion as that for Theorem 1. The presentation preliminary and proof details are provided in Appendix \ref{appen: alter HBL} to Appendix \ref{appen: thm HBL}. By Theorem \ref{thm: HBL}, one can see that the MUST algorithm also enjoys the same linear speedup with respect to $E$ and $N$.


\section{Experiment Results}
\label{section: simulation results}
In this section, we present the numerical performance of the proposed LSGT  and MUST algorithms.

\vspace{-0.2cm}
\subsection{Evaluation of LSGT}
\subsubsection{\bf{\underline{\emph{ Experiment Setup}}}} \label{sec: experiment HL}
In the experiments,  we set a $20$-agent connected network where they exchange information via the  mixing matrix $\Wb$ obtained by the max-degree rule \cite{pu2021distributed,sayed2014adaptive}. Assume these agents are connected by a random nework graph, which is generated as in \cite{yildiz2008coding}.

We consider the image classification task of $10$ handwritten digits based on the MNIST dataset \cite{MNIST}. It contains $60000$ training images and $10000$ testing images. The dimension of each image is $28\times 28$ which is vectorized into a $1 \times {784}$ vector. These data samples are partitioned to the agents' local dataset in the IID or non-IID fashion.

 $\bullet$ {\bf The IID setting} \cite{briggs2020federated}:
The $60000$ training samples are first shuffled and then assigned evenly to the agents. Thus, each agent holds $3000$ samples with the same class distribution.

 $\bullet$ {\bf The non-IID setting} \cite{mcmahan2017communication}:
We consider a pathological non-IID case where each agent receives $3000$ samples of at most $2$ digits. The training samples are sorted by the digit labels from $9$ to $0$; then they are divided into $40$ subsets so that each subset has $1500$ samples; next every $2$ subsets are randomly assigned to each agent.

To classify the images from MNIST dataset, we consider a 2-layer deep neural network (DNN) with one $30$-neuron hidden layer \cite{wang2021quantized}. The activation function of the hidden and output layer are Rectified Linear Unit (ReLU) and softmax function respectively. The training loss function is cross entropy \cite{Goodfellow-et-al-2016}. The experiment results are averaged over $5$ independent trials.

\vspace{0.2cm}
\subsubsection{\bf{\underline{\emph{ Results and Discussion}}}}

\begin{figure}[!t]
\begin{minipage}[a]{1.0\linewidth}
  \centering
  \epsfig{figure=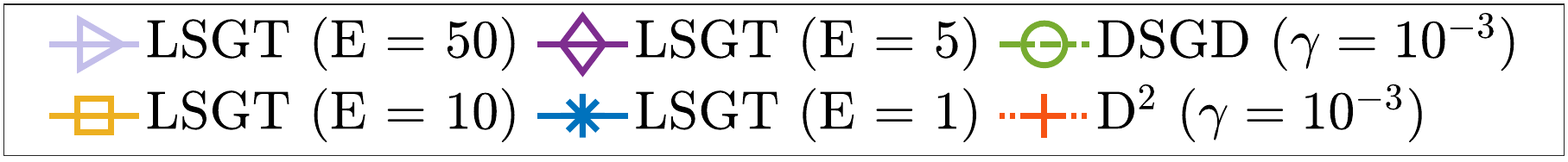,width=3.5in}
\end{minipage}
\begin{minipage}[b]{1.0\linewidth}
  \centering
  \epsfig{figure=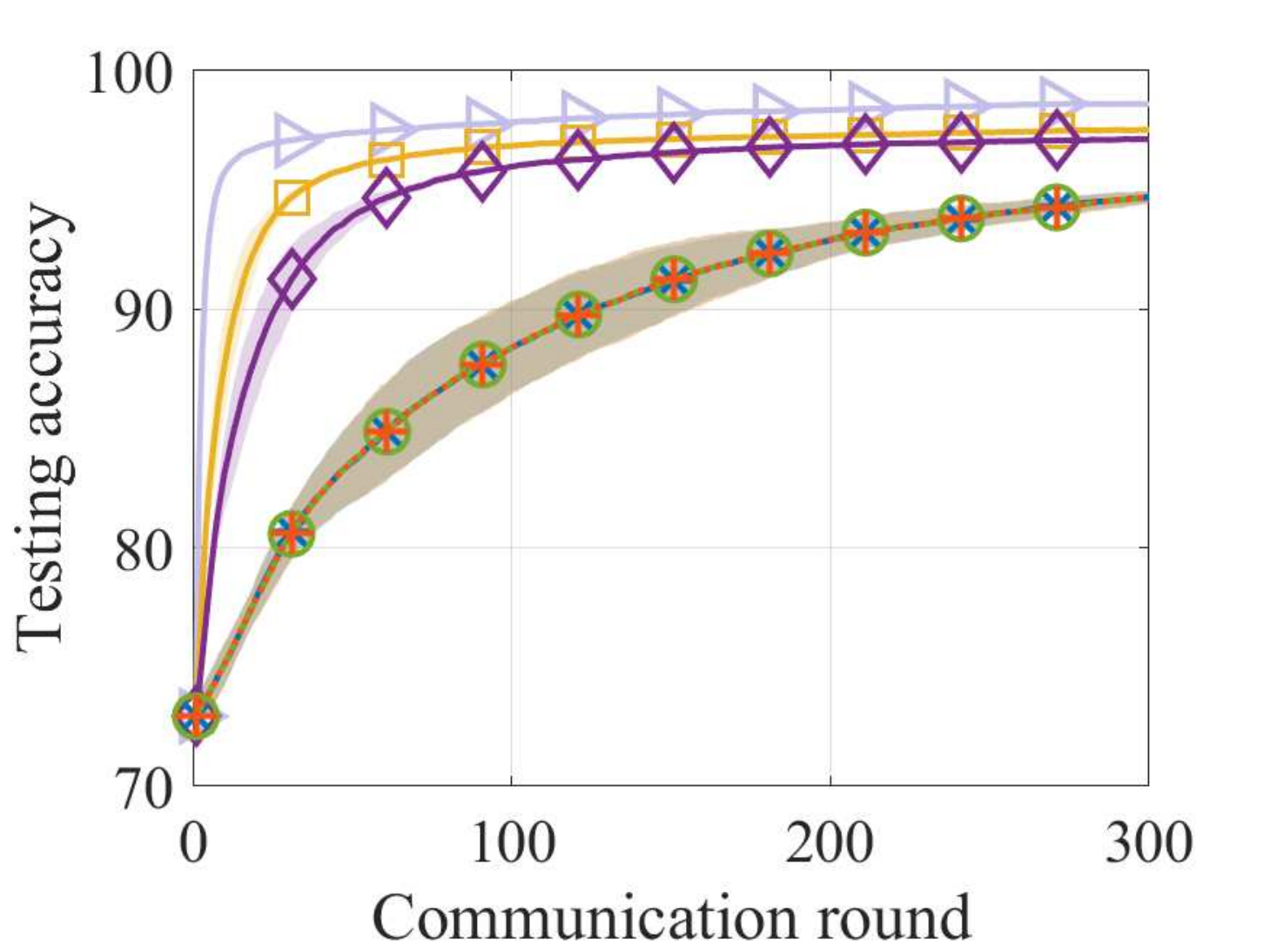,width=2.5in}
 \epsfig{figure=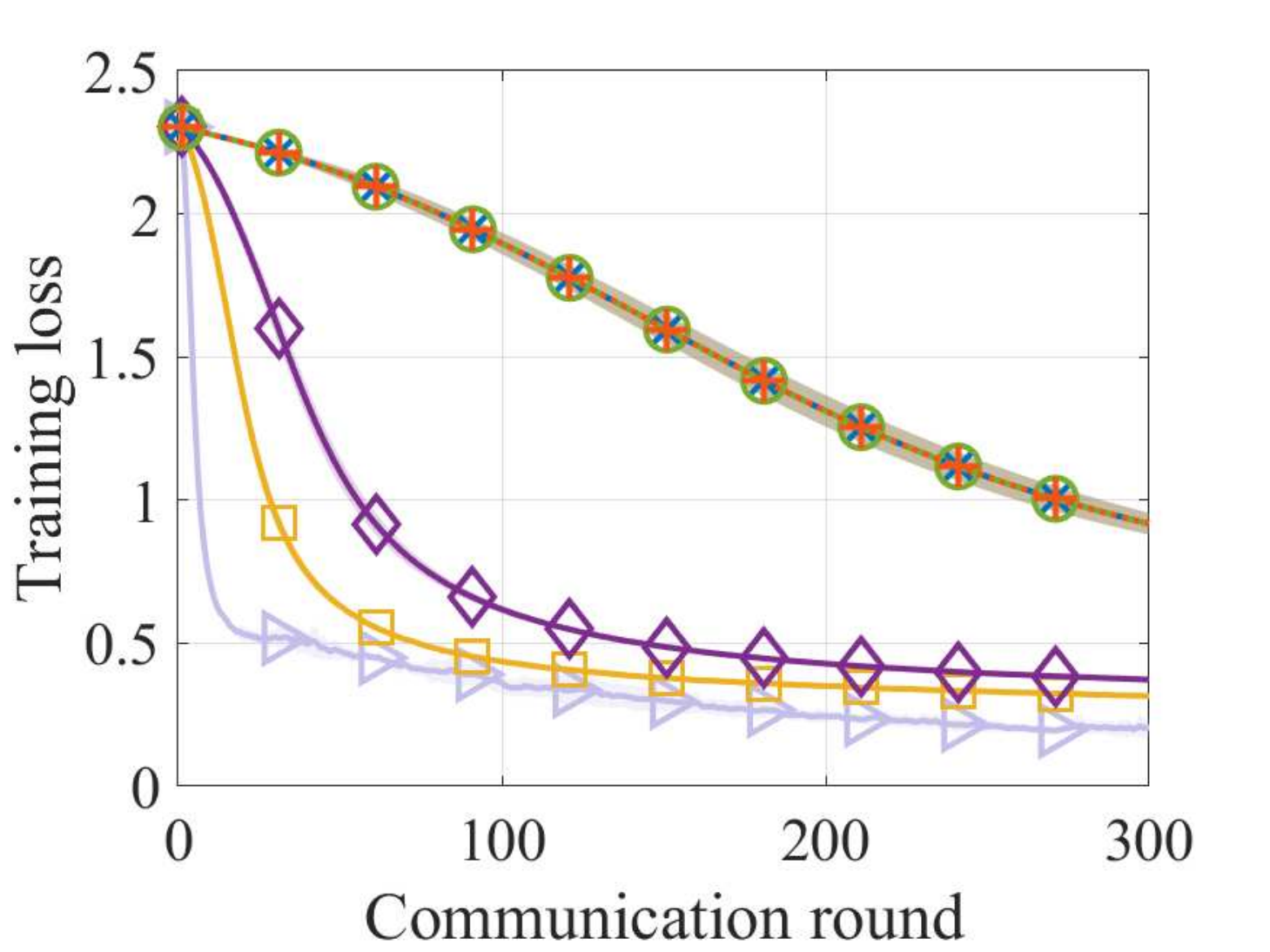,width=2.5in}
  \centerline{\small{(a) The IID setting.}}\medskip
\end{minipage}
\begin{minipage}[a]{1.0\linewidth}
  \centering
  \epsfig{figure=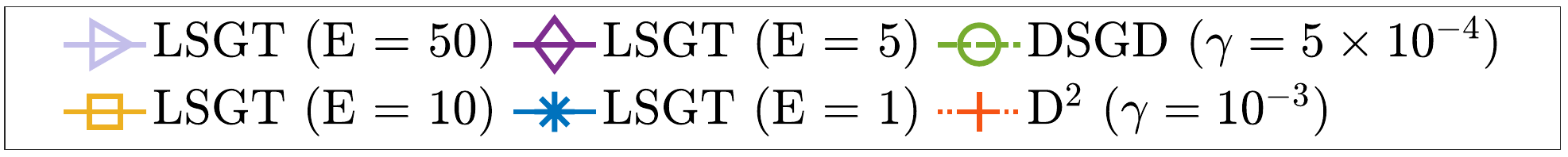,width=3.8in}
\end{minipage}
\begin{minipage}[b]{1.0\linewidth}
  \centering
  \epsfig{figure=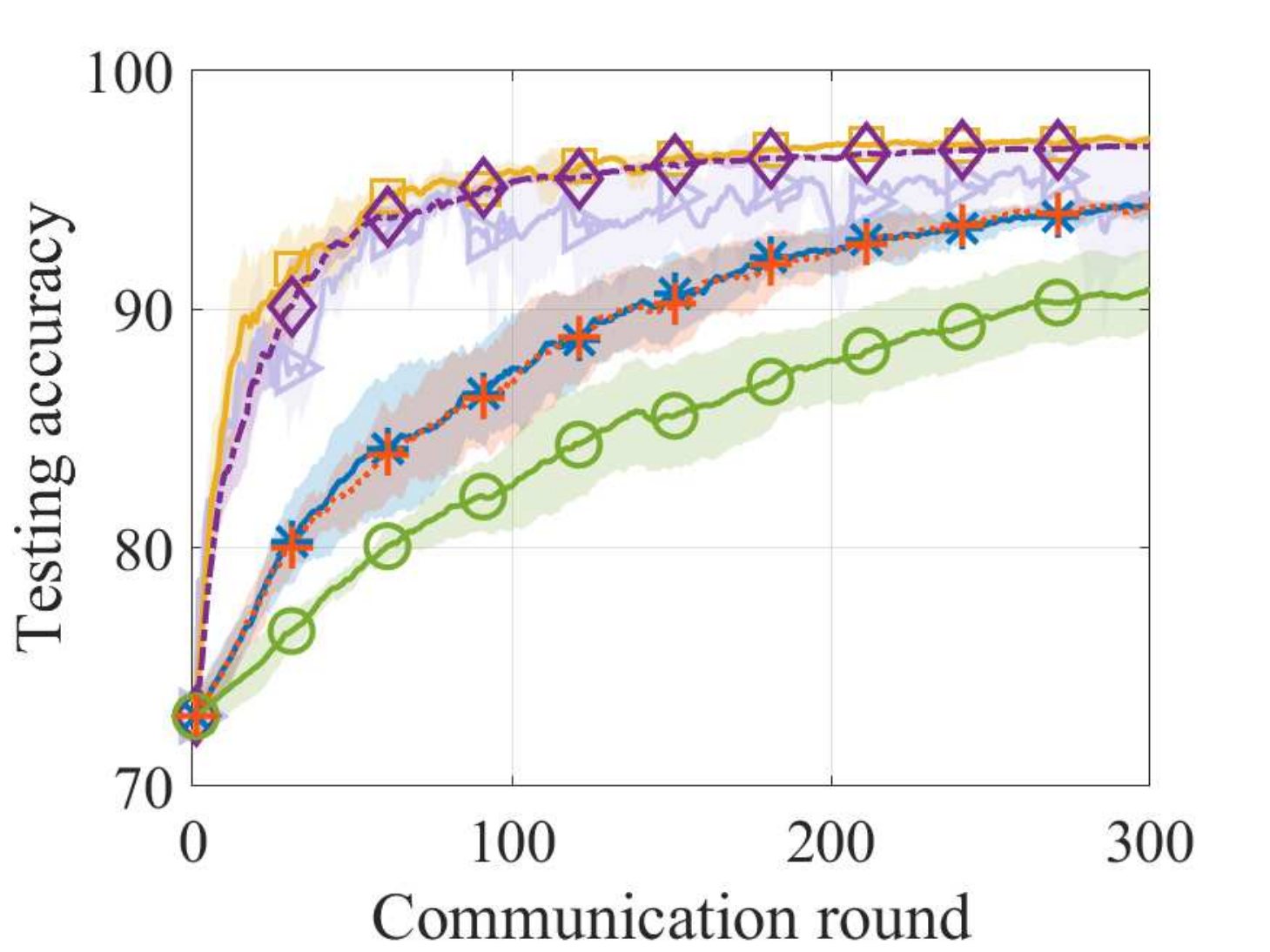,width=2.5in}
 \epsfig{figure=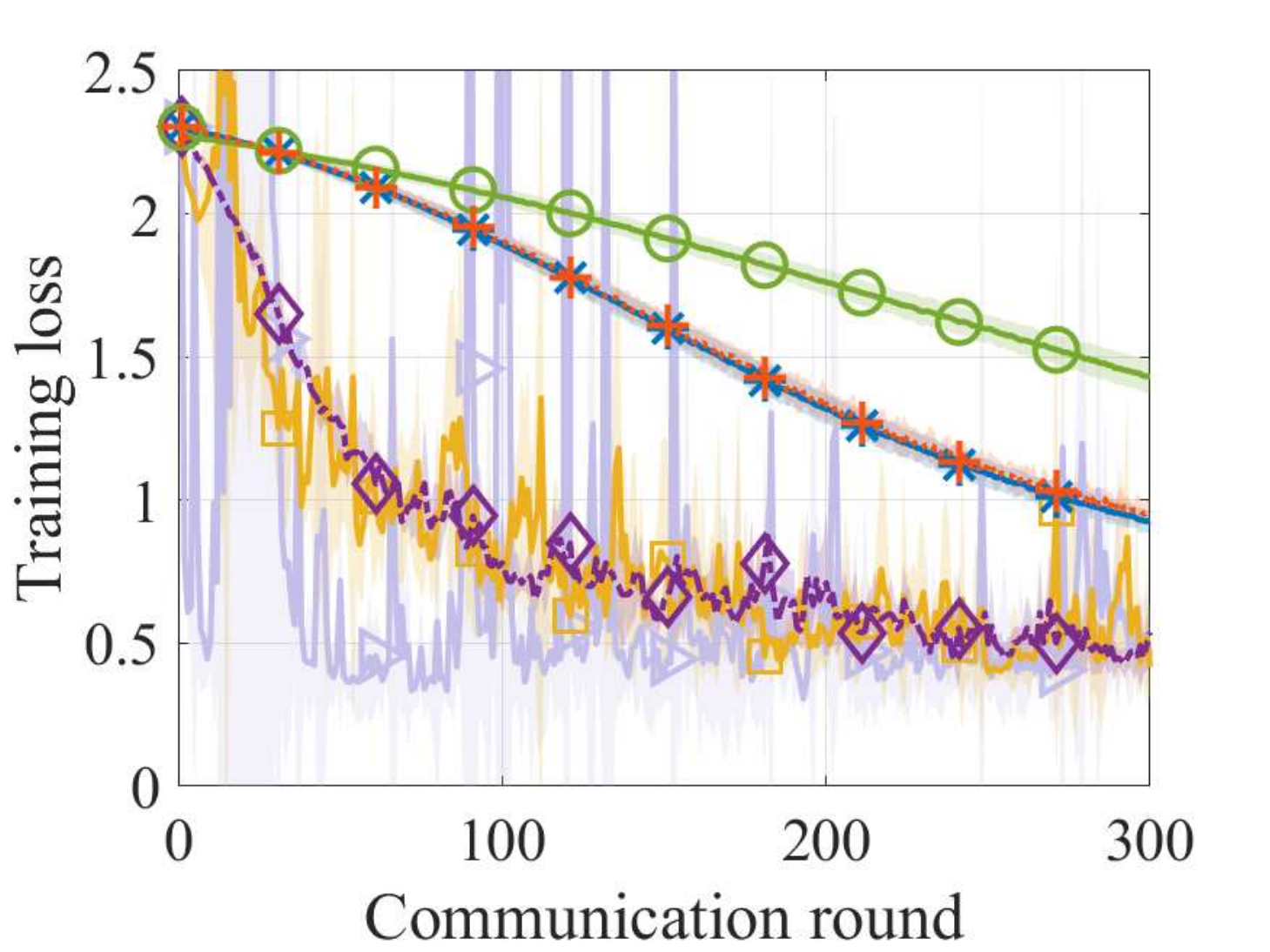,width=2.5in}
  \centerline{\small{(b) The non-IID setting.}}\medskip
\end{minipage}
\caption{\small{Convergence curves of the proposed LSGT algorithm with different $E$ ($\gamma = 10^{-3}$, random graph).}} \label{fig: Q}
\end{figure}

The influences of local iteration number $E$, stepsize $\gamma$, and network topology $\lambda_{w}$ are respectively studied as follows.

$\bullet$ {\bf Impact of $E$:} We set local updates $E=1, 5, 10, 50 $ to investigate the impact of local iteration on the performance of the proposed LSGT algorithm. As shown in Fig.~\ref{fig: Q}(a), the training loss and testing accuracy converges faster as $E$ grows for IID case. It is because in the IID setting the data among decentralized agents are in the same distribution thus more local updates help to learn a better common model. On the other hand, one can see from Fig.~\ref{fig: Q}(b) that the performance of LSGT first improves from $E=1$ to $E=10$, then gets worse when $E=50$. The reason is that a large $E$ enlarges the model discrepancy among agents with non-IID data. More specific relationship between the communication round and local updates $E$ is further investigated in Table~\ref{tbl: Q}. It depicts the communication rounds required by the LSGT algorithm with different $E$ to achieve testing accuracy $85\%,90\%, 95\%$. More local updates reduce the required communication rounds in the IID setting, while the needed communication rounds first decreases and then increases as $E$ grows in the non-IID case. It is consistent with the insight from Theorem \ref{thm: FL} that $E\leq (\frac{T}{N^5})^{\frac{1}{3}}$ should not be too larger. Besides, in Fig.~\ref{fig: Q}, the LSGT algorithm with $E=5$ is superior to DSGD and $\mathrm{D}^2$ methods for both IID and non-IID settings.

\begin{table}[!t]
    \centering
    \caption{{\small \\Communication rounds to achieve a certain testing accuracy.}}
    \tabcolsep=0.1cm
    \label{tbl: Q}
    \belowrulesep=0pt
    \aboverulesep=0pt
    \begin{tabular}{c|cccc|cccc}
        \hline
        \multirow{2}{*}{acc.} & \multicolumn{4}{c|}{The IID setting}  &  \multicolumn{4}{c}{ The non-IID setting} \\
        & $E\!=\!1$ & $E\!=\!5$ & $E\!=\!10$ &$E\!=\!50$ &
        $E\!=\!1$ & $E\!=\!5$ & $E\!=\!10$ &$E\!=\!50$\\
        \hline
        $85\%$ & $62$ & $14$ & $8$ & $\mathbf{3}$ &
        $75$ & $16$ & $\mathbf{10}$ & $11$\\
        $90\%$ & $126$ & $26$ & $14$ & $\mathbf{4} $&
        $136$ & $30$ & $\mathbf{21}$ & $41$\\
        $95\%$ & $335$ & $68$ & $35$ & $\mathbf{8}$ &
        $378$ & $91$ & $\mathbf{65}$ & $191$\\
        \hline
    \end{tabular}\\
    ``acc." denotes testing accuracy.
    \vspace{-0.2cm}
\end{table}

\begin{figure}[!t]
\begin{minipage}[a]{1.0\linewidth}
  \centering
  \epsfig{figure=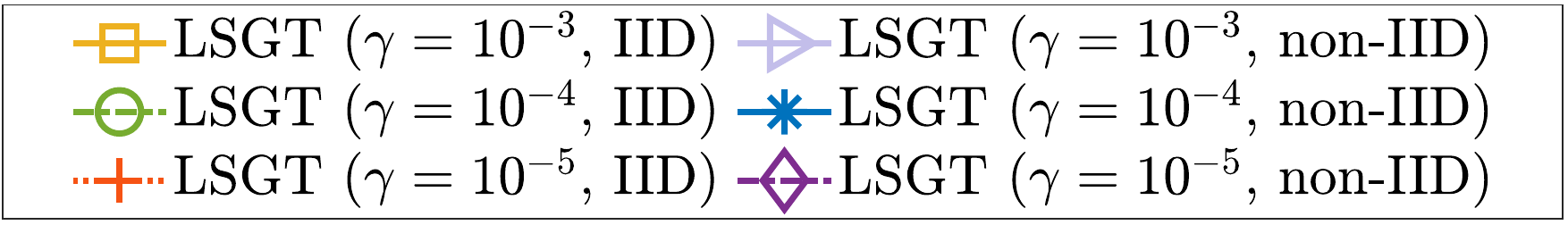,width=3.5in}
\end{minipage}
\begin{minipage}[b]{1.0\linewidth}
  \centering
  \epsfig{figure=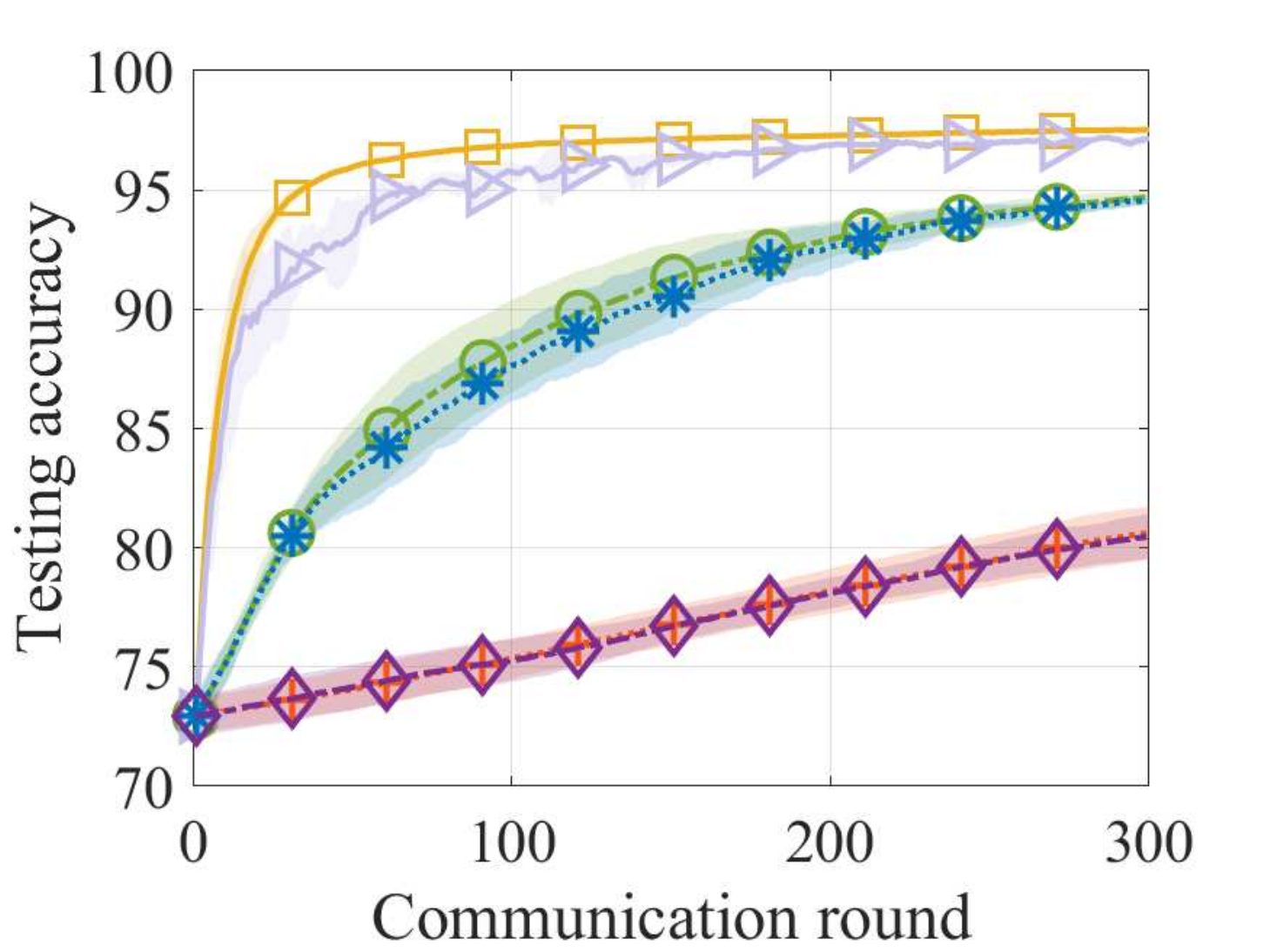,width=2.5in}
 \epsfig{figure=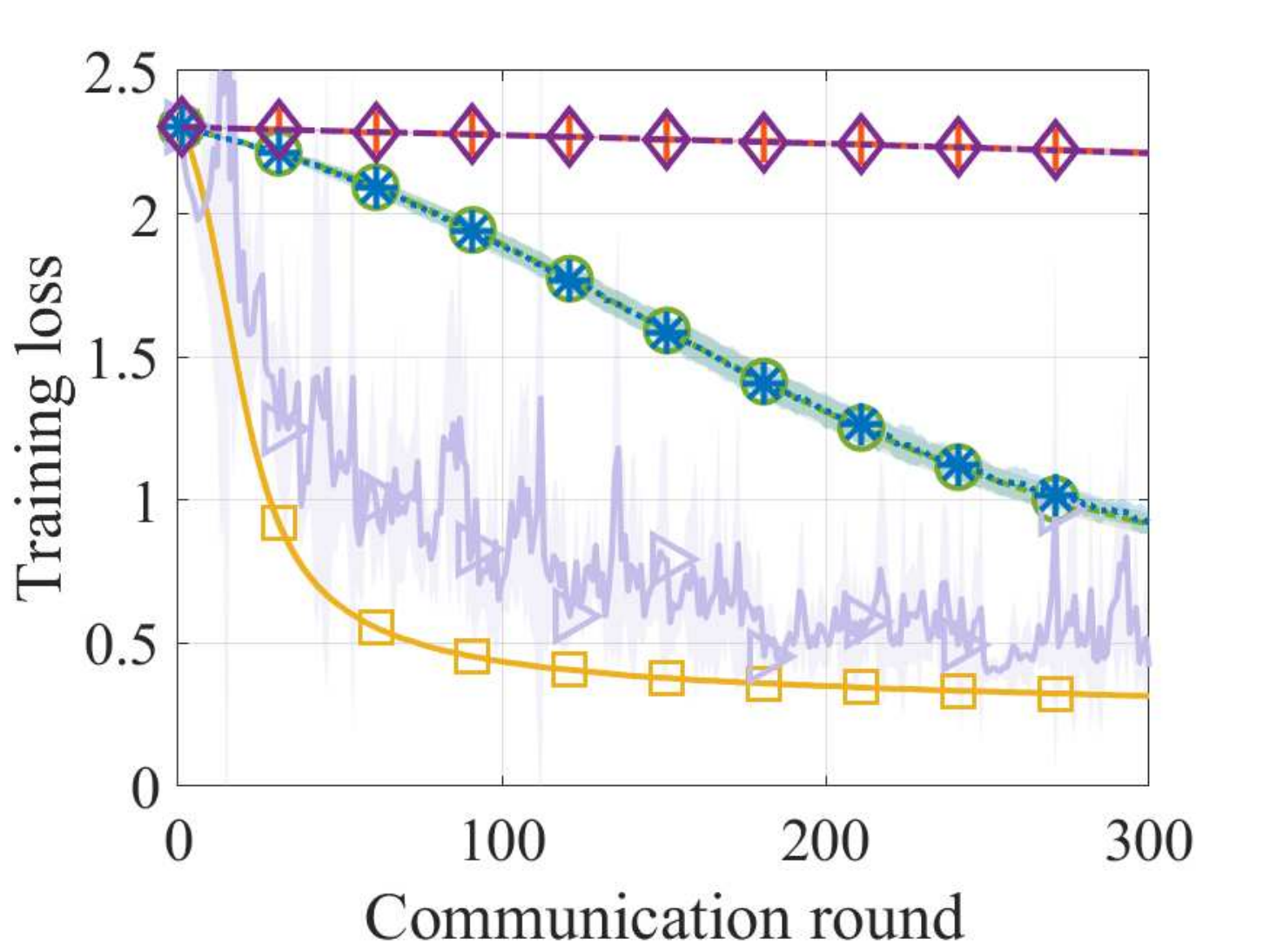,width=2.5in}
\end{minipage}
\caption{\small{Convergence curves of the proposed LSGT algorithm with different stepsize ($E=10$, random graph).}} \label{fig: step}
\end{figure}

\begin{figure}[!t]
\begin{minipage}[a]{1.0\linewidth}
  \centering
  \epsfig{figure=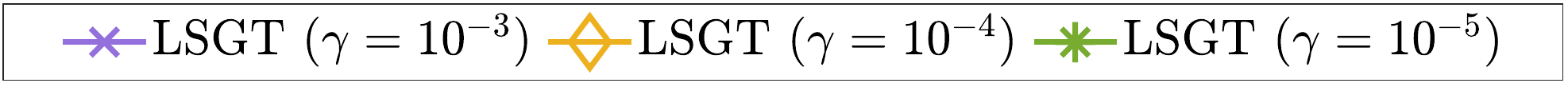,width=3.8in}
\end{minipage}
\begin{minipage}[b]{1.0\linewidth}
  \centering
  \epsfig{figure=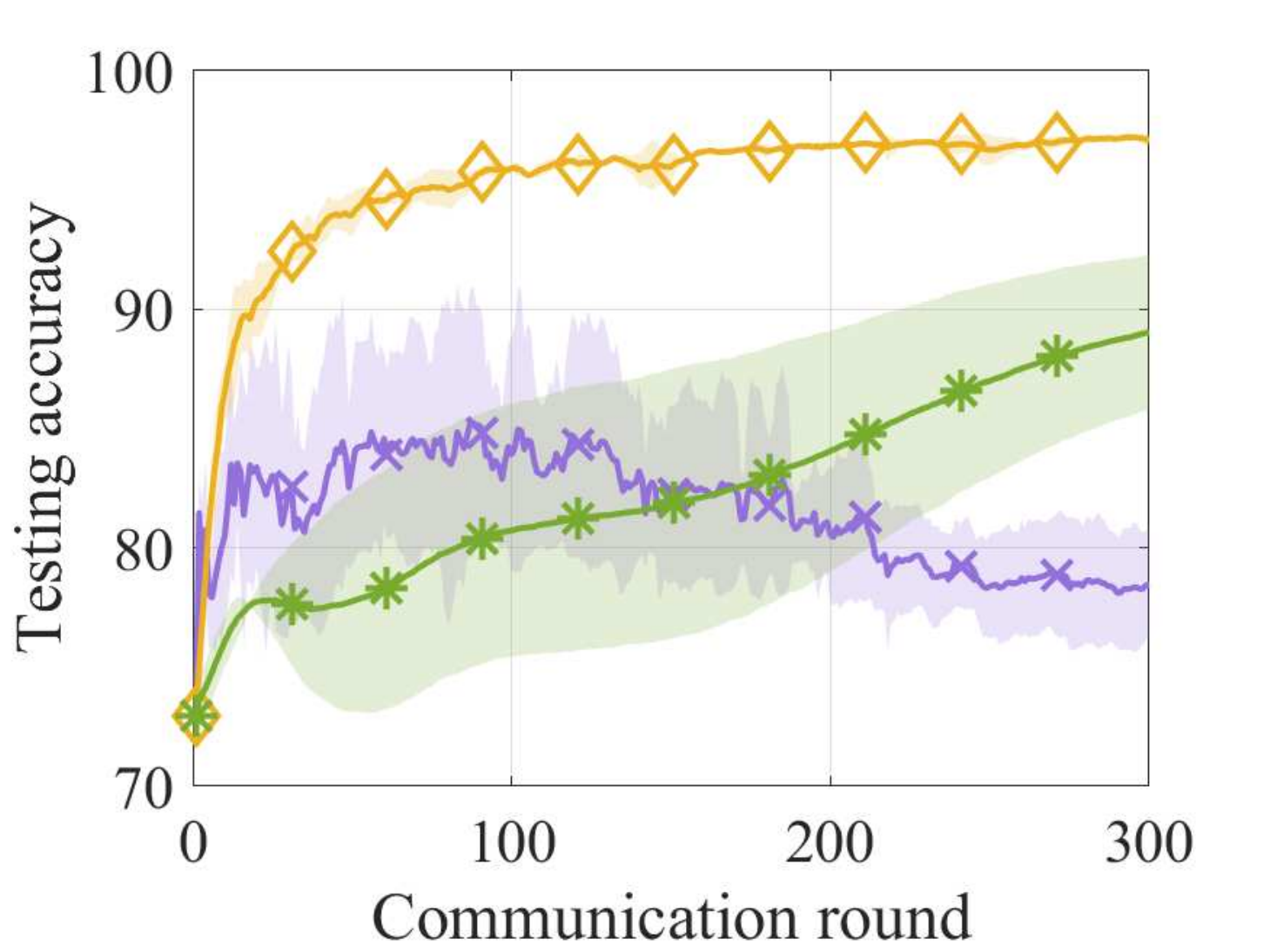,width=2.5in}
 \epsfig{figure=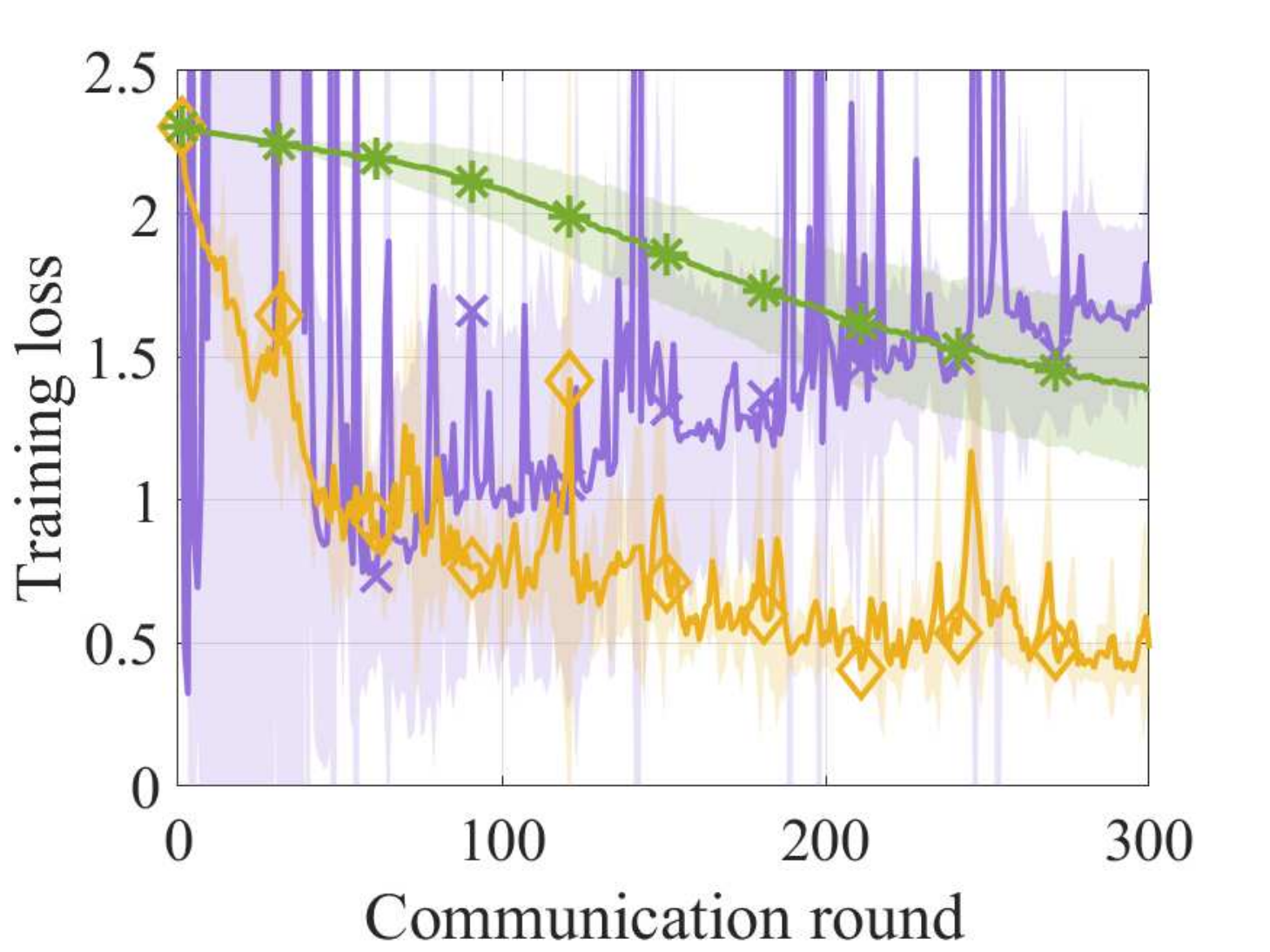,width=2.5in}
\end{minipage}
\caption{\small{Convergence curves of the proposed LSGT algorithm with different stepsize in the non-IID setting ($E=100$, random graph).}} \label{fig: Q100}
\vspace{-0.2cm}
\end{figure}

$\bullet$ {\bf Impact of $\gamma$:} It is found from Fig.~\ref{fig: step} and Fig.~\ref{fig: Q100} that for the IID and non-IID settings, the proposed LSGT algorithm with a smaller stepsize such as $\gamma = 10^{-5}$ converges slower in both training loss and testing accuracy, which is consistent with the analysis in Theorem \ref{thm: FL}. In Fig.~\ref{fig: Q100}, it should be pointed out that the LSGT algorithm with $E=100, \gamma = 10^{-3}$ cannot converge well. Recalling the stepsize setting $\gamma =  \sqrt{{N}/{ ET }}$ from Corollary \ref{coro: coro 2 Q larger than 1}, for a large $E=100$ we can consider reducing the stepsize. As observed in Fig.~\ref{fig: Q100}, when $E=100, \gamma = 10^{-4}$, the LSGT algorithm converges to a higher testing accuracy.

\begin{figure}[!t]
\begin{minipage}[a]{1.0\linewidth}
  \centering
  \epsfig{figure=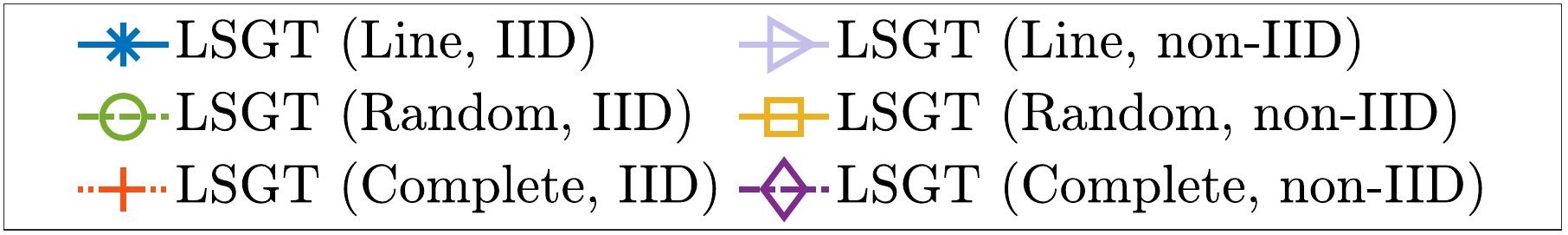,width=3.5in}
\end{minipage}
\begin{minipage}[b]{1.0\linewidth}
  \centering
  \epsfig{figure=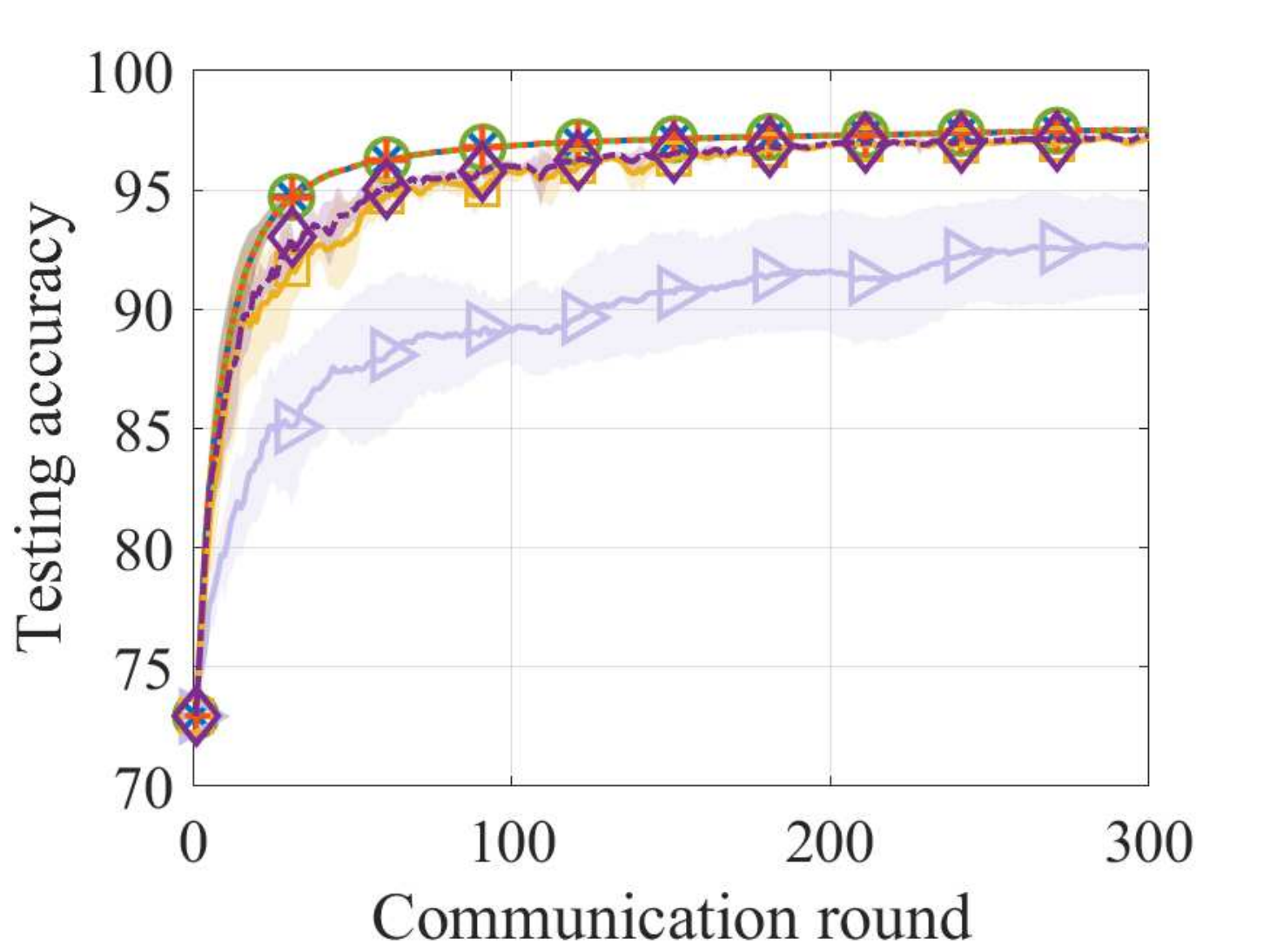,width=2.5in}
 \epsfig{figure=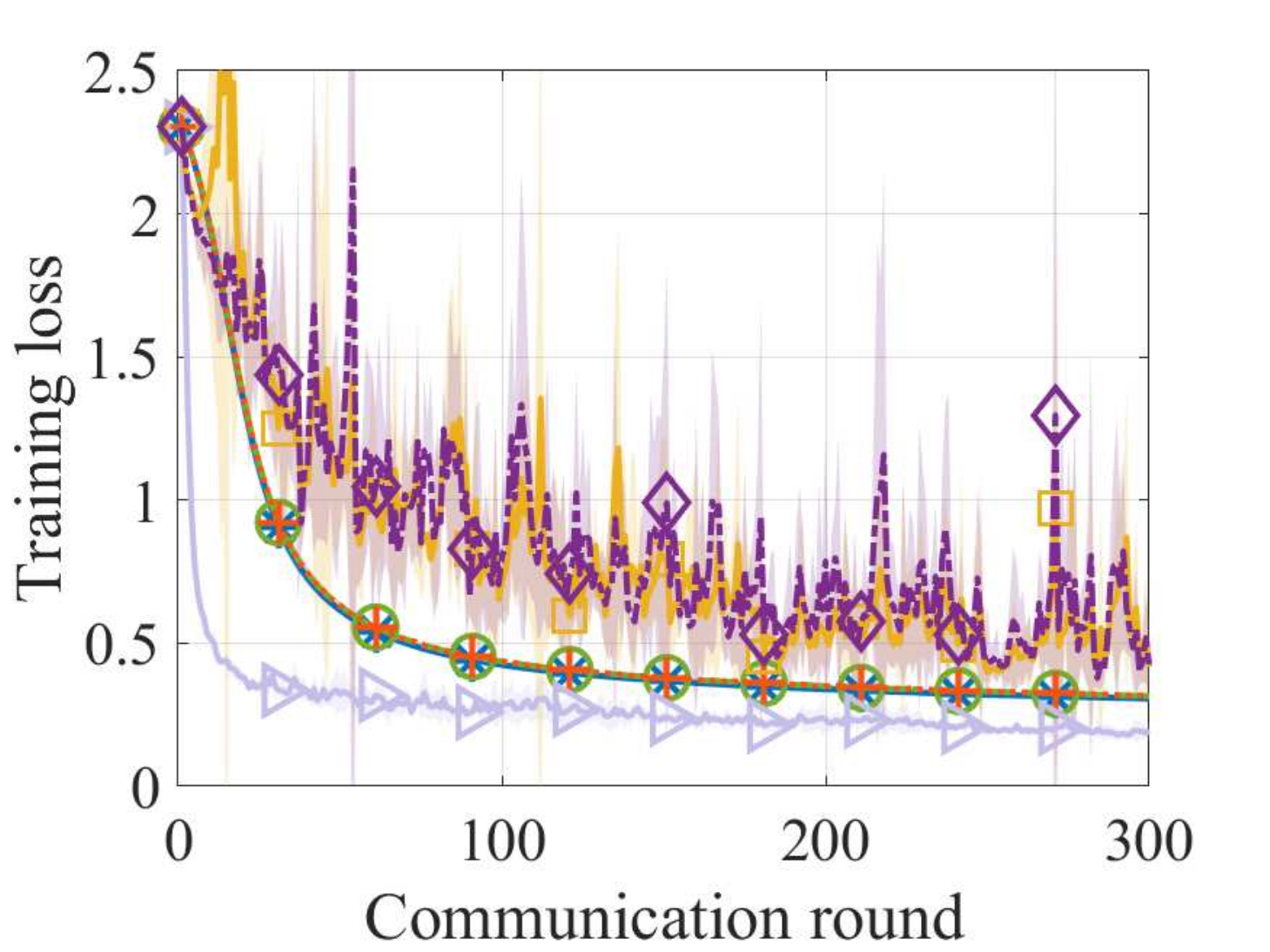,width=2.5in}
\end{minipage}
\caption{\small{Convergence curves of the proposed LSGT algorithm with different topologies ($E=10$, $\gamma = 10^{-3}$).}} \label{fig: topo}
\vspace{-0.2cm}
\end{figure}
$\bullet$ {\bf Impact of $\lambda_w$:} To investigate the influence of the network connectivity,
we consider $3$ graphs with order $\lambda_{w}(\text{line}) > \lambda_{w}(\text{random}) >\lambda_{w}(\text{complete})$, whose connectivity order is reversed. In Fig.~\ref{fig: topo}, for the IID setting, the network topology has a minor impact on the convergence performance. In the non-IID case, when the network has better connectivity (complete graph), the proposed method achieves a higher testing accuracy. It is in agreement with the analysis in Remark \ref{rmk: tau FL}. Interestingly, the proposed LSGT algorithm with stronger connectivity converges to a larger training loss. The reason is that as the connectivity of the network topology grows, the training model gets closer to the global model. When the training model is applied into local data, it may not work well owing to heterogeneity of local data \cite{kulkarni2020survey}.

\subsection{Evaluation of MUST algorithm}
\begin{figure}[!t]
\begin{minipage}[a]{1.0\linewidth}
  \centering
  \epsfig{figure=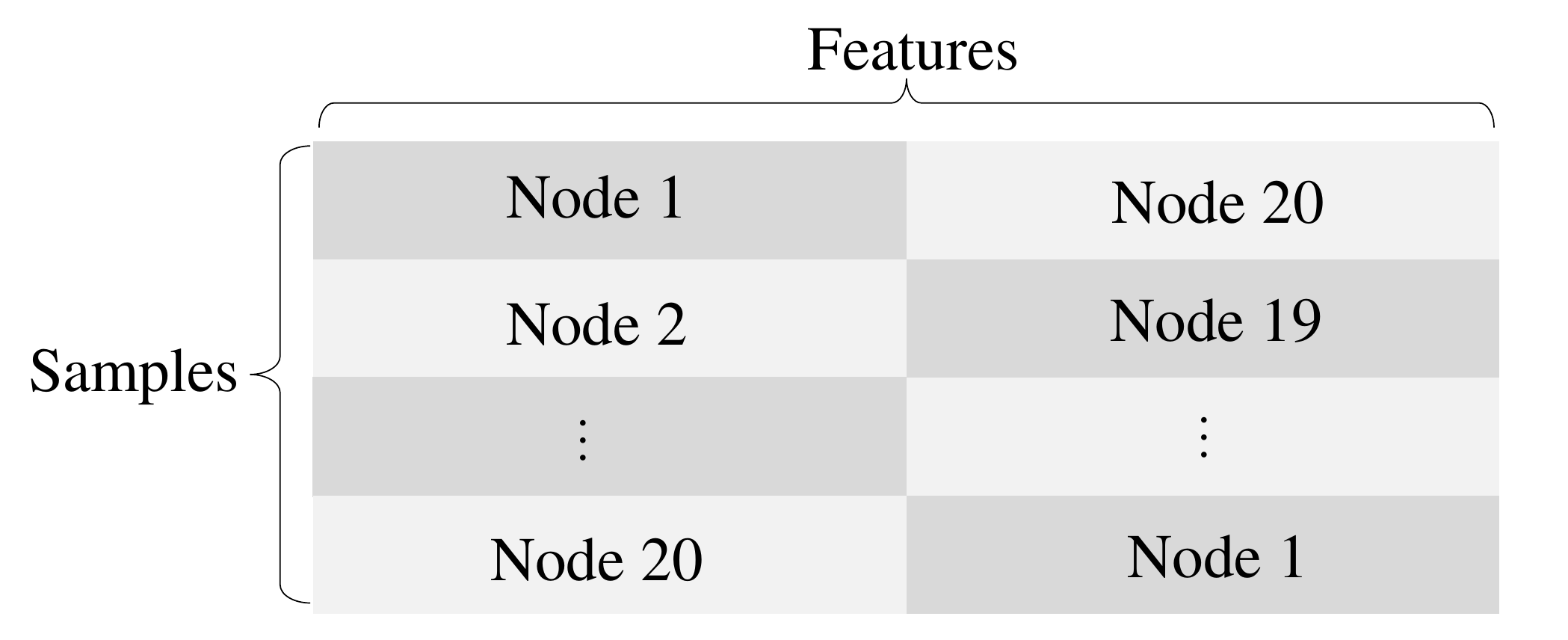,width=3.5in}
\end{minipage}
\caption{\small{Data partition for the hybrid data case.}} \label{fig: dp}
\vspace{-0.2cm}
\end{figure}
\begin{figure}[!t]
\begin{minipage}[a]{1.0\linewidth}
  \centering
  \epsfig{figure=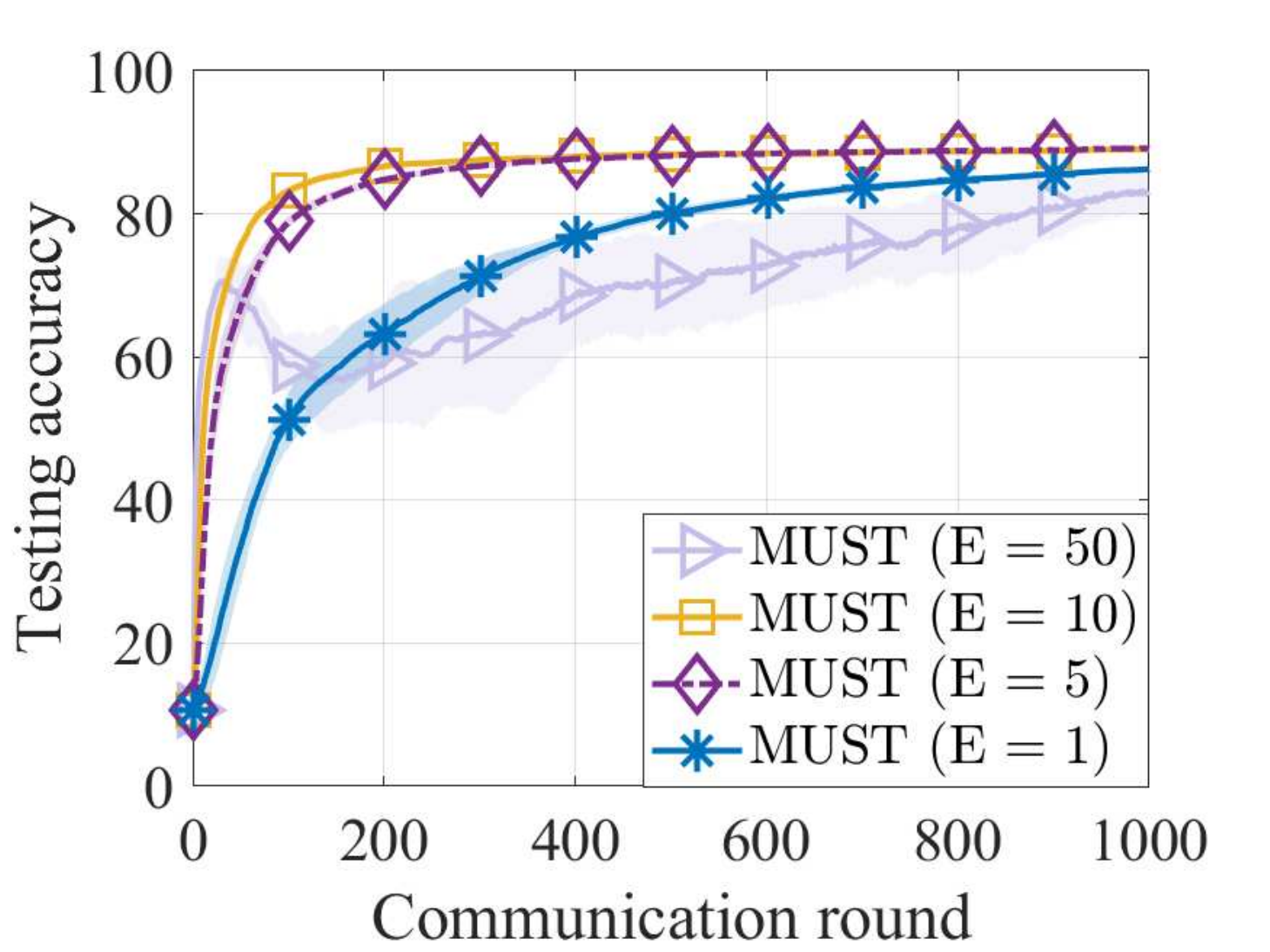,width=2.5in}
 \epsfig{figure=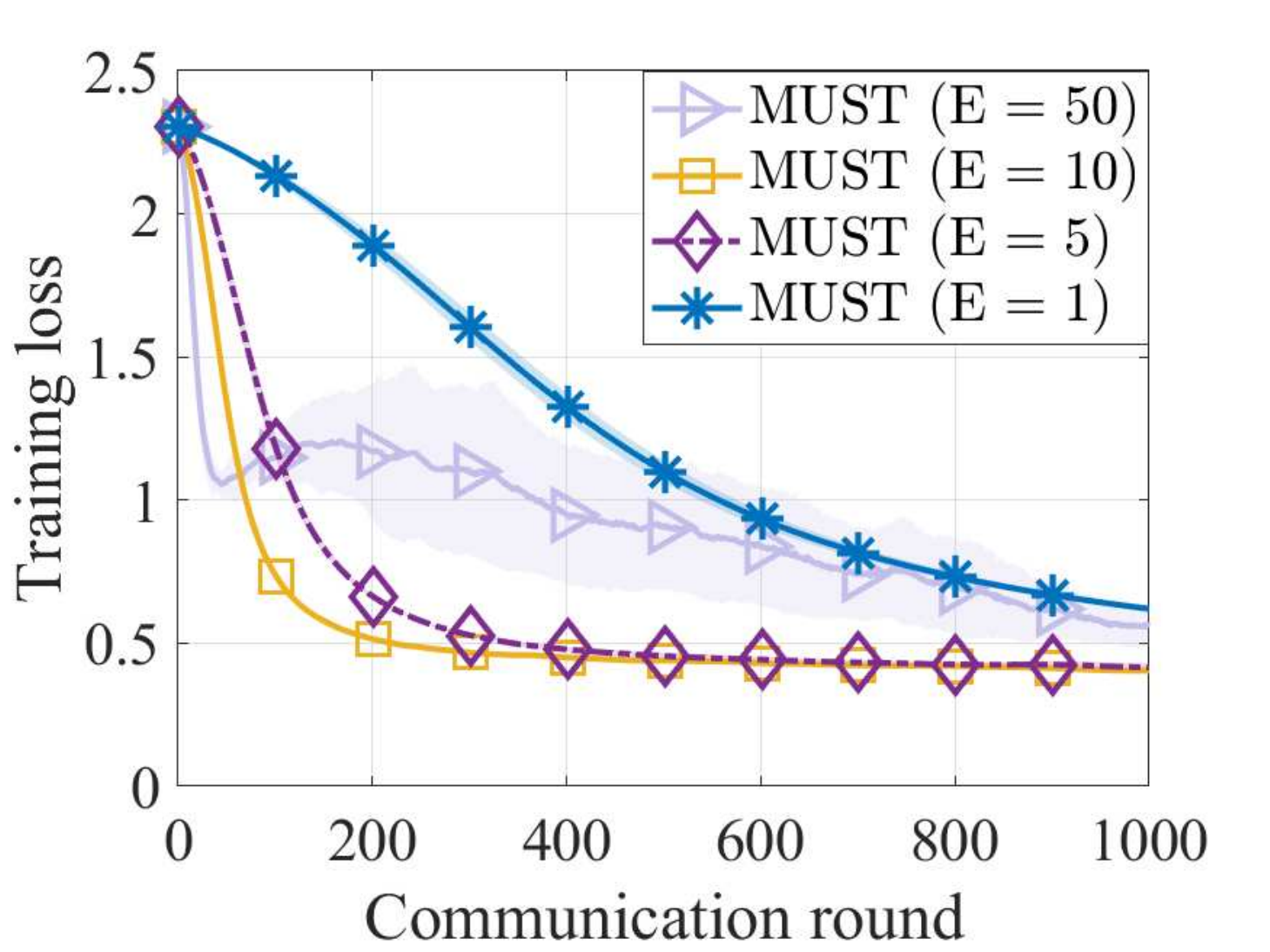,width=2.5in}
\end{minipage}
\caption{\small{Convergence curves of the proposed MUST algorithm with different $E$ ($\gamma = 10^{-3}$, random graph).}} \label{fig: HBL Q}
\vspace{-0.2cm}
\end{figure}
To examine the effectiveness of the proposed MUST algorithm, we consider the same random $20$-agent network and DNN as in Sec.~\ref{sec: experiment HL}. In the hybrid data setting, the training samples are shuffled first and then partitioned into $3000$ subsets with size $20$. Then, each subsection is distributed as in Fig.~\ref{fig: dp}. Specifically, in each $20$-sample subset, each image is divided into $2$ patches. The first patches of each subset are assigned to $20$ agents following a positive sequence, and the second patches are allocated in an inverted order. Thus, each agent has part of samples with incomplete features.

We evaluate the performance of the MUST algorithm with respect to the number of local iterations. It is set to $E=1, 5, 10, 50$. As shown in Fig.~\ref{fig: HBL Q}, the proposed MUST method performs robustly over the hybrid data. Particularly, as $E$ increases $1$ to $10$, the testing accuracy and training loss obtained by the MUST algorithm improves faster with the iterations. Under a larger $E=50$, as expected, the MUST algorithm cannot learn an effective DNN to classify the digits and even fluctuates around $\text{acc.}~ = 80\%$. 

\vspace{-0.2cm}
\section{Conclusion}\label{section: conclusions}
In this paper, to improve the communication efficiency of the existing GT method for solving problem \eqref{prob: HL 2}, we have proposed a new LSGT algorithm (Algorithm \ref{alg: LSGT}) which incorporates the local SGD technique into the stochastic GT method. Theoretically, we have built the convergence conditions of the LSGT algorithm (Theorem \ref{thm: FL}) and shown that it can benefit the linear speedup with the number of local SGD updates $E$ (Corollary \ref{coro: coro 2 Q larger than 1}). This is a strong contrast to the existing GT methods which either did not consider multiple steps of local SGD or not fully characterize the merit of local SGD for the GT methods. As an extension, we have also extended the gradient and variable tracking idea to develop the MUST algorithm (Algorithm \ref{alg: MUST}), for solving the hybrid-data learning problem \eqref{pro: HFL rewrite}.
The presented experiment results have shown that the proposed LSGT method yields better learning performance than one-step stochastic GT method, and the MUST algorithm performs robustly over heterogeneous hybrid data.
It is worthwhile to point out that the proposed LSGT and MUST methods, to the best of our knowledge, are the first stochastic GT algorithm with multiple local SGD updates for decentralized learning. In the future, we plan to extend this framework to time-varying communication network topologies \cite{scutari2019distributed} and that with compression \cite{liao2022compressed}.

\appendix

\section{Alternative Expression of Algorithm \ref{alg: LSGT}}
For ease of analysis, we can conclude the \eqref{eq: LSGT q 2} for $n\in [N]$ at $q$-th iteration of the $r$-th round as
\begin{subequations}\label{eq: update r q}
\begin{align}
&\yb_n^{r,q}  =\sum_{m=1}^N W_{n,m} \yb_m^r - \gamma \sum_{m=1}^N W_{n,m} \vb_m^r - \gamma   \sum_{k=1}^{q-1}\vb_n^{r, k}  \label{eq: y update r q},  \\
& \vb_n^{r,q}  =\sum_{m=1}^N W_{n,m} \vb_m^r
+   \gb_n^{r, q}
-   \gb_n^{r}   , \label{eq: v update r q}
\end{align}
\end{subequations}
where $\gb_n^r  \triangleq  \frac{1}{|\mathcal{I}|}\sum_{\xi   \in  \mathcal{I}_n^{r} }  g_n ( \yb_n^{r}, \xi  )$ defined in \eqref{eq: h def}, $q\in [E]$, and $r\in [T]$.
Substituting $q = E$ into \eqref{eq: update r q}, and by $\yb_n^{r+1}  = \yb_n^{r,E}$ and $\vb_n^{r+1}=  \vb_n^{r,E}$ in Algorithm \ref{alg: LSGT}, we have
\begin{subequations}
\begin{align}
 &\yb_n^{r+1}  =  \sum_{m=1}^N\! W_{n,m} \yb_m^r \!- \!\gamma  \sum_{m=1}^N \!W_{n,m} \vb_m^r \!- \!\gamma  \sum_{q=1}^{E\!-\!1}\vb_n^{r, q}  \label{eq: y update r+1}, \\
&\vb_n^{r+1}=   \sum_{m=1}^N W_{n,m} \vb_m^r
+    \gb_n^{r+1}
-   \gb_n^{r}  , \label{eq: v update r+1}
\end{align}
\end{subequations}
for any $ n\in[N]$.

Making average on the above two equalities, we obtain
\begin{subequations}
\begin{align}
 &\bar\yb^{r+1}  =   \bar\yb^r - \gamma   \bar\vb^r - \gamma   \cdot \frac{1}{N}\sum_{n=1}^N \sum_{q=1}^{E-1}\vb_n^{r, q}  \label{eq: y update ave}, \\
&\bar\vb^{r+1}=   \bar\vb^r
+   \bar{\gb}_n^{r+1}
-   \bar{\gb}_n^{r}  ,   \label{eq: v update ave}
\end{align}
\end{subequations}
where $ \bar{\gb}^{r} \triangleq \textstyle \frac{1}{N} \sum_{n=1}^N \gb_n^r$.

Similar to  \cite[Eqn. (27)]{YingTSP2019}, one can use induction to show:
\vspace{-0.2cm}
\begin{Lemma} \label{lem: ave = true sum}
  For all $ n\in [N]$, initializing $\vb_n^0 =  \gb_n^0$, by \eqref{eq: v update r+1}, we have
    \begin{align}
        & \bar \vb^r \triangleq  \frac{1}{N}  \sum _{n=1}^N \vb^r
       , ~~    \bar \vb^r=  \bar{\gb}^{r} , \forall r\geq 0.
        \label{eq: ave = true sum}
    \end{align}
\end{Lemma}
Then, based on Assumption \ref{assum: stochastic}, we have
\begin{align}
    & \mathbb{E}[\bar\vb^r] =  \mathbb{E}[\bar\gb^r] = \frac{1}{N} \sum_{n=1}^N \nabla f_n(\yb_n^r),\label{eq: unbiased diff bar v true gr} \\
    & \mathbb{E}[\|\bar\vb^r\|^2]  = \mathbb{E}[\|\bar\gb^r\|^2]\notag\\
    & = \mathbb{E}\bigg[\bigg\|\bar\gb^r - \frac{1}{N} \sum_{n=1}^N \nabla f_n(\yb_n^r)  + \frac{1}{N} \sum_{n=1}^N \nabla f_n(\yb_n^r)\bigg\|^2\bigg]\notag\\
    &  \leq  2\mathbb{E}\bigg[\bigg\|\bar\gb^r \!- \!\frac{1}{N} \sum_{n=1}^N\! \nabla f_n(\yb_n^r) \bigg\|^2\bigg]
        \! +\!2\mathbb{E}\bigg[\bigg\| \frac{1}{N} \sum_{n=1}^N \!\nabla f_n(\yb_n^r)\bigg\|^2\bigg]\notag\\
     &= \frac{2 \sigma ^2}{N |\mathcal{I}|}  + 2\mathbb{E}\bigg[\bigg\| \frac{1}{N} \sum_{n=1}^N \nabla f_n(\yb_n^r)\bigg\|^2\bigg].\label{eq: diff bar v true gr}
\end{align}

\section{Proof of Lemma \ref{lem: y r q}}\label{appen: y r q}
In the following, we prove \eqref{eq: y r q} and \eqref{eq: v r q}, respectively.

\noindent$\bullet$ {\bf Proof of \eqref{eq: y r q}:} We have
\begin{align}\label{eq: y rq 1}
&\sum_{n=1}^N \sum_{q=1}^{E-1}\mathbb{E} [\|\yb_n^{r, q} - \bar\yb^r\|^2 ]\notag\\
& \overset{\mathrm{(i)}}{=}\! \sum_{n=1}^N \! \sum_{q=1}^{E-1}\!\mathbb{E}\bigg[\bigg\|
    \sum_{m=1}^N \!W_{n,m} (\yb_m^r -\gamma \vb_m^r ) \!- \!  \bar\yb^r
     \!- \!\gamma    \sum_{k=1}^{q\!-\!1}\vb_n^{r, k}
    \bigg\|^2\bigg]\notag\\
& =  \sum_{n=1}^N \sum_{q=1}^{E\!-\!1}\mathbb{E}\bigg[\bigg\|
    \sum_{m=1}^N W_{n,m} \yb_m^r \! -\! \bar\yb^r  \! - \!\gamma     \sum_{m=1}^N \! \!W_{n,m} \vb_m^r \!- \!\gamma   \bar\vb^r
    - \gamma    \sum_{k=1}^{q-1}(\vb_n^{r, k}-\bar\vb^r)  - q \gamma   \bar\vb^r
    \bigg\|^2\bigg]\notag\\
& \overset{\mathrm{(ii)}}{\leq} 4 \sum_{n=1}^N \sum_{q=1}^{E-1}\mathbb{E}\bigg[\bigg\|
    \sum_{m=1}^N W_{n,m} \yb_m^r- \bar\yb^r
    \bigg\|^2\bigg]
    + 4   \gamma^2 \sum_{n=1}^N \sum_{q=1}^{Q-1}\mathbb{E}\bigg[\bigg\|
    \sum_{m=1}^N W_{n,m} \vb_m^r- \bar\vb^r
    \bigg\|^2\bigg]\notag\\
    &~~~
    + 4   \gamma^2 \sum_{n=1}^N \sum_{q=1}^{Q-1}\mathbb{E}\bigg[\bigg\|  \sum_{k=1}^{q-1}(\vb_n^{r, k}- \bar\vb^r)\bigg\|^2\bigg]
    + 4   \gamma^2 \sum_{n=1}^N \sum_{q=1}^{Q-1}q^2 \mathbb{E} [ \|  \bar\vb^r \|^2 ]\notag\\
& \leq 4 (E-1)\mathbb{E}\bigg[\bigg\| \Wb  \Yb^r - \mathbf{1}(\bar\yb^r)^\top
    \bigg\|_F^2\bigg] ~
    +4(E-1)\lambda_w^2   \gamma^2 \mathbb{E}\bigg[\bigg\| \Vb^r - \mathbf{1}(\bar\vb^r)^\top
    \bigg\|_F^2\bigg]\notag\\
    &~~~
    + 4  \gamma^2 \sum_{n=1}^N \sum_{q=1}^{E-1}(q-1)\sum_{k=1}^{q-1}\mathbb{E} [ \| \vb_n^{r, k}- \bar\vb^r \|^2 ]
    + 4  \gamma^2  N \sum_{q=1}^{E-1}q^2 \mathbb{E} [\|  \bar\vb^r \|^2]\notag\\
& \overset{\mathrm{(iii)}}{\leq} 4(E-1)\lambda_{w}^2 \phi_{y}^r + 4(E-1)\lambda_w^2  \gamma^2\phi_{v}^r
        +4   \gamma^2 \frac{(E-1)E(2E-1)}{6} N \mathbb{E} [ \|  \bar\vb^r \|^2 ]\notag\\
    & ~~~ +4  \gamma^2 \sum_{n=1}^N \sum_{q=1}^{E-1} \frac{(q+(E-1))(E-q)}{2} \mathbb{E} [ \| \vb_n^{r, q}- \bar\vb^r \|^2 ]\notag\\
& \overset{\mathrm{(iv)}}{\leq} 4(E-1)\lambda_{w}^2 \phi_{y}^r + 4(E-1)\lambda_w^2   \gamma^2\phi_{v}^r
        + 2E^2(E-1)N \gamma^2 \mathbb{E} [ \|  \bar\vb^r \|^2 ]\notag\\
    & ~~~+ 2E(E-1)   \gamma^2 \sum_{n=1}^N \sum_{q=1}^{E-1} \mathbb{E} [ \| \vb_n^{r, q}- \bar\vb^r \|^2 ],
\end{align}
where $\mathrm{(i)}$ is by \eqref{eq: y update r q}; $\mathrm{(ii)}$ is from Jensen's inequality $\textstyle\frac{\sum  _{i=1}^n f(y_i)}{n} \geq \textstyle f \Big(\frac{\sum _{i=1}^n y_i}{n} \Big)$ with $f(\cdot) = \|\cdot\|^2$; the first term of $\mathrm{(iii)}$ is obtained by the fact that $\|\Wb\Yb^r - \mathbf{1}({\bar\yb^r})^\top\|_F^2 = \textstyle \|(\Wb - \mathbf{1}\mathbf{1}^\top/N)[\Yb^r -  \mathbf{1}({\bar\yb^r})^\top]\|_F^2 \leq \lambda_{w}^2 \| \Yb^r -  \mathbf{1}({\bar\yb^r})^\top\|_F^2 $; the second term of $\mathrm{(iii)}$ comes from $\textstyle \sum _{i=1}^n i^2 = \frac{n(n+1)(2n+1)}{6}$; the last term of $\mathrm{(iii)}$ results from the easy-to-prove fact that
\begin{equation}
\sum \limits_{q=2}^E (q-1)\sum \limits_{k=1}^{q-1} y^k \leq \sum \limits_{q=1}^E \frac{(q+(E-1))(E-q)}{2} y^q,
\end{equation}
for $ y^k>0,k=1, \ldots,E-1$; the last term in $\mathrm{(iv)}$ is due to the fact that $(q+E-1)(E-q)\leq E(E-1)$.

Next, the upper bound of ${\sum}_{n=1}^N {\sum}_{q=1}^{E-1} \mathbb{E} [ \| \vb_n^{r, q}- \bar\vb^r \|^2 ]$ can be derived as
\begin{align}\label{eq: v rq 1}
& \sum_{n=1}^N \sum_{q=1}^{E-1} \mathbb{E} [ \| \vb_n^{r, q}- \bar\vb^r \|^2 ]\notag\\
&  \overset{\mathrm{(i)}}{=} \sum_{n=1}^N \sum_{q=1}^{E-1} \mathbb{E}\left[\left\|
    \sum_{m=1}^N W_{n,m} \vb_m^r-  \bar\vb^r
        +   (\gb_n^{r,q} -  \gb_n^r)
\right\|^2\right]\notag\\
&\overset{\mathrm{(ii)}}{\leq} 2 \sum_{n=1}^N \sum_{q=1}^{E-1} \! \mathbb{E} \bigg[ \bigg\|
    \sum_{m=1}^N W_{n,m} \vb_m^r-  \bar\vb^r \bigg\|^2 \bigg]
    + 2     \! \sum_{n=1}^N  \!\sum_{q=1}^{E-1} \! \mathbb{E}  [  \|
    \gb_n^{r,q} \!- \!\nabla f_n(\yb_n^{r, q})
    \! + \! \nabla f_n(\yb_n^{r, q})  \!- \! \nabla f_n(\bar\yb^r )
    \notag\\
           &~~~ ~~~
    + \nabla f_n(\bar\yb^r ) - \nabla f_n(\yb_n^{r })
    + \nabla f_n(\yb_n^{r })-  \gb_n^r
         \|^2 ]\notag\\
&\overset{\mathrm{(iii)}}{\leq} 2(E-1)\lambda_{w}^2 \phi_{v}^r + 16(E-1)N \frac{\sigma^2}{|\mathcal{I}|} + 8(E-1 ) L^2  \phi_{y}^r
        + 8   L^2 \sum_{n=1}^N \sum_{q=1}^{E-1} \mathbb{E} [ \|
        \yb_n^{r, q} - \bar\yb^{r }  \|^2 ],
\end{align}
where $\mathrm{(i)}$ is by \eqref{eq: v update r q}; $\mathrm{(ii)}$ is due to Jensen's inequality;  and $\mathrm{(iii)}$ is owing to Assumption \ref{assum: Lipschitz} and Assumption \ref{assum: stochastic}.

By inserting \eqref{eq: v rq 1} into \eqref{eq: y rq 1}, we have
\begin{align}\label{eq: y rq 2}
&\sum_{n=1}^N \sum_{q=1}^{E-1}\mathbb{E}\left[\|\yb_n^{r, q} - \bar\yb^r\|^2\right]\notag\\
& \leq \left[ 4(E-1)\lambda_{w}^2 + 16E(E-1)^2 L^2 \gamma^2 \right]\phi_{y}^r +4(E-1)\lambda_{w}^2[1+ E(E-1)  ] \gamma^2\phi_{v}^r\notag\\
    &~~~
    +  2E^2(E-1)N\gamma^2 \mathbb{E} [ \| \bar\vb^r \|^2 ]
    + 32E(E-1)^2 N  \gamma^2  \frac{\sigma^2}{|\mathcal{I}|} \notag\\
    &~~~
    +  16E(E\!-\!1)  L^2 \gamma^2 \sum_{n=1}^N \!\sum_{q=1}^{E\!-\!1} \! \mathbb{E}\!\left[\left\|
                 \yb_n^{r, q}\! -\! \bar\yb^{r } \!\right\|^2\right].
\end{align}

Let $\gamma \leq \textstyle \frac{1}{32 E    L  }$,
, i.e., $\frac{1}{1-  16E(E-1)   L^2\gamma^2 }\leq 2$. Then, letting $\gamma \leq \frac{\lambda_w^2}{2EL  }$ to simplify the first term in right hand side (RHS) of \eqref{eq: y rq 2}, and rearranging \eqref{eq: y rq 2}, we have
\begin{align}\label{eq: y r q final}
&\sum_{n=1}^N \sum_{q=1}^{E-1}\mathbb{E}\left[\|\yb_n^{r, q} - \bar\yb^r\|^2\right]\notag\\
& \leq  16(E-1)\lambda_{w}^2 \phi_{y}^r +8(E-1)\lambda_{w}^2[1+ E(E-1)  ] \gamma^2\phi_{v}^r\notag\\
    &
    +  4E^2(E-1)N\gamma^2 \mathbb{E}  [  \|  \bar \vb^r \|^2  ]
    + 64E(E-1)^2 N \gamma^2 \frac{\sigma^2}{|\mathcal{I}|} .
\end{align}
    %

\noindent $\bullet$ {\bf Proof of \eqref{eq: v r q}:}
By inserting \eqref{eq: y r q final} into \eqref{eq: v rq 1}, we obtain
\begin{align}
& \sum_{n=1}^N \sum_{q=1}^{E-1} \mathbb{E}\left[\left\| \vb_n^{r, q}- \bar\vb^r\right\|^2\right]\notag\\
& \leq 8(E-1) L^2 (1+ 16\lambda_{w}^2 )  \phi_{y}^r +\left\{ 2(E-1)\lambda_{w}^2 +  64(E-1)\lambda_{w}^2[1+ E(E-1)  ] \right. \left. L^2  N^2 \gamma^2 \right\} \phi_{v}^r
    \notag\\
    &~~~ +32   E^2(E-1) N  L^2\gamma^2  \mathbb{E} [ \|  \bar \vb ^r  \|^2 ]
    +16(E-1)N \left[ 1+  32E(E-1)   L^2\gamma^2  \right] \frac{\sigma^2}{|\mathcal{I}|} \label{eq: v r q final 1}\\
& \leq   8(E-1) L^2 (1+ 16 \lambda_{w}^2 )   \phi_{y}^r  +4(E-1)\lambda_{w}^2   \phi_{v}^r\! +\! 32   E^2(E\!-\!1) N  L^2\gamma^2  \mathbb{E}  [  \|  \bar \vb^r \|^2  ]
   \notag\\
    &~~~ \!+\!32(E\!-\!1)N   \frac{\sigma^2}{|\mathcal{I}|} ,\label{eq: v r q final}
\end{align}
where the $2$nd term and the $3$rd term of \eqref{eq: v r q final} reduces from  the $2$nd term and the $3$rd term of \eqref{eq: v r q final 1} by perspectively letting $\gamma  \leq
      \frac{ 1 }{  8EN L } $ and $\gamma  \leq
      \frac{ 1 }{  6E  L } $.

Thus, solving the intersection of the above conditions on $\gamma$
$
   \gamma\!   \leq \!
   \frac{  \lambda_{w}^2} {32 EN L  }$,
we can obtain the desired results in Lemma \ref{lem: y r q}.
 \hfill $\blacksquare$

%
\section{Proof of Lemma \ref{lem: y css}} \label{appen: y css}
By \eqref{eq: y update r+1} and \eqref{eq: y update ave}, we have
\begin{align}\label{eq: y css last 2}
&  \phi_{y}^{r+1} = \mathbb{E}\bigg[\bigg\|\Yb^{r+1} - \mathbf{1}(\bar\yb^{r+1})^\top\bigg\|_F^2\bigg]\notag\\
& =  \sum_{n=1}^N \mathbb{E}\bigg[\bigg\|  \sum_{m=1}^N W_{n,m} \yb_m^r - \bar\yb^r - \gamma   \bigg(\sum_{m=1}^N W_{n,m}\vb_m^r -\bar\vb^r\bigg)
 - \gamma    \frac{1}{N}\sum_{t=1}^N\sum_{q=1}^{E-1} (\vb_n^{r, q}- \vb_t^{r, q})\bigg\|^2\bigg]\notag\\
& \overset{\mathrm{(i)}}{\leq} (1+\delta_1)\lambda_{w}^2\phi_{y}^r
        + 2\bigg(1+ \frac{1}{\delta_1}\bigg)\lambda_{w}^2 \gamma^2\phi_{v}^r\notag\\
    & ~~~+ 2\bigg(1+ \frac{1}{\delta_1}\bigg) \gamma^2  \sum_{n=1}^N \frac{1}{N}\sum_{t=1}^N(E-1)\sum_{q=1}^{E-1}
            \mathbb{E}\bigg[\bigg\| \vb_n^{r, q}-\bar\vb^r
               +\bar\vb^r - \vb_t^{r, q}\bigg\|^2\bigg] \notag\\
&\overset{\mathrm{(ii)}}{\leq} (1+\delta_1)\lambda_{w}^2\phi_{y}^r
        + 2\bigg(1+ \frac{1}{\delta_1}\bigg) \lambda_{w}^2\gamma^2\phi_{v}^r
        + 4\bigg(1+ \frac{1}{\delta_1}\bigg)(E-1) \gamma^2   \sum_{n=1}^N \sum_{q=1}^{E-1}
            \mathbb{E} [ \| \vb_n^{r, q}-\bar\vb^r  \|^2 ] \notag\\
&\leq \bigg\{(1\!+\!\delta_1)\lambda_{w}^2
  \!  + \! 32\bigg(1\!+ \!\frac{1}{\delta_1}\bigg)L^2
     [ (  1\!+ \!16 \lambda_{w}^2)(E\!-\!1)^2]\gamma^2\bigg\}\phi_{y}^r  + 2\bigg(1+ \frac{1}{\delta_1}\bigg)\lambda_{w}^2[1+8(E-1)^2] \gamma^2   \phi_{v}^r\notag\\
    &~~~ \!+ \!256 \bigg(1 \!+ \! \frac{1}{\delta_1}\bigg)E^2(E \!- \!1)^2 N  L^2 \gamma^4 \mathbb{E} \! \bigg[\bigg \|  \frac{1}{N} \sum_{n=1}^N \!\nabla f_n(\yb_n^r)\bigg\|^2 \bigg] \notag\\
    &~~~  \!+  \!128\bigg(1+ \frac{1}{\delta_1}\bigg)(E \!- \!1)^2 \gamma^2( N \!+ \!2E^2L^2\gamma^2) \frac{\sigma^2}{|\mathcal{I}|},
\end{align}
where $\mathrm{(i)}$ is based on Young's inequality with $\delta_1>0$; $\mathrm{(ii)}$ is by Jensen's inequality; and the last inequality is by \eqref{eq: v r q final} in Lemma \ref{lem: y r q} and \eqref{eq: diff bar v true gr}.

Choose $\delta_1 = \textstyle \frac{1-\lambda_{w}^2}{8\lambda_{w}^2}$, and let $\gamma$ satisfy
\begin{align}\label{condi: zeta y css}
    \gamma   \leq \frac{(1-\lambda_{w}^2) ^2}{32 L^2(1+7\lambda_{w}^2)[8(1+16\lambda_{w}^2)(E-1)^2  ]},
\end{align}
which makes the $1$st term in RHS of \eqref{eq: y css last 2} reduce to $\textstyle \frac{1+3\lambda_{w}^2}{4}\phi_y^r$, and $N + 2E^2L^2\gamma^2 \leq 2N$ in the last term of \eqref{eq: y css last 2}'s RHS.
Then \eqref{eq: y css last 2} reduces to the desired results.
\hfill $\blacksquare$

\section{Proof of Lemma \ref{lem: v css}} \label{appen: v css}
For ease of analysis, recall $ \gb_{n}^r \triangleq \frac{1}{|\mathcal{I}|} \sum_{\xi \in \mathcal{I}_n^{r}} g_n(\yb_n^{r}, \xi ) $,  and  define $ {\Gb}^r = [ {\gb}_1^r, \ldots,  {\gb}_N^r]^\top, \bar { {\gb}}^r = \frac{1}{N}\sum_{n=1}^N {\gb}_n^r$.
Then, one can write \eqref{eq: v update r+1} and \eqref{eq: v update ave} as the following compact form
\begin{align}
& \Vb^{r+1} = \Wb\Vb^r +  {\Gb}^{r+1} -   {\Gb}^r ,\label{eq: v update r+1 compact}\\
& \bar\vb^{r+1} = \bar\vb^r +   \bar { {\gb}}^{r+1}-   \bar {{\gb}}^r .\label{eq: v update ave compact}
\end{align}
By \eqref{eq: v update r+1 compact} and \eqref{eq: v update ave compact}, we have
\begin{align}\label{eq: v css 1}
&  \phi_{v}^{r+1} = \mathbb{E}\left[\left\|\Vb^{r+1} - \mathbf{1}(\bar\vb^{r+1})^\top\right\|_F^2\right]\notag\\
& = \mathbb{E}\left[\left\|  \Wb\Vb^r- \mathbf{1}(\bar\vb^r)^\top +  \bigg(\mathbf{I} - \frac{1}{N} \mathbf{1}\mathbf{1}^\top\bigg)( {\Gb}^{r+1} -  {\Gb}^r)\right\|_F^2\right]\notag\\
& \overset{\mathrm{(i)}}{\leq} (1+\delta_2)\lambda_{w}^2 \phi_{v}^r
    +  \bigg(1+\frac{1}{\delta_2}\bigg)  \mathbb{E}\left[\left\| {\Gb}^{r+1} -  {\Gb}^r\right\|_F^2\right]\notag\\
& = (1+\delta_2)\lambda_{w}^2 \phi_{v}^r
    + \bigg(1+\frac{1}{\delta_2}\bigg) \sum_{n=1}^N  \mathbb{E} [ \|   \gb_n^{r+1} - \gb_n^r
     \|^2 ]\notag\\
& =  (1+\delta_2)\lambda_{w}^2 \phi_{v}^r
    + \bigg(1+\frac{1}{\delta_2}\bigg) \sum_{n=1}^N  \mathbb{E} [ \|
     \gb_n^{r+1}  - \nabla f_n(\yb_n^{r+1}) \notag\\
        &~~~~~~
      + \nabla_{\yb_n}f_n(\yb_n^{r+1}) - \nabla_{\yb_n}f_n(\yb_n^r)
      + \nabla_{\yb_n}f_n(\yb_n^r)- \gb_n^r
    \|^2 ]\notag\\
&\overset{\mathrm{(ii)}}{\leq}(1+\delta_2)\lambda_{w}^2 \phi_{v}^r
        + 6\bigg(1+\frac{1}{\delta_2}\bigg)N   \frac{\sigma^2}{|\mathcal{I}|}
         + 3\bigg(1+\frac{1}{\delta_2}\bigg)  L^2 \sum_{n=1}^N  \mathbb{E} [ \| \yb_n^{r+1} - \yb_n^r  \|^2 ],
\end{align}
where $\mathrm{(i)}$ is by Young's inequality with $\delta_2>0$ and the fact that $\textstyle \|\mathbf{I} - \frac{1}{N} \mathbf{1}\mathbf{1}^\top\|=1$; $\mathrm{(ii)}$ is due to Assumption \ref{assum: stochastic}.

Next, we solve the upper bound of $\textstyle \sum_{n=1}^N  \mathbb{E}[\| \yb_n^{r+1} - \yb_n^r \|^2]$. By \eqref{eq: y update r+1}, we have
\begin{align}
& \sum_{n=1}^N  \mathbb{E} [ \| \yb_n^{r+1} - \yb_n^r  \|^2 ]\notag\\
& \leq \sum_{n=1}^N \mathbb{E}\bigg[\bigg\|\sum_{m=1}^N W_{n,m} \yb_m^r - \bar \yb^r - (\yb_n^r -\bar \yb^r)
   \! - \!\gamma   \bigg(\sum_{m=1}^N \!W_{n,m} \vb_m^r\!-  \! \bar\vb^r\bigg) \! - \!\gamma   \sum_{q=1}^{E\!-\!1}(\vb_n^{r, q}\!-\!\bar\vb^r) \!-\! E \gamma   \bar\vb^r\bigg\|^2\bigg]\notag\\
& \overset{\mathrm{(i)}}{\leq} 5(1+\lambda_{w}^2) \phi_{y}^r + 5 \lambda_{w}^2 \gamma^2\phi_{v}^r
+ 5E^2N \gamma^2 \mathbb{E} [ \| \bar\vb^r \|^2 ]
        + 5(E-1) \gamma^2\sum_{n=1}^N \sum_{q=1}^{E-1}\mathbb{E} [ \| \vb_n^{r, q}-\bar\vb^r \|^2 ]  \notag\\
& \overset{\mathrm{(ii)}}{\leq}  [ 5(1+\lambda_{w}^2)\!  +\!40 (E-1)^2(1\!+ \!16 \lambda_{w}^2)  L^2 \gamma^2       ]\phi_{y}^r
 + 5 \lambda_{w}^2\gamma^2 [ 1+ 4(E-1)^2]\phi_{v}^r \notag\\
    &~~~ +5E^2N \gamma^2[1 + 32 (E-1)^2  L^2\gamma^2    ]\mathbb{E} [ \|  \bar \vb ^r  \|^2 ]
    +160(E-1)^2 N  \gamma^2   \frac{\sigma^2}{|\mathcal{I}|} \label{eq: y r+1- y r pro} \\
&\overset{\mathrm{(iii)}}{\leq}    10(1+\lambda_{w}^2) \phi_{y}^r   + 5  \lambda_{w}^2[1 +4(E-1)^2 ] \gamma^2\phi_{v}^r +10E^2N \gamma^2  \mathbb{E} [ \|  \bar  \vb^r  \|^2 ]
    +160(E-1)^2 N  \gamma^2   \frac{\sigma^2}{|\mathcal{I}|} ,\label{eq: y r+1- y r}
\end{align}
where $\mathrm{(i)}$ is by Jensen's inequality; $\mathrm{(ii)}$ is obtained by applying \eqref{eq: v r q} in Lemma \ref{lem: y r q}; and $\mathrm{(iii)}$ is by letting
    $\gamma \leq \frac{1}{136 E  L  }$,
so that the $1$st and the $3$rd terms of \eqref{eq: y r+1- y r} reduce from $1$st and the $3$rd terms of \eqref{eq: y r+1- y r pro}.

By inserting \eqref{eq: y r+1- y r} and \eqref{eq: diff bar v true gr} into \eqref{eq: v css 1}, we obtain
\begin{align}\label{eq: c css last two}
   \phi_{v}^{r+1}
& \leq \! 30\bigg(1+\frac{1}{\delta_2}\bigg)(1+\lambda_{w}^2) L^2 \phi_{y}^r
  \!+\! \bigg[\!(1+\delta_2) \!\lambda_{w}^2
            \!+\! 15 \! \bigg(1+\frac{1}{\delta_2}\bigg)\!\lambda_{w}^2[ 1+4(E-1)^2 ]\! L^2\gamma^2   \bigg]\!\phi_{v}^r\notag\\
    &~~~\! + 60\bigg(1+\frac{1}{\delta_2}\bigg)E^2N L^2\gamma^2
    \mathbb{E}\bigg[\bigg\|  \frac{1}{N}\sum_{n=1}^N \nabla f(\yb_n^r)\bigg\|^2\bigg]\notag\\
    &~~~\!+\! 6\bigg(1\!+\!\frac{1}{\delta_2}\bigg)\!N  \bigg\{1\! + \! 80(E\!-\!1)^2  L^2\gamma^2\! +\! \frac{5E^2L^2\gamma^2 }{N}\bigg\} \!  \frac{\sigma^2}{|\mathcal{I}|}.
\end{align}
By choosing $\delta_2 =\textstyle \frac{1-\lambda_{w}^2}{8\lambda_{w}^2}$, and
letting
    $\gamma\leq  \frac{(1-\lambda_{w}^2)}{3600EN L }$,
so that the $2$nd and $4$th terms of \eqref{eq: c css last two} reduces to those of \eqref{eq: v css}.
\hfill $\blacksquare$
\section{Proof of Lemma \ref{lem: descent lemma for F HL}} \label{appen: descent lemma for F HL}
Denote the gradient over the $ \bar \yb^r$ as $
 \nabla F(\bar \yb^{r}) \triangleq \frac{1}{N}\sum_{n=1}^N \nabla f_n(\bar\yb^r).$
Due to the Lipschitz smoothness in Assumption \ref{assum: Lipschitz}, we have
\begin{align}\label{eq: descent original}
  \mathbb{E}  [ F(\bar \yb^{r+1})]
 & \leq  \mathbb{E}  [ F(\bar \yb^{r})]
 +  \mathbb{E}  [ \langle \nabla F(\bar \yb^{r}) , \bar\yb^{r+1}\!-\!\bar\yb^r  \rangle ] + \frac{L}{2}\mathbb{E} [\|\bar\yb^{r+1} - \bar\yb^r\|^2 ].
\end{align}
Here, $\textstyle \mathbb{E}\left[  \langle \nabla F(\bar \yb^{r}), \bar\yb^{r+1}-\bar\yb^r \rangle\right]$ can be simplified as
\begin{align}\label{eq: descent 1}
&  \mathbb{E}\bigg[ \langle  \nabla F(\bar \yb^{r}), \bar\yb^{r+1}-\bar\yb^r \rangle\bigg] \notag\\
& \overset{\mathrm{(i)}}{=} \mathbb{E}\bigg[\bigg \langle
  \nabla F(\bar \yb^{r}), -\gamma  \bar \vb^r -\gamma   \frac{1}{N} \sum_{n=1}^N \sum_{q=1}^{E-1}\vb_n^{r, q}\bigg \rangle\bigg] \notag\\
& \overset{\mathrm{(ii)}}{=} -\gamma \mathbb{E}\bigg[\bigg \langle
   \nabla F(\bar \yb^{r}) ,   \bar \vb^r +  \frac{1}{N} \sum_{n=1}^N \sum_{q=1}^{E-1}
         \sum_{m=1}^N W_{n,m} \vb_m^r \bigg \rangle\bigg] -\gamma \mathbb{E} \bigg[  \bigg\langle
       \nabla F(\bar \yb^{r}) ,  \frac{1}{N}  \sum_{n=1}^N \sum_{q=1}^{E-1}\gb_n^{r, q} \bigg\rangle \bigg]\notag\\
    &~~~  +\gamma (E-1) \mathbb{E} \bigg[  \bigg\langle
        \nabla F(\bar \yb^{r}) ,
          \frac{1}{N}\sum_{n=1}^N  \gb_n^{r}
   \bigg\rangle \bigg]\notag\\
& \overset{\mathrm{(iii)}}{=} -\gamma E  \mathbb{E}\bigg [ \bigg \langle
   \nabla F(\bar \yb^{r}),   \bar \vb^r  \bigg\rangle \bigg]  -\gamma \mathbb{E} \bigg[  \bigg\langle
       \nabla F(\bar \yb^{r}) ,   \frac{1}{N} \sum_{n=1}^N \sum_{q=1}^{E-1}\nabla f_n(\yb_n^{r,q})  \bigg\rangle \bigg]\notag\\
    &~~~ +\gamma (E-1) \mathbb{E}\bigg[  \bigg\langle
        \nabla F(\bar \yb^{r}) , \frac{1}{N} \sum_{n=1}^N  \nabla f_n(\yb_n^{r })
         \bigg\rangle \bigg]\notag\\
& \overset{\mathrm{(iv)}}{=}-\gamma     \mathbb{E}\bigg[\bigg  \langle
  \nabla F(\bar \yb^{r}) ,\frac{1}{N}\sum_{n=1}^N  \nabla f_n(\yb_n^{r })  \bigg\rangle \bigg]  -\gamma   \sum_{q=1}^{E-1}
            \mathbb{E}\bigg[ \bigg\langle
               \nabla F(\bar \yb^{r}) ,\frac{1}{N}\sum_{n=1}^N\nabla f_n(\yb_n^{r,q})  \bigg \rangle \bigg] \notag\\
& \overset{\mathrm{(v)}}{=} -\gamma
    \bigg\{ \frac{1}{2}\mathbb{E}  [\|\nabla F(\bar \yb^{r})\|^2] + \frac{1}{2}\mathbb{E}\bigg[\bigg\|\frac{1}{N} \sum_{n=1}^N  \nabla f_n(\yb_n^{r }) \bigg\|^2\bigg]
     -\frac{1}{2} \mathbb{E} \bigg[\bigg\|\nabla F(\bar \yb^{r})-  \frac{1}{N}\sum_{n=1}^N  \nabla f_n(\yb_n^{r })\bigg\|^2\bigg]\bigg\}\notag\\
    &~~~ -\gamma   \bigg\{ \frac{E-1}{2}\mathbb{E}  [\|\nabla F(\bar \yb^{r})\|^2 ]
     + \frac{1}{2}\sum_{q=1}^{E-1}\mathbb{E} \bigg[\bigg\|\frac{1}{N} \sum_{n=1}^N  \nabla f_n(\yb_n^{r,q }) \bigg\|^2\bigg]
    \notag\\
    &~~~ ~~~ - \frac{1}{2}\sum_{q=1}^{E-1} \mathbb{E}\bigg[\bigg\|\nabla F(\bar \yb^{r}) - \frac{1}{N} \sum_{n=1}^N  \nabla f_n(\yb_n^{r ,q})\bigg\|^2\bigg]\bigg\}\notag\\
& \leq  -\frac{ \gamma  E}{2}
                  \|\nabla F(\bar \yb^{r})\|^2
        -\frac{ \gamma }{2} \mathbb{E} \bigg[ \bigg\|\frac{1}{N} \sum_{n=1}^N  \nabla f_n(\yb_n^{r }) \bigg\|^2\bigg] + \frac{ \gamma L^2  }{2N}   \phi_{y}^r
         +\frac{ \gamma     L^2}{2N} \sum_{n=1}^N \sum_{q=1}^{E-1}
       \mathbb{E}  [ \|\yb_n^{r ,q} - \bar\yb^{r }\|^2] ,
\end{align}
where $\mathrm{(i)}$ is by \eqref{eq: y update ave};  $\mathrm{(ii)}$ is due to \eqref{eq: v update r q}; the first term of $\mathrm{(iii)}$ is by the doubly stochastic property of $\Wb$; the second and the third terms of $\mathrm{(iii)}$ are obtained from the unbiased gradient in Assumption \ref{assum: stochastic}; $\mathrm{(iv)}$ is owing to \eqref{eq: unbiased diff bar v true gr}; $\mathrm{(v)}$ is by the fact that $\langle \ab, \bb\rangle = \textstyle \frac{1}{2}(\|\ab\|^2+ \|\bb\|^2-\|\ab-\bb\|^2)$ for some vectors $\ab, \bb$; and the last inequality is by omitting the negative term $\textstyle- \frac{ \gamma}{2}\sum_{q=1}^{E-1}\left\| \sum_{n=1}^N  \nabla f_n(\yb_n^{r,q }) \right\|^2 $.

Next, we solve the upper bound of $ \textstyle\frac{L}{2}\mathbb{E}\left[\|\bar\yb^{r+1} - \bar\yb^r\|^2\right]$ in \eqref{eq: descent original}. We have from \eqref{eq: y update ave} that,
\begin{align}\label{eq: decsent Last term}
& \frac{L}{2}\mathbb{E}\left[\|\bar\yb^{r+1} - \bar\yb^r\|^2\right]\notag\\
&=\frac{L \gamma^2}{2} \mathbb{E}\left[\left\|  \bar \vb^r +  \frac{1}{N} \sum_{n=1}^N \sum_{q=1}^{E-1}\vb_n^{r, q}\right\|^2\right]\notag\\
& \overset{\mathrm{(i)}}{=} \frac{L \gamma^2}{2} \mathbb{E}\left[\left\| E \bar \vb^r  +  \frac{1}{N} \sum_{n=1}^N \sum_{q=1}^{E-1} (\gb_n^{r, q}-\gb_n^r) \right\|^2\right]\notag\\
&  \overset{\mathrm{(ii)}}{\leq}    L\gamma^2
        \mathbb{E}\left[\left\| E\bigg( \bar \gb^r - \frac{1}{N} \sum_{n=1}^N \nabla f(\yb_n^r)\bigg)+ \frac{E}{N} \sum_{n=1}^N \nabla f(\yb_n^r)\right\|^2\right]\notag\\
&~~~
    +  3 L   \gamma^2   \mathbb{E}\bigg[\bigg\|  \frac{1}{N}\sum_{n=1}^N \sum_{q=1}^{E\!-\!1}  [\gb_n^{r, q}\!-\! \nabla f_n(\yb_n^{r,q})  ]  \bigg\|^2\bigg]
    + 3 L   \gamma^2   \mathbb{E}\bigg[\bigg\| \frac{1}{N} \sum_{n=1}^N \sum_{q=1}^{E-1} [  \gb_n^r \!-\! \nabla f_n(\yb_n^{r})]   \bigg\|^2\bigg]\notag\\
    &~~~
    + 3 L   \gamma^2   \mathbb{E}\bigg[\bigg\|  \frac{1}{N}\sum_{n=1}^N \sum_{q=1}^{E-1} [\nabla f_n(\yb_n^{r,q}) -\nabla f_n(\yb_n^{r })] \bigg\|^2\bigg]\notag\\
&  \overset{\mathrm{(iii)}}{\leq} 2EL\gamma^2
        \mathbb{E}\left[\left\|    \bar \gb^r - \frac{1}{N} \sum_{n=1}^N \nabla f(\yb_n^r)\right\|^2\right]
         + 2E^2L\gamma^2  \mathbb{E}\left[\left\|\frac{1}{N}\sum_{n=1}^N  \nabla f(\yb_n^r)\right\|^2\right]
    \notag\\
&~~~+  3 L   \gamma^2  \frac{1}{N^2}   \sum_{n=1}^N \sum_{q=1}^{E\!-\!1}\mathbb{E} [\| \gb_n^{r, q}\!-\! \nabla f_n(\yb_n^{r,q})      \|^2]
           +  3 L   \gamma^2   \frac{1}{N^2} \sum_{n=1}^N \sum_{q=1}^{E\!-\!1}\mathbb{E} [\| \gb_n^r \!-\! \nabla f_n(\yb_n^{r})\|^2 ]\notag\\
    &~~~
    + 3 (E-1)NL^3   \gamma^2  \frac{1}{N }  \sum_{n=1}^N \sum_{q=1}^{E-1}\mathbb{E} [ \|  \yb_n^{r,q}  - \yb_n^{r }  \|^2 ]\notag\\
&  \overset{\mathrm{(iv)}}{\leq} 2E^2  L\gamma^2
         \mathbb{E}\left[\left\|\frac{1}{N}\sum_{n=1}^N  \nabla f(\yb_n^r)\right\|^2\right]
        +  2(4E-3) L   \gamma^2    \frac{ \sigma^2}{N |\mathcal{I}|}
       \notag\\
    &~~~ + \frac{6 (E-1)^2  L^3   \gamma^2}{N}   \bigg[\phi_y^r +\sum_{n=1}^N \sum_{q=1}^{E-1}\mathbb{E} [ \|  \yb_n^{r,q}  - \bar\yb^{r }  \|^2 ]\bigg],
\end{align}%
where $\mathrm{(i)}$ is by  \eqref{eq: v update r q}; $\mathrm{(ii)}$ is due to Jensen's inequality; the second and the third terms of $\mathrm{(iii)}$ are obtained by the fact $\textstyle\|\sum_{n=1}^N (\xb_n - \mathbb{E}[\xb_n]\|^2  = \sum_{n=1}^N \|\xb_n - \mathbb{E}[\xb_n])\|^2, \forall \xb,$ due to Assumption \ref{assum: stochastic}; the last term of $\mathrm{(iii)}$ is owing to Assumption \ref{assum: Lipschitz}; $\mathrm{(iv)}$ is by the bounded variance in Assumption \ref{assum: stochastic}.

Summing up \eqref{eq: descent 1} and \eqref{eq: decsent Last term},  we obtain
\begin{align}
    &  \mathbb{E}\left[  \langle  \bar \nabla_{\yb_n}^r, \bar\yb^{r+1}-\bar\yb^r  \rangle\right]
    +   \frac{L}{2}\mathbb{E}\left[\|\bar\yb^{r+1} - \bar\yb^r\|^2\right]\notag\\
& \leq  -\frac{ \gamma  E}{2}
                  \|\bar \nabla_{\yb}^r\|^2
        -\frac{ \gamma }{2}[1- 4E^2  L\gamma  ] \mathbb{E} \bigg[ \bigg\| \frac{1}{N}\sum_{n=1}^N  \nabla f_n(\yb_n^{r }) \bigg\|^2\bigg]
        + \frac{ \gamma L^2  }{2N}   [1+ 12(E-1)^2 L\gamma  ]\phi_{y}^r
      \notag\\
    &~~~   +\frac{ \gamma     L^2}{2N} [1+ 12(E-1)^2 L\gamma ]
    \sum_{n=1}^N \sum_{q=1}^{E-1}
       \mathbb{E}  [ \|\yb_n^{r ,q} - \bar\yb^{r }\|^2]
     +  2 (4E-3) L   \gamma^2   \frac{ \sigma^2}{N |\mathcal{I}|}\label{descent sum 2 pro 1} \\
& \overset{\mathrm{(i)}}{\leq }  -\frac{ \gamma  E}{2}
                  \|\bar \nabla_{\yb}^r\|^2
        -\frac{ \gamma }{2}[1- 4E^2  L\gamma  ] \mathbb{E} \bigg[ \bigg\|  \frac{1}{N}\sum_{n=1}^N  \nabla f_n(\yb_n^{r }) \bigg\|^2\bigg] \notag\\
    &~~~
        +    \frac{  L^2 \gamma }{N}  \phi_{y}^r
         + \frac{  L^2 \gamma }{N} \sum_{n=1}^N \sum_{q=1}^{E-1} \mathbb{E}  [ \|\yb_n^{r ,q} - \bar\yb^{r }\|^2]
    +  2 (4E-3)  L   \gamma^2    \frac{ \sigma^2}{N |\mathcal{I}|} \label{descent sum 2 pro 2}\\
& \overset{\mathrm{(ii)}}{\leq }  -\frac{ \gamma  E}{2}
                  \|\bar \nabla_{\yb}^r\|^2
        -\frac{ \gamma }{2}\bigg[1- 4E^2  L\gamma
            - 16 E^2(E-1) L^2 \gamma^2
        \bigg] \mathbb{E} \bigg[ \bigg\| \frac{1}{N}\sum_{n=1}^N  \nabla f_n(\yb_n^{r }) \bigg\|^2\bigg] \notag\\
    &~~~+   \frac{L^2 \gamma}{N}  [1+  16(E-1)\lambda_{w}^2] \phi_{y}^r + 8(E-1)\lambda_{w}^2[1+ E(E-1)  ]  \frac{  L^2 \gamma^3}{N}   \phi_{v}^r\notag\\
    &~~~+\frac{2   L   \gamma^2}{N}     \bigg[ 4E-3 + 32E(E-1)^2   N  L \gamma   \bigg]\frac{\sigma^2}{|\mathcal{I}|} \label{descent sum 2 pro 3}\\
& \overset{\mathrm{(iii)}}{\leq } -\frac{ \gamma  E}{2}
                  \|\bar \nabla_{\yb}^r\|^2
        -\frac{ \gamma }{2}\bigg[1- 2E^2  L\gamma
             - 8 E^2(E-1) L^2 \gamma^2
        \bigg] \mathbb{E} \bigg[ \bigg\| \frac{1}{N}\sum_{n=1}^N  \nabla f_n(\yb_n^{r }) \bigg\|^2\bigg] \notag\\
    &~~~+   \frac{L^2 \gamma}{N}  [1+  16(E-1)\lambda_{w}^2] \phi_{y}^r
    +10EL   \gamma^2    \frac{\sigma^2}{N |\mathcal{I}|} + 8(E-1)\lambda_{w}^2 [1+ E(E-1) ]  \frac{ L^2 \gamma^3}{N}   \phi_{v}^r,\label{descent sum 2}
\end{align}
where $\mathrm{(i)}$ and $\mathrm{(iii)}$ are by letting
\begin{align}
\gamma \leq \min \bigg\{\frac{1}{ 12E^2 L   }, \frac{1}{32E^2  N  L     }\bigg\}= \frac{1}{32E^2  N  L     },
\end{align}
so that the $3$rd term of \eqref{descent sum 2 pro 1} reduces to the $3$rd term of \eqref{descent sum 2 pro 2} and the last term of \eqref{descent sum 2 pro 3} reduces to the last term of \eqref{descent sum 2};
and  $\mathrm{(ii)}$ is due to \eqref{eq: v r q} in Lemma \ref{lem: y r q}.

Thus, after combining \eqref{eq: descent original} and \eqref{descent sum 2}, we have desired results. \hfill $\blacksquare$
\section{Proof of the fact $\rho(\Ab)< 1$}\label{appen:  phi iterates FL}
Assume $\bsb=[s_1, s_2]^\top \in \mathbb{R}^{2}$. One can verify $\Ab\bsb <\bsb$ by equivalently solving
\begin{align}
   & \begin{cases}
   \Ab_{1, 1}s_1 + \Ab_{1, 2}s_2 < s_1, \\
     \Ab_{2, 1}s_1 + \Ab_{2, 2}s_2 < s_2,
    \end{cases} \notag\\
\Leftarrow
    & \begin{cases}
   \Ab_{1, 2}s_2 < (1-\Ab_{1, 1})s_1,   \\
      \Ab_{2, 1}s_1 < (1-\Ab_{2, 2})s_2  ,
    \end{cases} \notag\\
\Leftarrow
    & \begin{cases}
    2  \frac{(1+7\lambda_{w}^2)  }{1-\lambda_{w}^2} \left[1 + 8(E-1)^2   \lambda_{w}^2  \right] s_2\gamma^2  <  \frac{3(1-\lambda_{w}^2) }{4} s_1 ,  \\
      30\frac{(1+7\lambda_{w}^2)(1+\lambda_{w}^2)}{1-\lambda_{w}^2}N^2 L^2 s_1 <\frac{3(1-\lambda_{w}^2) }{4}s_2,
    \end{cases} \notag\\
\Leftarrow
    & \begin{cases}
   \gamma^2 \leq       \frac{3(1 \!- \!\lambda_{w}^2)^2 s_1}{16(1 \!+ \!7\lambda_{w}^2) [1  \!+ \! 4(E \!- \!1)^2   \lambda_{w}^2  ]  s_2   } \!<   \!  \frac{3(1 \!- \!\lambda_{w}^2)^2 s_1}{8(1 \!+ \!7\lambda_{w}^2) [1  \!+ \! 8(E \!- \!1)^2   \lambda_{w}^2   ]  S_2   }  ,  \\
      s_1 <\frac{3(1-\lambda_{w}^2) }{120\frac{(1+7\lambda_{w}^2)(1+\lambda_{w}^2)}{1-\lambda_{w}^2}N^2 L^2 }s_2.
    \end{cases}
\end{align}%
Thus, for $\bsb$ and $\gamma$ satisfying
\begin{align*}
    s_1& <\frac{ (1-\lambda_{w}^2)^2 }{40 (1+7\lambda_{w}^2)(1+\lambda_{w}^2) L^2}s_2,\\
    \gamma & \leq       \frac{3(1-\lambda_{w}^2) s_1}{16(1+7\lambda_{w}^2) [1 + 8(E-1)   \lambda_{w}  ]  s_2   },
\end{align*}
we have $ \Ab\bsb <\bsb$. Then, one can obtain that $\rho(\Ab)<1$.
\hfill $\blacksquare$

\section{Proof of Lemma \ref{lem: det I -G}}\label{appen: det I -G}
Following the rule that $(\Ib - \Ab)^{-1} = (\Ib - \Ab)^*/\det (\Ib- \Ab)$, we divide the derivation into two steps, including solving the  determinant and the adjoint matrix $ (\Ib - \Ab)^*$,  respectively.

\noindent $\bullet$ {\bf Solve the determinant $\det (\Ib- \Ab)$:}

Recalling $\Ab$ in \eqref{eq: G def}, we have
\begin{align}\label{eq: det I-G}
 & \det (\mathbf{I} - \Ab) \notag\\
 & = (1-A_{1, 1})(1-A_{2,2})-A_{1, 2}A_{2,1}\notag\\
 & = \frac{9(1-\lambda_{w}^2)^2}{16}   - 60  \frac{(1+7\lambda_{w}^2)^2 (1+\lambda_{w}^2) }{(1-\lambda_{w}^2)^2} [1
    + 8(Q-1)^2   \lambda_{w}^2    ]  L^2\gamma^2  \notag\\
&  \geq  \frac{(1-\lambda_{w}^2)^2}{8},
\end{align}
where the last inequality is obtained by letting
$$\gamma  \leq \textstyle \frac{3(1-\lambda_{w}^2)^2}{320(1+7\lambda_{w}^2) (1+\lambda_{w}^2)[1
    + 8(Q-1)    \lambda_{w}  ]  L   }.$$

\noindent $\bullet$ {\bf Solve the determinant $ (\Ib- \Ab)^*$:}
\begin{align}\label{eq: adjoint I-G}
 &  (\mathbf{I} - \Ab)^*
  = \begin{bmatrix}
        1-G_{2, 2} & A_{1,2}\\
        A_{2,1} & 1-A_{1,1}
    \end{bmatrix}
  \leq  \begin{bmatrix}
        1 &  \frac{ 2(1+7\lambda_{w}^2)  [1
    + 8(Q-1)^2   \lambda_{w}^2   ] \gamma^2 }{1-\lambda_{w}^2}\\
\frac{ 30(1+7\lambda_{w}^2)(1+\lambda_{w}^2) L^2 }{1-\lambda_{w}^2}& 1
    \end{bmatrix},
\end{align}
where the last inequality is element-wise.

By combining \eqref{eq: det I-G} and \eqref{eq: adjoint I-G}, we have the desired results.
\hfill $\blacksquare$

%
\section{Alternative Expression of Algorithm \ref{alg: MUST}}\label{appen: alter HBL}
For ease of analysis, we can conclude the \eqref{eq: Alg2 q2} for $n\in [N]$ at $q$-th iteration of the $r$-th round as
\begin{subequations}\label{eq: r q update HBL}
\begin{align}
    & \thetab_n^{r, q}   =  \sum_{m=1}^N  \Wb_{n, m}\thetab_m^r
        - \alpha\gb_{\theta,n}^{r }
        - \alpha \sum_{k=1}^{q-1} \gb_{\theta,n}^{r,k} ,\label{eq: update theta r q HBL}\\
    & \xb_n^{r, q }= \sum_{m=1}^N  \Wb_{n, m}\xb_m^r
        - \beta  \sum_{m=1}^N  \Wb_{n, m} \ub_m^r
        - \beta  \sum_{k=1}^{q-1}  \ub_n^{r, k},\label{eq: update x r q HBL}\\
   & \zb_{n, i}^{r, q}\! = \!  \sum_{m=1}^N \! \Wb_{n, m} \zb_{m, i}^r \!+\! N\Bb_{n, i} \bigg( \xb_n^{r, q}
        \!-\!   \sum_{m=1}^N \!\Wb_{n, m} \xb_m^{r}\bigg),  \forall i\in [S], \label{eq: update z r q HBL}\\
      & \ub_n^{r, q} =   \sum_{m=1}^N  \Wb_{n, m} \ub_m^r
        + \gb_{x,n}^{r, q} - \gb_{x,n}^{r}.\label{eq: update u r q HBL}
\end{align}
\end{subequations}
By inserting $q=E$ into \eqref{eq: r q update HBL}, we have
\begin{subequations}\label{eq: full update HBL}
\begin{align}
    & \thetab_n^{r+1}   =  \sum_{m=1}^N  \Wb_{n, m}\thetab_m^r
        - \alpha\gb_{\theta,n}^{r }
        - \alpha \sum_{q=1}^{E-1} \gb_{\theta,n}^{r,q} ,\label{eq: full update theta HBL}\\
    & \xb_n^{r+1} =  \sum_{m=1}^N  \Wb_{n, m}\xb_m^r
        - \beta    \sum_{m=1}^N  \Wb_{n, m} \ub_m^r
        - \beta  \sum_{q=1}^{E-1}  \ub_n^{r, q},\label{eq: full update x HBL}\\
   & \zb_{n, i}^{r+1} \!= \!  \sum_{m=1}^N \! \Wb_{n, m} \zb_{m, i}^r \!+\! N\Bb_{n, i}    \bigg( \xb_n^{r, q}
        \!-\!    \sum_{m=1}^N \!\Wb_{n, m} \xb_m^{r} \bigg), \forall i\in [S], \label{eq: full update z HBL}\\
      & \ub_n^{r+1} =   \sum_{m=1}^N  \Wb_{n, m} \ub_m^r
        + \gb_{x,n}^{r+1} - \gb_{x,n}^{r}.\label{eq: full update u HBL}
\end{align}
\end{subequations}
By making an average on the above subequations, we obtain
\begin{subequations}\label{eq: bar update HBL}
\begin{align}
    & \bar \thetab^{r+1}   =   \bar \thetab^{r}
        - \alpha \bar \gb_{\theta}^{r }
        - \alpha \frac{1}{N}\sum_{n=1}^N \sum_{q=1}^{E-1} \gb_{\theta,n}^{r,q} ,\label{eq: theta bar update HBL}\\
    & \bar\xb^{r+1} = \bar\xb^{r}
        - \beta  \bar \ub^r
        - \beta \frac{1}{N}\sum_{n=1}^N  \sum_{q=1}^{E-1}  \ub_n^{r, q},\label{eq: x bar update HBL}\\
   & \bar \zb_{ i}^{r+1}  =   \bar \zb_{ i}^{r}  +
    \sum_{n=1}^N \Bb_{n, i}
       \bigg( \xb_n^{r+1}
        \!-\!   \sum_{m=1}^N \!\Wb_{n, m} \xb_m^{r}\bigg),\!  \forall i\in [S], \label{eq: z bar update HBL}\\
      & \bar \ub^{r+1} =   \bar \ub^{r}
        + \bar \gb_{x }^{r+1} -\bar \gb_{x}^{r},\label{eq: u bar update HBL}
\end{align}
\end{subequations}
where $\bar \theta^r \triangleq \textstyle \frac{1}{N} \sum_{n=1}^N \thetab_n^{r }$, and $\bar \xb^r$, $\bar \zb^r$, $\bar \ub^r$,   $ \bar \gb_{\theta}^{r } $, $\bar \gb_{x}^{r } $ are defined in the same fashion.

Similar to Lemma \ref{lem: ave = true sum}, we have
\begin{Lemma} \label{lem: ave = true sum HBL}
  For all $ n\in [N]$, initializing $\zb_{n,i}^0 = N \Bb_{n, i} \xb_{n}^0, \forall i\in[S]$, $\ub_n^0 =  \gb_{x,n}^0$, by \eqref{eq: full update z HBL} and \eqref{eq: full update u HBL}, we have
  \begin{subequations}
    \begin{align}
        & \bar \zb_i^r =  \sum _{n=1}^N \Bb_{n,i}  \xb_n^r, \forall i\in[S],\label{eq: z ave = true sum HBL}\\
       &  \bar \ub^r=  \bar{\gb}_x^{r} , \forall r\geq 0.
        \label{eq: ave = true sum HBL}
    \end{align}
    \end{subequations}
\end{Lemma}

Then, similar to \eqref{eq: unbiased diff bar v true gr} and \eqref{eq: diff bar v true gr}, we have
\begin{align}
    & \mathbb{E}[\bar\ub^r] =  \mathbb{E}[\bar\gb_x^r]  =\frac{1}{N}   \sum_{n=1}^N \frac{1}{ S} \sum_{i=1}^S\Bb_{n, i}^\top \nabla_{\zb}f( \zb_{n, i}^r, \thetab_n^r) ,\label{eq: unbiased diff bar u true gr HBL} \\
    & \mathbb{E}[\|\bar\ub^r\|^2]  = \mathbb{E}[\|\bar\gb_x^r\|^2] \leq  \frac{2 \sigma }{N|\mathcal{I}|} + 2\mathbb{E}\bigg[\bigg\| \frac{1}{N}   \sum_{n=1}^N \frac{1}{ S} \sum_{i=1}^S\Bb_{n, i}^\top \nabla_{\zb}f( \zb_{n, i}^r, \thetab_n^r)\bigg\|^2\bigg].\label{eq: diff bar u true gr HBL}
\end{align}
To establish convergence, for $r$-th round, we
denote
$\Thetab^r= [\thetab_1^r , \cdots, \thetab_N^r ]^\top $
and define $\Xb^
r$, $\Zb^r$ and $\Ub^
r$ in the same fashion.
Based on these notations,  we define
\begin{align}
\hat{\phib}^{r} =\begin{bmatrix}
             \phi_{\theta}^{r}\\
             \phi_{x}^{r}\\
             \phi_{z}^{r}\\
             \phi_{w}^{r}
        \end{bmatrix} = \begin{bmatrix}
             \mathbb{E}\left[ \|\Thetab^{r} - \mathbf{1}(\bar{\thetab}^{r})^\top\|_F^2\right]\\
               \mathbb{E}\left[ \|\Xb^{r} - \mathbf{1}(\bar{\xb}^{r})^\top\|_F^2\right]\\
                 \mathbb{E}\left[ \|\Zb^{r} - \mathbf{1}(\bar{\zb}^{r})^\top\|_F^2\right]\\
                 \mathbb{E}\left[   \|\Ub^{r} - \mathbf{1}(\bar{\ub}^{r})^\top\|_F^2\right]
        \end{bmatrix}\in \mathbb{R}^4
\end{align}
as the consensus error matrix, where $\|\cdot\|_F$ denotes the matrix Frobenius norm. Moreover, the gradients over average $\bar \thetab^r$ and $\bar \xb^r$ are defined as
\begin{align}
& \bar \nabla_{\thetab}^r \triangleq  \nabla_{\thetab} \tilde{F}(\bar \thetab^r, \bar {\xb}^r)=  \frac{1}{NS}\sum_{i=1}^S\nabla_{\thetab} f\left( \Bb_{i}\bar\xb^r, \bar\thetab^r\right),\\
&\bar \nabla_{\thetab}^r \triangleq
 \nabla_{\xb} \tilde{F}(\bar \thetab^r, \bar {\xb}^r)  = \frac{1}{NS}\sum_{i=1}^S \Bb_{ i}^\top \nabla_{  \Bb_{i}\bar\xb^r} f\left( \Bb_{i}\bar\xb^r, \bar\thetab^r\right).
\end{align}

\section{Proof of Theorem \ref{thm: HBL}} \label{appen: thm HBL}
\subsection{Key Lemmas}
\begin{Lemma} \label{lem: local step consensus HBL}
    (Local updates) Let Assumption  \ref{assum network} to Assumption \ref{assum: Lipschitz}, and Assumption \ref{assum: stochastic HBL} hold. For sufficiently small $\alpha$ and $\beta $ satisfying
    \begin{align}
& \alpha \leq \min \bigg\{ \frac{\lambda_w}{20EB_{\max}L}, \frac{\lambda_w}{4EL}\bigg\},
 \\
& \beta \leq  \min \bigg\{  \frac{\lambda_w}{28ENB_{\max} L}, \frac{\lambda_w}{56ENB_{\max}^2L},  \frac{\lambda_w}{272ENB_{\max}^3L}  \bigg\},
\end{align}
    we have the following dynamics of distance from local-step variable $\thetab_n^{r, q}, \xb_n^{r, q}, \zb_n^{r, q}, \ub_n^{r, q}$ to the corresponding average $\bar\thetab^r, \bar \xb^r, \bar\zb^r, \bar \ub^r$ generated by our proposed MUST algorithm
    \begin{subequations}
    \begin{align}
         &  \sum_{n=1}^N \sum_{q=1}^{E-1}\frac{1}{S}\sum_{i=1}^S  \left\|\zb_{n, i}^{r, q} - \bar \zb_i^r\right\|^2  +  \sum_{n=1}^N \sum_{q=1}^{E-1} \left\|\thetab_n^{r, q} - \bar \thetab^r\right\|^2 \notag \\
         &  \leq   7(E-1)\lambda_w^2 \phi_\theta^r
             + 200(E-1) N^2B_{\max}^2\lambda_w^2\phi_x^r
             + 32(E-1)  \lambda_w^2\frac{1}{S} \phi_z^r \notag\\
            &~~~ +128(E-1) [1+ E(E-1) ]  N^2B_{\max}^2\lambda_w^2\beta^2 \phi_u^r
             + 2E(E-1)N[ 5\alpha^2
                + 768   N B_{\max}^4 L^2\beta^2
            ]\frac{\sigma^2}{|\mathcal{I}|} \notag\\
            & ~~~+ 32E^2(E-1)N\alpha^2\mathbb{E} [\|\bar \nabla_{\thetab}^r  \|^2 ]
             + 48E^2(E-1)N^2B_{\max}^2 \beta^2 \mathbb{E}[  \|\bar\ub^r  \|^2] , \label{eq: local step consensus z theta HBL}\\
          &  \sum_{n=1}^N \sum_{q=1}^{E-1}  \left\|\ub_{n}^{r, q} - \bar \ub ^r\right\|^2\notag\\
            &
            \leq  12 (E-1) (1+ 7 \lambda_w^2)B_{\max}^2 L^2 \phi_{\theta}^r
             + 2400  (E-1) N^2B_{\max}^4 L^2 \lambda_w^2\phi_x^r\notag\\
    &~~~+ 12(E-1) (1+ 32 \lambda_w^2) B_{\max}^2 L^2  \frac{1}{S} \phi_z^r
    +3(E-1)\lambda_w^2 \phi_u^r\notag\\
    &~~~+24 E(E-1)N B_{\max}^2 L^2  [ 5\alpha^2
            + 768   N B_{\max}^4 L^2\beta^2
        ]\frac{\sigma^2}{|\mathcal{I}|}
        + 384 E^2(E-1)NB_{\max}^2 L^2 \alpha^2\mathbb{E} [\|\bar \nabla_{\thetab}^r  \|^2 ]\notag\\
    &~~~  +  576E^2(E-1)N^2B_{\max}^4L^2 \beta^2 \mathbb{E}[  \|\bar\ub^r  \|^2],
        \label{eq: local step consensus w HBL}
        \end{align}
        and
        \begin{align}
            & \sum_{n=1}^N \sum_{q=1}^{E-1}  \left\|\xb_{n}^{r, q} - \bar \xb^r \right\|^2\notag\\
            &  \leq   24 E(E\!-\!1)^2(1\!+\!7\lambda_w^2) B_{\max}^2 L^2\beta^2
        \phi_{\theta}^r
        \! +\! 5 (E\!-\!1)\lambda_w^2   \phi_x^r
         \!+\!24 E(E\!\!-1)^2 (1\!+\!32\lambda_w^2)B_{\max}^2 L^2 \frac{1}{S}\beta^2  \phi_z^r\notag\\
    &~~~ + 5(E-1)\beta^2 [1+ E(E-1) ] \lambda_w^2  \phi_u^r
     + 25 E(E-1)^2  N \beta^2 \frac{\sigma^2}{|\mathcal{I}|} \notag\\
    &~~~
    + 3 E^2(E-1)\beta^2 \mathbb{E}[  \|\bar\ub^r  \|^2]
    +768 E^3(E-1)^2N B_{\max}^2 L^2\alpha^2\beta^2\mathbb{E} [\|\bar \nabla_{\thetab}^r  \|^2 ].
         \label{eq: local step consensus x HBL}
    \end{align}
    \end{subequations}
\end{Lemma}

\begin{Lemma} \label{lem: consensus HBL}
    (Consensus error) Let Assumption  \ref{assum network} to Assumption \ref{assum: Lipschitz}, and Assumption \ref{assum: stochastic HBL} hold. For sufficiently small $\alpha$ and $\beta$ satisfying
    \begin{align}
& \alpha \leq \frac{1}{144EL}, ~~\beta \leq  \frac{1}{109 ENB_{\max}^2 L} ,
\end{align}
    we have the following contraction property of the iterates generated by MUST algorithm
    \begin{subequations}
    \begin{align}
         & \phi_\theta^{r+1} = \mathbb{E}\left[ \left\| \Thetab^{r+1} - \mathbf{1}(\bar\thetab^{r+1})^\top \right\|^2 \right] \notag\\
         & \leq\frac{1+3\lambda_w^2}{4}  \phi_\theta^r
         + \!80  \frac{1+7\lambda_w^2}{1-\lambda_w^2} NB_{\max}^2 L^2[ E^2 \!+\!  100  (E\!-\!1) \lambda_w^2N] \alpha^2
             \phi_x^r\notag\\
            &~~~ + 40\frac{1+7\lambda_w^2}{1-\lambda_w^2}[1+32(E-1)\lambda_w^2 ]  \frac{1}{S}  L^2 \alpha^2 \phi_z^r
            +5120 \frac{1+7\lambda_w^2}{1-\lambda_w^2}(E-1)[1+E(E-1) ]
                \notag\\
                &~~~ ~~~N^2B_{\max}^2L^2\alpha^2\lambda_w^2\beta^2 \phi_u^r
                + 48 \frac{1+7\lambda_w^2}{1-\lambda_w^2}EN \alpha^2  \frac{\sigma^2}{|\mathcal{I}|}
                \notag\\
            &~~~+ 88 \frac{1+7\lambda_w^2}{1-\lambda_w^2}E^2\alpha^2  \mathbb{E} [\|\bar \nabla_{\thetab}^r  \|^2 ]\notag\\
            &~~~  + 1920 \frac{1+7\lambda_w^2}{1-\lambda_w^2} E^2(E-1)N^2B_{\max}^2 L^2\alpha^2\beta^2 \mathbb{E}[  \|\bar\ub^r  \|^2],\label{eq: theta consensus}
            \end{align}
            and
            \begin{align}
         & \phi_x^r = \mathbb{E}\left[ \left\|  \Xb^r - \mathbf{1}(\bar\xb^r)^\top\right\|^2\right]\notag\\
         & \leq   96\frac{1+7\lambda_w^2}{1-\lambda_w^2} (E-1)^2 (1+ 7
        \lambda_w^2)B_{\max}^2 L^2 \beta^2 \phi_{\theta}^r  + \frac{1+3\lambda_w^2}{4} \phi_x^r\notag\\
            &~~~ +  96\frac{1+7\lambda_w^2}{1-\lambda_w^2} (E-1)^2(1+ 32 \lambda_w^2) B_{\max}^2 L^2  \frac{1}{S} \beta^2
                   \phi_z^r
                   +  2\frac{1+7\lambda_w^2}{1-\lambda_w^2}\beta^2
            [1+ 12 (E-1)^2  \lambda_w^2]  \phi_u^r \notag\\
            &~~~ +  192\frac{1+7\lambda_w^2}{1-\lambda_w^2}   E(E-1)^2N B_{\max}^2 L^2 [ 5\alpha^2 + 768   N B_{\max}^4 L^2\beta^2 ] \beta^2\frac{\sigma^2}{|\mathcal{I}|}\notag\\
            &~~~ + 3072 \frac{1+7\lambda_w^2}{1-\lambda_w^2} E^2 (E-1)^2 NB_{\max}^2 L^2 \alpha^2\beta^2
        \mathbb{E} [\|\bar \nabla_{\thetab}^r  \|^2 ]\notag\\
        &~~~ +4608 \frac{1+7\lambda_w^2}{1-\lambda_w^2} E^2(E-1)^2N^2B_{\max}^4L^2 \beta^4 \mathbb{E}[  \|\bar\ub^r  \|^2].\label{eq: x consensus}
    \end{align}
    \end{subequations}
\end{Lemma}

\begin{Lemma} \label{lem: tracking HBL}
    (Tracking error) Let Assumption  \ref{assum network} to Assumption \ref{assum: Lipschitz}, and Assumption \ref{assum: stochastic HBL} hold. For sufficiently small $\alpha$ and $\beta$ satisfying
    \begin{align}
& \alpha \leq \frac{1}{78EL}, ~\beta \leq \min \bigg\{\frac{1}{712 ENB_{\max}^2 L}, \frac{1}{EN\lambda_w}\bigg\},
\end{align}
    we have the following contraction property of the iterates generated by MUST algorithm
    \begin{subequations}
    \begin{align}
         & \phi_z^{r+1} = \mathbb{E}\left[ \left\| \Zb^{r+1} - \mathbf{1}(\bar\zb^{r+1})^\top \right\|^2 \right] \notag\\
         & \leq  144 \frac{1+7\lambda_w^2}{1-\lambda_w^2}(E-1)^2(1+ 7 \lambda_w^2) N^2 S
            B_{\max}^4  L^2 \beta^2  \phi_{\theta}^r
            \notag\\
            &~~~+ 28800\frac{1+7\lambda_w^2}{1-\lambda_w^2} (E-1)^2N^4 S B_{\max}^6 L^2\lambda_w^2\beta^2\phi_x^r \notag\\
            &~~~ +  \frac{1+3\lambda_w^2}{4} \phi_z^r
            + 12\frac{1+7\lambda_w^2}{1-\lambda_w^2}[1 + 3 (E-1)^2  \lambda_w^2] N^2 S B_{\max}^2 \beta^2\phi_u^r\notag\\
            &~~~ + 288\frac{1+7\lambda_w^2}{1-\lambda_w^2}  E(E\!-\!1)^2N^3 SB_{\max}^4 L^2  [ 5\alpha^2
                               \!
                               +\! 768   N B_{\max}^4 L^2\beta^2
                            ]\beta^2\frac{\sigma^2}{|\mathcal{I}|}\notag\\
            &~~~+ 4608\frac{1+7\lambda_w^2}{1-\lambda_w^2}   E^2(E-1)^2N^3SB_{\max}^4 L^2 \alpha^2\beta^2
                    \mathbb{E} [\|\bar \nabla_{\thetab}^r  \|^2 ]
                    +  24\frac{1+7\lambda_w^2}{1-\lambda_w^2}
                     E^2N^3 S  B_{\max}^2 \beta^2  \mathbb{E}[  \|\bar\ub^r  \|^2].\label{eq: z consensus HBL}
                     \end{align}
                     and
                     \begin{align}
         & \phi_u^{r+1} = \mathbb{E}\left[ \left\|  \Ub^{r+1} - \mathbf{1}(\bar\ub^{r+1})^\top\right\|^2\right]\notag\\
         & \leq 49 \frac{1+7\lambda_w^2}{1-\lambda_w^2}(1+\lambda_w^2) B_{\max}^2L^2 \phi_\theta^r  \notag\\
            &~~~+ 12\frac{1+7\lambda_w^2}{1-\lambda_w^2}
                \bigg\{
                    21600(E-1)^2N^4B_{\max}^6L^2\lambda_w^2 \beta^2
                    + 30NB_{\max}^2L^2 [E^2 + 160 (E-1)\lambda_w^2N]\alpha^2
                \bigg\}  \phi_x^r\notag\\
            &~~~+ 49\frac{1+7\lambda_w^2}{1-\lambda_w^2} (1+\lambda_w^2)  B_{\max}^2L^2  \frac{1}{S}\phi_z^r
            + \frac{1+3\lambda_w^2}{4} \phi_u^r+25 \frac{1+7\lambda_w^2}{1-\lambda_w^2} N\frac{\sigma^2}{|\mathcal{I}|}
            \notag\\
            &~~~+ 397\frac{1+7\lambda_w^2}{1-\lambda_w^2} E^2 B_{\max}^2L^2  \alpha^2 \mathbb{E} [\|\bar \nabla_{\thetab}^r  \|^2 ]  +217 \frac{1+7\lambda_w^2}{1-\lambda_w^2} E^2N^3B_{\max}^4L^2  \beta^2  \mathbb{E}[  \|\bar\ub^r  \|^2],\label{eq: u consensus HBL}
    \end{align}
    \end{subequations}
\end{Lemma}

We establish a key descent inequality that characterizes the decrease property of the global objective function $F$ over the iterations from $r=1$ to some $T\geq1$.

Define
\begin{align}
\tilde{F}(\bar\xb^{r}, \bar\thetab^{r} )\triangleq \mathbb{E}\bigg[ \frac{1}{NS} \sum_{i=1}^S f\bigg(\sum_{n=1}^N  \Bb_{n, i} \bar\xb^{r},  \bar\thetab^{r}  \bigg) \bigg].
\end{align}
\begin{Lemma} \label{lem: descent lemma HBL}
(Descent lemma) Let Assumption \ref{assum network} to  Assumption \ref{assum: Lipschitz}, and Assumption \ref{assum: stochastic HBL} hold. For $\forall r\geq0$, we have
\begin{align}\label{eq: f HBL}
     &\tilde{F}(\bar\xb^{r+1}, \bar\thetab^{r+1} ) - \tilde{F}(\bar\xb^{r}, \bar\thetab^{r} )\notag\\
        & \leq - \!\frac{\alpha E}{2}\bigg\{
                1 -  10 E L \alpha
                - \bigg[ \frac{2\alpha L^2}{ N^3} +\frac{2\beta
                    B_{\max}^2 L^2}{N } \bigg]
                   64E(E-1)N\alpha
            \bigg\}   \mathbb{E}[ \|\bar \nabla_{\thetab}^r\|^2]
            -   \frac{\beta E}{2} \mathbb{E}[ \|\bar \nabla_{\xb}^r\|^2]\notag\\
         &~~~ - \frac{\beta}{2}    \bigg\{
                    1- 10E^2L\beta
                    - \bigg[ \frac{2\alpha L^2}{ N^3}
             +\frac{2\beta B_{\max}^2 L^2}{N } \bigg]
                    96E (E-1)N^2B_{\max}^2 \beta
         \bigg\}  \mathbb{E}  [ \| \bar \ub^r\|^2  ] \notag\\
    &~~~+ \bigg\{ \frac{2\alpha L^2}{ N^3}
             +\frac{2\beta B_{\max}^2 L^2}{N } \bigg\}
             [1+ 7(E-1)\lambda_w^2] \phi_\theta^r \notag\\
    &~~~ + \bigg\{
             \frac{2\alpha L^2}{ N^3}
             +\frac{2\beta B_{\max}^2 L^2}{N }
        \bigg\}  [E + 200(E-1) \lambda_w^2] N^2B_{\max}^2\phi_x^r \notag\\
    &~~~+ \bigg\{ \frac{2\alpha L^2}{ SN^3}
             +\frac{2\beta B_{\max}^2 L^2}{SN } \bigg\}
             [1+ 32(E-1)\lambda_w^2] \phi_z^r \notag\\
    &~~~ + \bigg\{
             \frac{2\alpha L^2}{ N^3}
             +\frac{2\beta B_{\max}^2 L^2}{N }
        \bigg\} 128(E-1) [1+ E(E-1) ]\lambda_w^2 N^2B_{\max}^2\beta^2 \phi_u^r \notag\\
    &~~~ +\frac{3EL\alpha^2  +9(E-1)L\beta^2}{N}\frac{\sigma^2}{|\mathcal{I}|} ,
\end{align}
\end{Lemma}
\subsection{LTI system}
Based on  Lemma \ref{lem: consensus HBL} and Lemma \ref{lem: tracking HBL}, we now summarize the iterate contraction into an LTI system which characterizes the convergence of consensus and gradient tracking processes in hybrid data setting.

Here we define the system matrix $\hat{\Ab}$, coefficient matrix $\hat{\Cb}$ of vector $\hat{\eb}^r$ as
\begin{subequations}
\begin{align}
& \hat{\Ab} \triangleq \begin{bmatrix}
    a_1 & a_2\alpha^2 & a_3\alpha^2 & a_4\alpha^2\beta^2 \\
    b_1\beta^2 & b_2  & b_3\beta^2 & b_4\beta^2 \\
    c_1\beta^2 & c_2\beta^2  & c_3  & c_4 \beta^2 \\
    d_1 & d_{21}\alpha^2 + d_{22}\beta^2 & d_3  & d_4
    \end{bmatrix},\\
&\hat{\Cb} \triangleq \begin{bmatrix}
    a_5\alpha^2 & a_6\alpha^2 & a_7 \alpha^2\beta^2 \\
    b_5\beta^2 & b_6 \alpha^2\beta^2  & b_7\beta^4   \\
    c_5\beta^2 & c_6  \alpha^2\beta^2 & c_7 \beta^2  \\
    d_5 & d_6\alpha^2 & d_7\beta^2
    \end{bmatrix},\\
&\hat{\eb}^r  \triangleq \begin{bmatrix}
    \frac{\sigma^2}{|\mathcal{I}|} \\
    \mathbb{E}[\|\bar\nabla_{\thetab}^r\|^2] \\
    \mathbb{E}[\|\bar\ub^r\|^2]
    \end{bmatrix},
\end{align}
\end{subequations}
where the parameters $a_1$-$a_7$, $b_1$-$b_7$, $c_1$-$c_7$, $d_1$-$d_7$ can be correspondingly determined in Lemma \ref{lem: consensus HBL} and Lemma \ref{lem: tracking HBL}.

Then, one can concisely write the LTI system as
\begin{align}
\hat{ \phib}^{r+1} \leq \hat{\Ab}\phib^r + \hat{\Cb}\hat{\eb}^r,
\end{align}
where the inequality is element-wise.

It is easy to determine that there exists a positive vector $\hat{\bsb} = [\hat{s}_1, \ldots, \hat{s}_4]^\top \in \mathbb{R}^4$ satisfying
\begin{align}
\hat{s}_1 \leq \frac{1-d_4}{5d_1}\hat{s}_4,~
\hat{s}_3 \leq \frac{1-d_4}{5d_3}\hat{s}_4,
\end{align}
such that $\Ab \hat{\bsb} < \hat{\bsb}$ as long as
\begin{align}\label{eq: condi alp bet HBL}
\alpha &\leq \min\bigg\{\frac{\sqrt{(1-d_4)\hat{s}_4}}{\sqrt{5d_{21}\hat{s}_2}},
\frac{\sqrt{(1-a_1)\hat{s}_1}}{\sqrt{a_{2}\hat{s}_2+ a_{3}\hat{s}_3+a_{4}\hat{s}_4 +1}}
\bigg\}\\
\beta & \leq \min\bigg\{
\frac{\sqrt{(1-d_4)\hat{s}_4}}{\sqrt{5d_{22}\hat{s}_2}},
 \frac{\sqrt{(1-b_2)\hat{s}_2}}{\sqrt{b_{1}\hat{s}_1+ b_{3}\hat{s}_3+b_{4}\hat{s}_4 +1}},
 \frac{\sqrt{(1-c_3)\hat{s}_3}}{\sqrt{c_{1}\hat{s}_1+ c_{2}\hat{s}_2+c_{4}\hat{s}_4 +1}}
\bigg\}.
\end{align}
Thus, under \eqref{eq: condi alp bet HBL},
the spectral radius of $\hat{\Ab}$ satisfies $\rho(\hat{\Ab})<1$
according to \cite[Corollary 8.1.29]{horn2012matrix}.
Then, we solve the inverse of system matrix $\hat{\Ab}$ in the following Lemma.

\begin{Lemma}\label{lem: det I -G HBL}
 For sufficiently small $\alpha$ and $\beta$ satisfying
 \begin{align}
& \alpha \leq   \min \bigg\{
      \frac{1}{\sqrt{d_{21}}} ,
    \frac{1}{\sqrt{d_{22}}} ,
   \frac{\sqrt{c_2d_3}}{ \sqrt{a_3 c_2  d_1}},
        \frac{\sqrt{c_2d_3}}{ \sqrt{ a_2c_1  d_3 }},
    \frac{b_4}{4\sqrt{a_4}}, \frac{b_4}{4\sqrt{a_3}}
\bigg\}\!,\\
&\beta  \leq   \min \bigg\{
 \frac{\sqrt{c_1}}{2\sqrt{b_1c_2}},
        \frac{\sqrt{c_1}}{2\sqrt{b_1c_4}},
 \frac{\sqrt{a_3}}{2 \sqrt{a_2b_3}}, \frac{\sqrt{a_3}}{2\sqrt{ a_4 d_3}},  \frac{\sqrt{a_3}}{2\sqrt{ a_4 b_3}},
  \frac{\sqrt{a_3}}{2\sqrt{a_2 b_4d_3}}, \frac{1}{\sqrt{d_{22}}},
 \frac{\sqrt{b_4}}{2\sqrt{ b_3 c_4}},
         \frac{\sqrt{b_4}}{2\sqrt{ b_3 c_1}},\notag\\
            &~~~ ~~~ ~~~~~~
          \frac{\sqrt{b_4}}{2\sqrt{ b_1 c_4}}
\bigg\},
\end{align}
 the determinant $(\Ib- \hat{\Ab})^{-1}$ has an element-wise upper bound
\begin{align}\label{eq: upper bound det I-G HBL}
    &(\mathbf{I}-\Ab)^{-1} \leq
    \frac{1 }{(1-\lambda_{w}^2)^4}\hat{\Ab}^{*},
\end{align}
where
\begin{align}
 &   \hat{\Ab}^* \triangleq\!
   \!\begin{bmatrix}
        1  &  a_2c_4d_3\alpha^2\beta^2  & 2a_3\alpha^2 & a_4b_3c_2\alpha^2\beta^6 \\
        b_1c_4d_3\beta^4  & 1  & a_4b_3d_1\alpha^2\beta^4 & 2b_4\beta^2 \\
        ( c_4d_1 \!+\!2c_1)\beta^2
         &   a_4c_2 d_1\alpha^2\beta^4&1 & a_2b_1c_4\alpha^2\beta^4\\
        b_3c_2 d_1\beta^4
        & (2 c_2d_3\!+ \!  d_{22}) \beta^2\!+\!d_{21}\alpha^2 & a_2b_1d_3\alpha^2\beta^2 & 1
    \end{bmatrix}\!.\label{eq: A 1 2}
\end{align}
 \end{Lemma}

Similar to \eqref{eq: sum phi FL}, we have
\begin{align}\label{eq: sum phi HBL}
\sum_{r=0}^{T-1}\hat{\phib}^{r}
& \leq    (\mathbf{I} - \hat{\Ab})^{-1}\hat{\phib}^{0}
+ (\mathbf{I} - \hat{\Ab})^{-1}\sum_{r=0}^{T-1}\hat{\Cb}\hat{\eb}^{r}.
\end{align}
Let $\hat{\phi}_\theta^0 = \hat{\phi}_x^0 = 0$. Then, by inserting \eqref{eq: upper bound det I-G HBL} into \eqref{eq: sum phi HBL}, we have the following Lemma.

\begin{Lemma}\label{lem: sum r0 T-1 xthetazu}
For sufficiently small $\alpha$ and $\beta$ satisfying
\begin{align}
\alpha \leq &\min \bigg\{
    \frac{\sqrt{a_6}}{2\sqrt{a_2d_3}},
     \frac{\sqrt{a_6}}{4\sqrt{a_3}},
    \frac{\sqrt{a_6}}{2\sqrt{a_4d_6}},
    \frac{\sqrt{a_7}}{2\sqrt{a_2d_3}},
    \frac{\sqrt{a_7}}{4\sqrt{a_4}},
    \frac{\sqrt{b_5}}{2\sqrt{a_5d_3}},
     \frac{\sqrt{b_5}}{ 2\sqrt{a_4d_1}},
    \frac{\sqrt{b_5}}{ 2\sqrt{a_6d_3}},
    \frac{\sqrt{b_5}}{ 2\sqrt{a_7d_3}} \notag\\
    & ~~~
     \frac{\sqrt{c_5}}{3\sqrt{a_5(c_4d_1+2c_1)}},
     \frac{\sqrt{c_5}}{ 3\sqrt{a_4d_1}},
    \frac{\sqrt{c_5}}{ 3\sqrt{a_2d_5}},
    \frac{\sqrt{c_6}}{ 2\sqrt{a_4d_1}},
    \frac{\sqrt{c_6}}{ 2\sqrt{a_1d_6}},
    \frac{\sqrt{c_7}}{ 3\sqrt{a_7}},
    \frac{\sqrt{c_7}}{ 3\sqrt{a_4d_1}},
    \frac{\sqrt{c_7}}{ 3\sqrt{a_2d_7}}\notag\\
    & ~~~
    \frac{\sqrt{d_5}}{2\sqrt{a_5 d_1 }},
     \frac{\sqrt{d_5}}{ 2\sqrt{d_{21}}},
    \frac{\sqrt{d_5}}{ 2\sqrt{a_2d_3}},
    \frac{\sqrt{d_6}}{ 2\sqrt{d_{21}}},
    \frac{\sqrt{d_6}}{ 2\sqrt{a_2d_3}},
    \frac{\sqrt{d_7}}{ 2\sqrt{a_7d_1}},
    \bigg(\frac{d_7}{b_7(2c_2d_3+d_{22})}\bigg)^{\frac{1}{4}},
    \frac{\sqrt{d_7}}{ 2\sqrt{d_{21}}},\notag\\
    &~~~
    \frac{\sqrt{d_7}}{ 2\sqrt{a_2d_3}}
\bigg\}, \\
  \beta \leq & \min \bigg\{
      \bigg( \frac{ a_5}{ a_2b_5c_4d_3}\bigg)^{\frac{1}{4}},
      \frac{\sqrt{a_5}}{4\sqrt{a_3c_5}},
      \bigg( \frac{ a_5}{ a_4b_3c_2 }\bigg)^{\frac{1}{6}},
      \frac{\sqrt{a_6}}{2\sqrt{b_6c_4}},
      \frac{\sqrt{a_6}}{4\sqrt{ c_6}},
      \frac{\sqrt{a_6}}{2\sqrt{b_3c_2}},
      \frac{\sqrt{a_7}}{2\sqrt{b_7c_4}},
       \frac{\sqrt{a_7}}{2\sqrt{b_3c_2}},\notag\\
       &~~~
       \frac{\sqrt{b_5}}{ 2\sqrt{b_1c_4}} ,
        \frac{\sqrt{b_5}}{ 2\sqrt{b_3c_5}},
          \frac{\sqrt{b_5}}{ 2\sqrt{b_3c_6}},
            \frac{\sqrt{b_5}}{ 2\sqrt{b_3c_7}}
        \frac{\sqrt{c_5}}{ 2\sqrt{b_5c_2}},
        \frac{\sqrt{c_5}}{ 3\sqrt{b_1c_4}},
        \frac{\sqrt{c_6}}{ 2\sqrt{b_6c_2}},
            \frac{\sqrt{c_6}}{ 2\sqrt{b_1c_4}},
            \frac{\sqrt{c_7}}{ 2\sqrt{b_1c_4}} ,\notag\\
    &~~~
        \frac{\sqrt{c_7}}{ 2\sqrt{c_4d_1+2c_1}},
        \frac{\sqrt{c_7}}{ 3\sqrt{b_7c_2}}
       \frac{\sqrt{d_5}}{2\sqrt{b_3 c_2 }}, \frac{\sqrt{d_5}}{ 2\sqrt{b_5}},
      \bigg(\frac{d_5}{2b_5(2c_2d_3+d_{22})}\bigg)^{\frac{1}{4}},
    \frac{\sqrt{d_5}}{ 2\sqrt{b_1c_5}},
    \bigg(\frac{d_6}{2a_6b_3c_2d_1 }\bigg)^{\frac{1}{4}},\notag\\
    &~~~
    \bigg(\frac{d_6}{2b_6 (2c_2d_3+d_{22}) }\bigg)^{\frac{1}{4}},
    \frac{\sqrt{d_6}}{ 2\sqrt{b_6}},
    \frac{\sqrt{d_6}}{ 2\sqrt{b_1c_6}},
    \frac{\sqrt{d_7}}{ 2\sqrt{b_3c_2}},
    \frac{\sqrt{d_7}}{ 2\sqrt{b_7}},
    \frac{\sqrt{d_7}}{ 2\sqrt{b_1c_7}}
    \bigg\},
\end{align}
we have
\begin{subequations}\label{eq: sum all HBL}
\begin{align}
\sum_{r=0}^{T-1}\phi_\theta^r &   \leq \frac{1}{(1-\lambda_{w}^2)^4} \bigg\{
        2a_5\alpha^2 T\frac{\sigma^2}{|\mathcal{I}|} + 2a_6\alpha^2 \sum_{r=0}^{T-1} \mathbb{E}[\|\bar\nabla_{\thetab}^r\|^2]
        + 2(a_7+a_3c_7)\alpha^2\beta^2\sum_{r=0}^{T-1} \mathbb{E}[\|\bar\ub^r\|^2] \notag\\
        &~~~ +2a_3\alpha^2\phi_z^0 + a_4b_3c_2\alpha^2\beta^6\phi_u^0 \bigg\},\label{eq: sum theta HBL}\\
\sum_{r=0}^{T-1}\phi_x^r &   \leq \frac{1}{(1-\lambda_{w}^2)^4} \bigg\{
        2(b_5+b_4d_5)\beta^2 T\frac{\sigma^2}{|\mathcal{I}|}
        + 2(b_6+b_4d_6)\alpha^2\beta^2 \sum_{r=0}^{T-1}\mathbb{E}[\|\bar\nabla_{\thetab}^r\|^2]  \notag\\
        &~~~ + 2(b_7+b_4d_7) \beta^4 \sum_{r=0}^{T-1}\mathbb{E}[\|\bar\ub^r\|^2] + a_4b_3d_1\alpha^2\beta^4 \phi_z^0
         + 2b_4 \beta^2\phi_u^0 \bigg\},\label{eq: sum x HBL}\\
  \sum_{r=0}^{T-1}\phi_z^r &   \leq \frac{1}{(1-\lambda_{w}^2)^4} \bigg\{
        2c_5\beta^2 T\frac{\sigma^2}{|\mathcal{I}|}
        + 2[a_6( c_4d_1 + 2c_1)+ c_6]\alpha^2\beta^2 \sum_{r=0}^{T-1}\mathbb{E}[\|\bar\nabla_{\thetab}^r\|^2]   \notag\\
        &~~~+ 2c_7 \beta^2 \sum_{r=0}^{T-1}\mathbb{E}[\|\bar\ub^r\|^2] +   \phi_z^0
         + a_2b_1c_4\alpha^2 \beta^4\phi_u^0 \bigg\},\label{eq: sum z HBL}
        \end{align}
        and
        \begin{align}
  \sum_{r=0}^{T-1}\phi_u^r &   \leq \frac{1}{(1-\lambda_{w}^2)^4} \bigg\{
        2d_5  T\frac{\sigma^2}{|\mathcal{I}|}
        + 2d_6\alpha^2 \sum_{r=0}^{T-1} \mathbb{E}[\|\bar\nabla_{\thetab}^r\|^2]
         + 2d_7 \beta^2\sum_{r=0}^{T-1} \mathbb{E}[\|\bar\ub^r\|^2]  \notag\\
         &~~~+ a_2b_1d_3 \alpha^2 \beta^2 \phi_z^0
         + \phi_u^0 \bigg\}.\label{eq: sum u HBL}
\end{align}
\end{subequations}
\end{Lemma}

\subsection{Proof of Theorem \ref{thm: HBL}}
By summing \eqref{eq: f HBL} up from $r=0$ to $r=T-1$, and then inserting \eqref{eq: sum all HBL} into the obtained results, we have
\begin{align}
     &\tilde{F}(\bar\xb^{T}, \bar\thetab^{T} ) - \tilde{F}(\bar\xb^{0}, \bar\thetab^{0} )\notag\\
        & \leq - \!\frac{\alpha E}{2}\bigg\{
                1 -  10 E L \alpha
                - \bigg[ \frac{2\alpha L^2}{ N^3} +\frac{2\beta
                    B_{\max}^2 L^2}{N } \bigg]
                   64E(E-1)N\alpha
            \bigg\}   \sum_{r=0}^{T-1}\mathbb{E}[ \|\bar \nabla_{\thetab}^r\|^2]
            \notag\\
         &~~~~-   \frac{\beta E}{2}  \sum_{r=0}^{T-1}\mathbb{E}[ \|\bar \nabla_{\xb}^r\|^2]\notag\\
         &~~~ - \frac{\beta}{2}    \bigg\{
                    1- 10E^2L\beta
                    - \bigg[ \frac{2\alpha L^2}{ N^3}
             +\frac{2\beta B_{\max}^2 L^2}{N } \bigg]
                     96E (E-1)N^2B_{\max}^2 \beta
         \bigg\}   \sum_{r=0}^{T-1}\mathbb{E}  [ \| \bar \ub^r\|^2  ] \notag\\
    &~~~+ \bigg\{ \frac{2\alpha L^2}{ N^3}
             +\frac{2\beta B_{\max}^2 L^2}{N } \bigg\}
             [1+ 7(E-1)\lambda_w^2]
             \bigg\{
                \frac{1}{(1-\lambda_{w}^2)^4} \{
                2a_5\alpha^2 T\frac{\sigma^2}{|\mathcal{I}|} + 2a_6\alpha^2 \sum_{r=0}^{T-1} \mathbb{E}[\|\bar\nabla_{\thetab}^r\|^2] \notag\\
                &~~~~~~ + 2(a_7+a_3c_7)\alpha^2\beta^2 \sum_{r=0}^{T-1}\mathbb{E}[\|\bar\ub^r\|^2]  +2a_3\alpha^2\phi_z^0 + a_4b_3c_2\alpha^2\beta^6\phi_u^0 \}
             \bigg\}\notag\\
    &~~~ + \bigg\{
             \frac{2\alpha L^2}{ N^3}
             +\frac{2\beta B_{\max}^2 L^2}{N }
        \bigg\}  [E + 200(E-1) \lambda_w^2] N^2B_{\max}^2
            \bigg\{ \frac{1}{(1-\lambda_{w}^2)^4} \{
                2(b_5+b_4d_5)\beta^2T \frac{\sigma^2}{|\mathcal{I}|} \notag\\
                &~~~~~~
                + 2(b_6+b_4d_6)\alpha^2\beta^2 \sum_{r=0}^{T-1}\mathbb{E}[\|\bar\nabla_{\thetab}^r\|^2]
               + 2(b_7+b_4d_7) \beta^4 \sum_{r=0}^{T-1}\mathbb{E}[\|\bar\ub^r\|^2] \notag\\
                &~~~~~~ + a_4b_3d_1\alpha^2\beta^2 \phi_z^0
             + 2b_4 \beta^2\phi_u^0 \}
             \bigg\} \notag\\
    &~~~+ \bigg\{ \frac{2\alpha L^2}{ SN^3}
             +\frac{2\beta B_{\max}^2 L^2}{SN } \bigg\}
             [1+ 32(E-1)\lambda_w^2]
             \bigg\{
                 \frac{1}{(1-\lambda_{w}^2)^4} \{
                2c_5\beta^2 T\frac{\sigma^2}{|\mathcal{I}|}
                + 2[a_6(  c_4d_1 + 2c_1)+ c_6]\notag\\
                &~~~ ~~~ ~~~\alpha^2\beta^2 \sum_{r=0}^{T-1}\mathbb{E}[\|\bar\nabla_{\thetab}^r\|^2]  + 2c_7 \beta^2 \sum_{r=0}^{T-1}\mathbb{E}[\|\bar\ub^r\|^2]
                +   \phi_z^0
                 + a_2b_1c_4\alpha^2 \beta^4\phi_u^0 \}
             \bigg\}\notag\\
    &~~~ + \bigg\{
             \frac{2\alpha L^2}{ N^3}
             +\frac{2\beta B_{\max}^2 L^2}{N }
        \bigg\} 128(E-1) [1+ E(E-1) ]\lambda_w^2 N^2B_{\max}^2\beta^2
         \bigg\{
             \frac{1}{(1-\lambda_{w}^2)^4} \{
            2d_5  T\frac{\sigma^2}{|\mathcal{I}|}
           \notag\\
            &~~~~~~ + 2d_6\alpha^2 \sum_{r=0}^{T-1} \mathbb{E}[\|\bar\nabla_{\thetab}^r\|^2] + 2d_7 \beta^2 \sum_{r=0}^{T-1}\mathbb{E}[\|\bar\ub^r\|^2]  + a_2b_1d_3 \alpha^2 \beta^2 \phi_z^0
             + \phi_u^0 \}
         \bigg\} \notag\\
    &~~~ +[3EL\alpha^2  +9(E-1)L\beta^2]T\frac{\sigma^2}{|\mathcal{I}|} \notag\\
& = - \!\frac{\alpha E}{2}\bigg\{
        1\! -  \!10 E L \alpha
        \!-\! \bigg( \frac{2\alpha L^2}{ N^3}\! +\!\frac{2\beta B_{\max}^2 L^2}{N } \bigg) 64E(E\!-\!1)N\alpha\notag\\
        &~~~~~~
         -\bigg( \frac{2\alpha L^2}{ N^3}
             +\frac{2\beta B_{\max}^2 L^2}{N } \bigg)
             [1+ 7(E-1)\lambda_w^2]
            \frac{4 a_6\alpha }{E(1-\lambda_{w}^2)^4}\notag\\
        &~~~~~~ 
         -\bigg(
             \frac{2\alpha L^2}{ N^3}
             +\frac{2\beta B_{\max}^2 L^2}{N }
        \bigg)  [E + 200(E-1) \lambda_w^2 ]
        N^2B_{\max}^2
            \frac{1}{(1-\lambda_{w}^2)^4}
               4(b_6+b_4d_6)\alpha \beta^2 \notag\\
        &~~~~~~ 
        -\bigg( \frac{2\alpha L^2}{ SN^3}
             +\frac{2\beta B_{\max}^2 L^2}{SN } \bigg)
             [1+ 32(E-1)\lambda_w^2]
                 \frac{1}{(1-\lambda_{w}^2)^4}
                  4[a_6( c_4d_1+  2c_1)+ c_6]\alpha \beta^2\notag\\
        &~~~~~~ 
        -\bigg(
             \frac{2\alpha L^2}{ N^3}
             +\frac{2\beta B_{\max}^2 L^2}{N }
        \bigg) 128(E-1) [1+ E(E-1) ]\lambda_w^2
        \frac{4 d_6 N^2B_{\max}^2\alpha \beta^2}{(1-\lambda_{w}^2)^4}
    \bigg\}   \sum_{r=0}^{T-1}\mathbb{E}[ \|\bar \nabla_{\thetab}^r\|^2]\notag\\
    &~~~ 
    - \frac{\beta}{2}    \bigg\{
        1- 10E^2L\beta - \bigg( \frac{2\alpha L^2}{ N^3} +\frac{2\beta B_{\max}^2 L^2}{N } \bigg) 96E (E-1)N^2B_{\max}^2 \beta\notag\\
    &~~~~~~ 
        -\bigg( \frac{2\alpha L^2}{ N^3}
             +\frac{2\beta B_{\max}^2 L^2}{N } \bigg)
        [1+ 7(E-1)\lambda_w^2]
              \frac{4(a_7+a_3c_7)\alpha^2\beta}{(1-\lambda_{w}^2)^4}\notag\\
    &~~~~~~ 
        -\bigg(
             \frac{2\alpha L^2}{ N^3}
             +\frac{2\beta B_{\max}^2 L^2}{N }
        \bigg)  [E + 200(E-1) \lambda_w^2] N^2B_{\max}^2
              \frac{4(b_7+b_4d_7) \beta^3 }{(1-\lambda_{w}^2)^4}  \notag\\
    &~~~~~~ 
        -\bigg( \frac{2\alpha L^2}{ SN^3}
             +\frac{2\beta B_{\max}^2 L^2}{SN } \bigg)
             [1+ 32(E-1)\lambda_w^2]  \frac{4 c_7 \beta}{(1-\lambda_{w}^2)^4}\notag\\
    &~~~~~~ 
        - \bigg(
             \frac{2\alpha L^2}{ N^3}
             +\frac{2\beta B_{\max}^2 L^2}{N }
        \bigg) 128(E-1) [1+ E(E-1) ]\lambda_w^2
        N^2B_{\max}^2\beta^2
             \frac{4 d_7 \beta^2}{(1-\lambda_{w}^2)^4}
         \bigg\}   \sum_{r=0}^{T-1}\mathbb{E}  [ \| \bar \ub^r\|^2  ] \notag\\
    &~~~ -   \frac{\beta E}{2}  \sum_{r=0}^{T-1}\mathbb{E}[ \|\bar \nabla_{\xb}^r\|^2]\notag\\
    &~~~ 
    +\bigg\{[3EL\alpha^2  +9(E-1)L\beta^2]
        + \bigg( \frac{2\alpha L^2}{ N^3} +\frac{2\beta B_{\max}^2 L^2}{N } \bigg) [1+ 7(E-1)\lambda_w^2]  \frac{2 a_5\alpha^2 }{(1-\lambda_{w}^2)^4}  \notag\\
    &~~~ ~~~ 
        + \bigg( \frac{2\alpha L^2}{ N^3} +\frac{2\beta B_{\max}^2 L^2}{N } \bigg)  [E + 200(E-1) \lambda_w^2]
         \frac{2(b_5+b_4d_5)N^2B_{\max}^2\beta^2}{(1-\lambda_{w}^2)^4}  \notag\\
    &~~~ ~~~ 
        +\bigg( \frac{2\alpha L^2}{ SN^3}
             +\frac{2\beta B_{\max}^2 L^2}{SN } \bigg)
             [1+ 32(E-1)\lambda_w^2]  \frac{2c_5\beta^2 }{(1-\lambda_{w}^2)^4}  \notag\\
    &~~~ ~~~ 
       + \bigg( \frac{2\alpha L^2}{ N^3} +\frac{2\beta B_{\max}^2 L^2}{N } \bigg) 128(E-1)
        [1+ E(E-1) ] \lambda_w^2\frac{2d_5  N^2B_{\max}^2\beta^2 }{(1-\lambda_{w}^2)^4}
        \bigg\}T\frac{\sigma^2}{|\mathcal{I}|} \notag\\
    &~~~ 
        +\bigg\{
            \bigg( \frac{2\alpha L^2}{ N^3} +\frac{2\beta B_{\max}^2 L^2}{N } \bigg) [1+ 7(E-1)\lambda_w^2]
            \frac{2 a_3\alpha^2}{(1-\lambda_{w}^2)^4} \notag\\
        &~~~ ~~~ 
            +\bigg( \frac{2\alpha L^2}{ N^3} +\frac{2\beta B_{\max}^2 L^2}{N } \bigg)  [E + 200(E-1) \lambda_w^2]
             \frac{  a_4b_3d_1 N^2B_{\max}^2\alpha^2\beta^2  }{(1-\lambda_{w}^2)^4}\notag\\
        &~~~ ~~~ 
            +\bigg(  \frac{2\alpha L^2}{ N^3} +\frac{2\beta B_{\max}^2 L^2}{N } \bigg) 128(E-1) [1+ E(E-1) ] \lambda_w^2
            \frac{  a_2b_1d_3 N^2B_{\max}^2 \alpha^2 \beta^4}{(1-\lambda_{w}^2)^4}\notag\\
        &~~~ ~~~ 
            + \bigg( \frac{2\alpha L^2}{ SN^3} +\frac{2\beta B_{\max}^2 L^2}{SN } \bigg) \frac{ [1+ 32(E-1)\lambda_w^2]}{(1-\lambda_{w}^2)^4}
        \bigg\}\phi_z^0 \notag\\
    &~~~ 
        +\bigg\{
            \bigg( \frac{2\alpha L^2}{ N^3} +\frac{2\beta B_{\max}^2 L^2}{N } \bigg) [1\!+ \!7(E\!-\!1)\lambda_w^2]
            \frac{ a_4b_3c_2\alpha^2\beta^6}{(1\!-\!\lambda_{w}^2)^4} \notag\\
        &~~~ ~~~ 
            + \bigg( \frac{2\alpha L^2}{ SN^3} \!+\!\frac{2\beta B_{\max}^2 L^2}{SN } \bigg) [1\!+ \!32(E-1)\lambda_w^2]  \frac{ a_2b_1c_4\alpha^2 \beta^4}{(1\!-\!\lambda_{w}^2)^4}\notag\\
        &~~~ ~~~ 
            +\bigg(  \frac{2\alpha L^2}{ N^3} +\frac{2\beta B_{\max}^2 L^2}{N } \bigg) 128(E-1) [1+ E(E-1) ]\lambda_w^2
             \frac{    N^2B_{\max}^2 \beta^2 }{(1-\lambda_{w}^2)^4}\notag\\
        &~~~ ~~~ 
            +\bigg( \frac{2\alpha L^2}{ N^3} +\frac{2\beta B_{\max}^2 L^2}{N } \bigg)  [E + 200(E-1)  \lambda_w^2]
            \frac{2b_4 N^2B_{\max}^2 \beta^2  }{(1-\lambda_{w}^2)^4}
        \bigg\}\phi_u^0 \label{eq: descent insert sum 1}\\
& \leq -\frac{\alpha E}{4} \sum_{r=0}^{T-1}\mathbb{E}[ \|\bar \nabla_{\thetab}^r\|^2]
    - \frac{\beta E}{4} \sum_{r=0}^{T-1}\mathbb{E}[ \|\bar \nabla_{\xb}^r\|^2]
     - \frac{\beta}{4} \sum_{r=0}^{T-1}\mathbb{E}  [ \| \bar \ub^r\|^2  ] \notag\\
    &~~~ 
        + [4EL\alpha^2  +10(E-1)L\beta^2]T\frac{\sigma^2}{|\mathcal{I}|} \notag\\
    &~~~ 
        +  \bigg( \frac{2\alpha L^2}{ N^3} +\frac{2\beta B_{\max}^2 L^2}{N } \bigg)\frac{1  }{(1-\lambda_{w}^2)^4}
         \bigg\{    [1+ 7(E-1)\lambda_w^2] 2 a_5\alpha^2
            + \frac{1}{S}[1+32(E-1)\lambda_w^2] 2c_5\beta^2
    \notag\\
    &~~~ ~~~
        +  [E + 200(E-1) \lambda_w^2]4b_4d_5 N^2B_{\max}^2\beta^2 \notag\\
         &~~~~~~    
       +  128(E-1) [1+ E(E-1)]\lambda_w^2 2d_5  N^2B_{\max}^2\beta^2
        \bigg\}T\frac{\sigma^2}{|\mathcal{I}|}  \notag\\
    &~~~
        + \bigg( \frac{2\alpha L^2}{ SN^3} +\frac{2\beta B_{\max}^2 L^2}{SN } \bigg) [1+ 32(E-1)\lambda_w^2]  \frac{2}{(1-\lambda_{w}^2)^3} \phi_z^0 \notag\\
    &~~~
        +  \bigg( \frac{2\alpha L^2}{ N^3} +\frac{2\beta B_{\max}^2 L^2}{N } \bigg)   \frac{  N^2B_{\max}^2 \beta^2  }{(1-\lambda_{w}^2)^4}
         \bigg\{2b_4 [E + 200(E-1) \lambda_w^2] +1
           \notag\\
         &~~~~~~~~~  + 128(E-1) [1+ E(E-1) ]\lambda_w^2 \bigg\}
         \phi_u^0,\label{eq: sum f HBL}
\end{align}
where the last inequality holds if $\alpha$ and $\beta$ satisfy
\begin{align}
\alpha & \leq  \min \bigg\{
   \frac{N^3}{4L^2},
   \frac{1}{120EL},
   \frac{1}{768 E^2 N},
   \frac{E(1-\lambda_w^2)^4}{48a_6 [1+7(E-1)\lambda_w^2]},
   \frac{(1\!-\!\lambda_w^2)^4}{12N^2B_{\max}^2 [E\!+\!200(E\!-\!1)\lambda_w^2]},\notag\\
   &~~~
   \frac{(1\!-\!\lambda_w^2)^4}{48[1\!+\!32(E\!-\!1)\lambda_w^2]},
   \frac{(1-\lambda_w^2)^4}{6144 d_6 N^2B_{\max}^2 (E-1) [1+E(E-1)\lambda_w^2]},
   \frac{(1\!-\!\lambda_w^2)^2}{48 \sqrt{a_7\!+\!a_3 c_7} [1\!+\!7(E\!-\!1)\lambda_w^2]},\notag\\
   &~~~
   \frac{\sqrt{[1\!+\!32(E\!-\!1)\lambda_w^2]} }{8a_3 [1\!+\!7(E\!-\!1)\lambda_w^2]},
    \sqrt{\frac{1}{a_4d_1N^2B_{\max}^2}},
    \sqrt{\frac{1}{a_2d_3N^2B_{\max}^2}},
    \frac{1}{\sqrt{a_2c_4}}
\bigg\},\\
 \beta& \leq \min \bigg\{
        \frac{N}{4B_{\max}^2 L^2},
        \frac{1}{2\sqrt{b_6\!+\!b_4d_6}},
        \frac{1}{\sqrt{a_6(c_4d_1 \!+\! 2c_1)\!+\!c_6}},
        \frac{1}{120E^2L},
        \frac{1}{1152E(E-1)N^2B_{\max}^2},\notag\\
   &~~~
        \frac{(1\!-\!\lambda_w^2) }{ \{48 (b_7+b_4d_7)N^2B_{\max}^2 [E+200(E-1)\lambda_w^2]\}^{\frac{1}{3}}},
        \frac{(1\!-\!\lambda_w^2)^4 }{48Sc_7[1\!+\!32(E\!-\!1)\lambda_w^2] },\notag\\
   &~~~
        \frac{(1\!-\!\lambda_w^2)^2 }{78NB_{\max} \sqrt{d_7(E-1)[1+E(E-1)\lambda_w^2]} },
        \sqrt{\frac{[1+32(E-1)\lambda_w^2]}{4b_3 [E+200(E-1)\lambda_w^2]}}, \notag\\
   &~~~
        \bigg(\frac{1+32(E-1)\lambda_w^2}{512 b_1(E-1)[1+E(E-1)]\lambda_w^2}\bigg)^{\frac{1}{4}},
         \bigg(\frac{1}{2a_4b_3c_2  [1+7E(E-1) \lambda_w^2]}\bigg)^{\frac{1}{4}}, \notag\\
   &~~~
         \sqrt{\frac{1}{2b_1S [1+32(E-1)\lambda_w^2]}}
    \bigg\},
\end{align}
so that in \eqref{eq: descent insert sum 1},
the last $6$ negative terms in the bracket before $\textstyle\sum_{r=0}^{T-1} \mathbb{E}[\|\bar\nabla_{\thetab}^r\|^2]$ are smaller than $\textstyle \frac{1}{12}$,
the last $6$ negative terms in the bracket before $\textstyle\sum_{r=0}^{T-1} \mathbb{E}[\|\bar\ub^r\|^2]$ are smaller than $\textstyle \frac{1}{12}$,
the first $3$ terms are smaller than the last term before $\phi_z^0$,
the first $2$ terms are smaller than $\textstyle \frac{1}{2}$ before $\phi_u^0$.

By rearranging the above inequality and further inserting $\tilde{F}(\bar\xb^{T}, \bar\thetab^{T} ) \geq \underline{F}$ in Assumption \ref{assum: lower bound}, we have
\begin{align}\label{eq: sum bar theta x HBL}
& \frac{\alpha E}{4} \sum_{r=0}^{T-1}\mathbb{E}[ \|\bar \nabla_{\thetab}^r\|^2] + \frac{\beta E}{4} \sum_{r=0}^{T-1}\mathbb{E}[ \|\bar \nabla_{\xb}^r\|^2] \notag\\
&\leq \tilde{F}(\bar \thetab^0, \bar\xb^{0} ) - \underline{F}
          - \frac{\beta}{4} \sum_{r=0}^{T-1}\mathbb{E}  [ \| \bar \ub^r\|^2  ]
        + [4EL\alpha^2  +10(E-1)L\beta^2]T\frac{\sigma^2}{|\mathcal{I}|} \notag\\
    &~~~ 
        +  \bigg( \frac{2\alpha L^2}{ N^3} +\frac{2\beta B_{\max}^2 L^2}{N } \bigg)\frac{ 1 }{(1-\lambda_{w}^2)^4}
         \bigg\{   [1+ 7(E-1)\lambda_w^2] 2 a_5\alpha^2
            + \frac{1}{S}[1+32(E-1)\lambda_w^2] 2c_5\beta^2
    \notag\\
    &~~~ ~~~
        +  [E + 200(E-1) \lambda_w^2]4b_4d_5 N^2B_{\max}^2\beta^2 \notag\\
         &~~~~~~    
       +  128(E-1) [1+ E(E-1)]\lambda_w^2 2d_5  N^2B_{\max}^2\beta^2
        \bigg\}T\frac{\sigma^2}{|\mathcal{I}|}  \notag\\
    &~~~
        + \bigg( \frac{2\alpha L^2}{ SN^3} +\frac{2\beta B_{\max}^2 L^2}{SN } \bigg) [1+ 32(E-1)\lambda_w^2]  \frac{2}{(1-\lambda_{w}^2)^3} \phi_z^0 \notag\\
    &~~~
        + \! \bigg( \frac{2\alpha L^2}{ N^3} \!+\!\frac{2\beta B_{\max}^2 L^2}{N } \bigg)   \frac{  N^2B_{\max}^2 \beta^2  }{(1\!-\!\lambda_{w}^2)^4} \bigg\{2b_4 [E \!+ \!200(E\!-\!1) \lambda_w^2] \!+\!1  \!\notag\\
         &~~~  ~~~ +\! 128(E\!-\!1) [1\!+\! E(E\!-\!1) ]\lambda_w^2 \bigg\}
         \phi_u^0.
\end{align}

Besides, we perspectively solve the average gradients with respect to $\thetab$ and $\xb$ as follows.
\begin{align}\label{eq: gr theta HBL}
& \frac{1}{N} \sum_{n=1}^N\mathbb{E}\bigg[\bigg\| \frac{1}{NS}\sum_{i=1}^S\nabla_{\thetab} f\bigg(\sum_{t=1}^N \Bb_{t,i}\xb_n^r, \thetab_n^r\bigg)\bigg\|^2\bigg] \notag\\
& \leq 2 \frac{1}{N} \sum_{n=1}^N\mathbb{E}\bigg[\bigg\| \frac{1}{NS}\sum_{i=1}^S\bigg[\nabla_{\thetab} f\bigg(\sum_{t=1}^N \Bb_{t,i}\xb_n^r, \thetab_n^r\bigg)
-\nabla_{\thetab} f\bigg(\sum_{t=1}^N \Bb_{t,i}\bar \xb^r, \bar \thetab^r\bigg) \bigg\|^2\bigg] + 2 \mathbb{E}[\|\bar \nabla_{\thetab}^r\|^2] \notag\\
&  \leq 2 \frac{1}{N} \sum_{n=1}^N\frac{1}{N^2S}\sum_{i=1}^S L^2 \mathbb{E}\bigg[\bigg\|  \sum_{t=1}^N \Bb_{t,i}(\xb_n^r - \bar \xb^r)\bigg\|^2\bigg]
    + 2 \frac{1}{N} \sum_{n=1}^N\frac{1}{N^2S}\sum_{i=1}^S L^2 \mathbb{E} [\| \thetab_n^r - \bar \thetab^r  \|^2 ] \notag\\
    &~~~+ 2 \mathbb{E}[\|\bar \nabla_{\thetab}^r\|^2] \notag\\
& \leq   \frac{ 2 L^2    B_{\max}^2}{N}   \phi_x^r
    +   \frac{2 L^2}{N^3 } \phi_\theta^r + 2 \mathbb{E}[\|\bar \nabla_{\thetab}^r\|^2].
\end{align}

Similarly, we have
\begin{align}\label{eq: gr x HBL}
&\frac{1}{N} \sum_{n=1}^N \mathbb{E}\bigg[\bigg\| \frac{1}{NS}\sum_{i=1}^S\sum_{t=1}^N \Bb_{t,i}^\top \nabla_{\Bb_i\xb} f\bigg(\sum_{t=1}^N \Bb_{t,i}\xb_n^r, \thetab_n^r)\bigg\|^2\bigg] \notag\\
& \leq 2 \frac{1}{N} \sum_{n=1}^N\mathbb{E}\bigg[\bigg\| \frac{1}{NS}\sum_{i=1}^S\sum_{t=1}^N\Bb_{t,i}^\top \bigg[\nabla_{\Bb_i\xb} f\bigg(\sum_{t=1}^N \Bb_{t,i}\xb_n^r, \thetab_n^r\bigg)
 -\nabla_{\Bb_i\xb} f\bigg(\sum_{t=1}^N \Bb_{t,i}\bar \xb^r, \bar \thetab^r\bigg) \bigg\|^2\bigg] \notag\\
         &~~~~+ 2 \mathbb{E}[\|\bar \nabla_{\xb}^r\|^2] \notag\\
&  \leq 2 \frac{1}{N} \sum_{n=1}^N\frac{1}{N^2S}\sum_{i=1}^S N\sum_{t=1}^N B_{\max}^2 L^2 \mathbb{E}\bigg[\bigg\|  \sum_{t=1}^N \Bb_{t,i}(\xb_n^r - \bar \xb^r)\bigg\|^2\bigg]\notag\\
    &~~~
    + \frac{1}{N} \sum_{n=1}^N\frac{1}{N^2S}\sum_{i=1}^S N\sum_{t=1}^N B_{\max}^2 L^2 \mathbb{E} [\| \thetab_n^r - \bar \thetab^r  \|^2 ] + 2 \mathbb{E}[\|\bar \nabla_{\xb}^r\|^2] \notag\\
& \leq  2 N B_{\max}^4 L^2  \phi_x^r
    + \frac{2 B_{\max}^2 L^2}{N} \phi_\theta^r + 2 \mathbb{E}[\|\bar \nabla_{\xb}^r\|^2]
\end{align}

Multiplying $\frac{\alpha E}{10}$ and $\frac{\beta E}{10}$ on the both sides of \eqref{eq: gr theta HBL} and \eqref{eq: gr x HBL} respectively, and then summing the results up from $r=0$ to $r=T-1$, we obtain
\begin{align}\label{eq: gr sum theta x HBL 1}
&\sum_{r=0}^{T-1}\frac{\alpha E}{10}\frac{1}{N} \sum_{n=1}^N \mathbb{E}\bigg[\bigg\| \frac{1}{NS}\sum_{i=1}^S\nabla_{\thetab} f\bigg(\sum_{t=1}^N \Bb_{t,i}\xb_n^r, \thetab_n^r\bigg)\bigg\|^2\bigg] \notag\\
&~~~ + \sum_{r=0}^{T-1} \frac{\beta E}{10}\frac{1}{N} \sum_{n=1}^N \mathbb{E}\bigg[\bigg\| \frac{1}{NS}\sum_{i=1}^S\sum_{t=1}^N \Bb_{t,i}^\top  \nabla_{\Bb_i\xb} f\bigg(\sum_{t=1}^N \Bb_{t,i}\xb_n^r, \thetab_n^r)\bigg\|^2\bigg] \notag\\
& \leq
    \frac{   L^2    B_{\max}^2}{5N}( \alpha+N^2B_{\max}^2 \beta)  \sum_{r=0}^{T-1}\phi_x^r
    +   \frac{  L^2}{5N^3 } (\alpha + N^2B_{\max}^2 \beta)\sum_{r=0}^{T-1} \phi_\theta^r
+ \frac{ \alpha E}{5}\sum_{r=0}^{T-1} \mathbb{E}[\|\bar \nabla_{\thetab}^r\|^2]
    \notag\\
    &~~~  + \frac{ \beta E}{5}\sum_{r=0}^{T-1}  \mathbb{E}[\|\bar \nabla_{\xb}^r\|^2]
\end{align}
By inserting \eqref{eq: sum theta HBL} and \eqref{eq: sum x HBL} into \eqref{eq: gr sum theta x HBL 1}, we have
\begin{align}
&\sum_{r=0}^{T-1}\frac{\alpha E}{10} \frac{1}{N} \sum_{n=1}^N\mathbb{E}\bigg[\bigg\| \frac{1}{NS}\sum_{i=1}^S\nabla_{\thetab} f\bigg(\sum_{t=1}^N \Bb_{t,i}\xb_n^r, \thetab_n^r\bigg)\bigg\|^2\bigg] \notag\\
&~~~ + \sum_{r=0}^{T-1} \frac{\beta E}{10}\frac{1}{N} \sum_{n=1}^N\mathbb{E}\bigg[\bigg\| \frac{1}{NS}\sum_{i=1}^S\sum_{t=1}^N \Bb_{t,i}^\top  \nabla_{\Bb_i\xb} f\bigg(\sum_{t=1}^N \Bb_{t,i}\xb_n^r, \thetab_n^r)\bigg\|^2\bigg] \notag\\
& \leq
    \frac{   L^2    B_{\max}^2}{5N}( \alpha+N^2B_{\max}^2 \beta)  \frac{1}{(1-\lambda_{w}^2)^4}\bigg\{
        2(b_5+b_4d_5)\beta^2 T\frac{\sigma^2}{|\mathcal{I}|}
        + 2(b_6+b_4d_6)\alpha^2\beta^2 \sum_{r=0}^{T-1}\mathbb{E}[\|\bar\nabla_{\thetab}^r\|^2] \notag\\
        &~~~ ~~~+ 2(b_7+b_4d_7) \beta^4 \sum_{r=0}^{T-1}\mathbb{E}[\|\bar\ub^r\|^2]  + a_4b_3d_1\alpha^2\beta^2 \phi_z^0
         + 2b_4 \beta^2\phi_u^0
    \bigg\}\notag\\
    &~~~
    +   \frac{  L^2}{5N^3 } (\alpha + N^2B_{\max}^2 \beta)
     \frac{1}{(1-\lambda_{w}^2)^4} \bigg\{
        2a_5\alpha^2 T\frac{\sigma^2}{|\mathcal{I}|} + 2a_6\alpha^2 \sum_{r=0}^{T-1} \mathbb{E}[\|\bar\nabla_{\thetab}^r\|^2] \notag\\
        &~~~~~~ + 2(a_7+a_3c_7)\alpha^2\beta^2\sum_{r=0}^{T-1} \mathbb{E}[\|\bar\ub^r\|^2]  +2a_3\alpha^2\phi_z^0 + a_4b_3c_2\alpha^2\beta^6\phi_u^0
    \bigg\}
\notag\\
    &~~~  + \frac{ \alpha E}{5}\sum_{r=0}^{T-1} \mathbb{E}[\|\bar \nabla_{\thetab}^r\|^2]
    + \frac{ \beta E}{5}\sum_{r=0}^{T-1}  \mathbb{E}[\|\bar \nabla_{\xb}^r\|^2]\notag\\
& \leq \frac{ \alpha E}{4}\sum_{r=0}^{T-1} \mathbb{E}[\|\bar \nabla_{\thetab}^r\|^2]
    + \frac{ \beta E}{4}\sum_{r=0}^{T-1}  \mathbb{E}[\|\bar \nabla_{\xb}^r\|^2] 
        + \bigg\{
            \frac{   L^2    B_{\max}^2}{5N}( \alpha+N^2B_{\max}^2 \beta)    \frac{1}{(1-\lambda_{w}^2)^4}   2(b_7+b_4d_7) \beta^4\notag\\
        & ~~~~~~ 
            + \frac{  L^2}{5N^3 } (\alpha + N^2B_{\max}^2 \beta)
            \frac{1}{(1-\lambda_{w}^2)^4}   2(a_7+a_3c_7)\alpha^2\beta^2
    \bigg\}\sum_{r=0}^{T-1}\mathbb{E}[\|\bar\ub^r\|^2]\notag\\
    & ~~~ 
    + \bigg\{\frac{   L^2    B_{\max}^2}{5N}( \alpha\!+\!N^2B_{\max}^2 \beta) \frac{1}{(1\!-\!\lambda_{w}^2)^4} 2(b_5\!+\!b_4d_5)\beta^2
     \notag\\
         &~~~~ ~~ +  \frac{  L^2}{5N^3 } (\alpha \!+ \!N^2B_{\max}^2 \beta)  \frac{1}{(1\!-\!\lambda_{w}^2)^4}   2a_5\alpha^2
    \bigg\} T\frac{\sigma^2}{|\mathcal{I}|}\notag\\
    &~~~
    + \bigg\{\frac{   L^2    B_{\max}^2}{5N}( \alpha+N^2B_{\max}^2 \beta) \frac{1}{(1-\lambda_{w}^2)^4}  a_4b_3d_1\alpha^2\beta^2
     \notag\\
         &~~~~ ~~ +  \frac{  L^2}{5N^3 } (\alpha + N^2B_{\max}^2 \beta)  \frac{1}{(1-\lambda_{w}^2)^4}   2a_3\alpha^2
    \bigg\}   \phi_z^0\notag\\
    &~~~
    + \bigg\{\frac{   L^2    B_{\max}^2}{5N}( \alpha+N^2B_{\max}^2 \beta) \frac{1}{(1-\lambda_{w}^2)^4} 2b_4 \beta^2
      \notag\\
         &~~~~ ~~+  \frac{  L^2}{5N^3 } (\alpha + N^2B_{\max}^2 \beta)  \frac{1}{(1-\lambda_{w}^2)^4}  a_4b_3c_2\alpha^2\beta^6
    \bigg\}   \phi_u^0 \label{eq: ave gr insert sum 1}\\
& \leq \frac{ \alpha E}{4}\sum_{r=0}^{T-1} \mathbb{E}[\|\bar \nabla_{\thetab}^r\|^2]
    + \frac{ \beta E}{4}\sum_{r=0}^{T-1}  \mathbb{E}[\|\bar \nabla_{\xb}^r\|^2]  
        + \bigg\{
            \frac{   L^2    B_{\max}^2}{5N}( \alpha+N^2B_{\max}^2 \beta)    \frac{1}{(1-\lambda_{w}^2)^4}   2(b_7+b_4d_7) \beta^4\notag\\
        & ~~~~~~ 
            + \frac{  L^2}{5N^3 } (\alpha + N^2B_{\max}^2 \beta)
            \frac{1}{(1-\lambda_{w}^2)^4}   2(a_7+a_3c_7)\alpha^2\beta^2
    \bigg\}\sum_{r=0}^{T-1}\mathbb{E}[\|\bar\ub^r\|^2]\notag\\
    & ~~~ 
    + \bigg\{\frac{   L^2    B_{\max}^2}{5N}( \alpha+N^2B_{\max}^2 \beta) \frac{1}{(1-\lambda_{w}^2)^4} 2(b_5+b_4d_5)\beta^2
    \notag\\
         &~~~~ ~~ +  \frac{  L^2}{5N^3 } (\alpha + N^2B_{\max}^2 \beta)  \frac{1}{(1-\lambda_{w}^2)^4}   2a_5\alpha^2
    \bigg\} T\frac{\sigma^2}{|\mathcal{I}|}\notag\\
    &~~~
    +      \frac{  L^2}{5N^3 } (\alpha + N^2B_{\max}^2 \beta)  \frac{1}{(1-\lambda_{w}^2)^4}   2a_3\alpha^2   \phi_z^0
    +  \frac{   L^2    B_{\max}^2}{5N}( \alpha+N^2B_{\max}^2 \beta) \frac{1}{(1-\lambda_{w}^2)^4} 2b_4 \beta^2    \phi_u^0,\label{eq: gr sum theta x HBL 2}
\end{align}
where the last inequality is by letting
\begin{align}
 \alpha \leq \frac{\sqrt{b_4}}{\sqrt{a_4}}, ~ \beta \leq \min \bigg\{\frac{\sqrt{a_3}}{\sqrt{a_4b_3d_1B_{\max}^2 N^2}}, \frac{\sqrt{b_4}}{\sqrt{b_3c_2N^2}}\bigg\},
\end{align}
so that in \eqref{eq: ave gr insert sum 1}, the first term is smaller than the second term before $\phi_z^0$, and the first term is smaller than the second term before $\phi_u^0$.

By inserting \eqref{eq: sum bar theta x HBL} into the above inequality, we get
\begin{align}
&\sum_{r=0}^{T-1}\frac{\alpha E}{10} \frac{1}{N} \sum_{n=1}^N \mathbb{E}\bigg[\bigg\| \frac{1}{NS}\sum_{i=1}^S\nabla_{\thetab} f\bigg(\sum_{t=1}^N \Bb_{t,i}\xb_n^r, \thetab_n^r\bigg)\bigg\|^2\bigg] \notag\\
&~~~ + \sum_{r=0}^{T-1} \frac{\beta E}{10}\frac{1}{N} \sum_{n=1}^N\mathbb{E}\bigg[\bigg\| \frac{1}{NS}\sum_{i=1}^S\sum_{t=1}^N \Bb_{t,i}^\top \nabla_{\Bb_i\xb} f\bigg(\sum_{t=1}^N \Bb_{t,i}\xb_n^r, \thetab_n^r)\bigg\|^2\bigg] \notag\\
& \leq   \tilde{F}(\bar \thetab^0, \bar\xb^{0} ) - \underline{F}  +
         \bigg\{
            \frac{   L^2    B_{\max}^2}{5N}( \alpha+N^2B_{\max}^2 \beta)    \frac{1}{(1-\lambda_{w}^2)^4}   2(b_7+b_4d_7) \beta^4\notag\\
        & ~~~~~~ 
            + \frac{  L^2}{5N^3 } (\alpha + N^2B_{\max}^2 \beta)
            \frac{1}{(1-\lambda_{w}^2)^4}   2(a_7+a_3c_7)\alpha^2\beta^2
           - \frac{\beta}{4}
    \bigg\}\sum_{r=0}^{T-1}\mathbb{E}[\|\bar\ub^r\|^2]\notag\\
    &~~~ 
        + [4EL\alpha^2  +10(E-1)L\beta^2]T\frac{\sigma^2}{|\mathcal{I}|} \notag\\
    &~~~ 
        + \bigg\{  \bigg( \frac{2\alpha L^2}{ N^3} +\frac{2\beta B_{\max}^2 L^2}{N } \bigg)\frac{ 1 }{(1-\lambda_{w}^2)^4}
         \bigg[   [1+ 7(E-1)\lambda_w^2] 2 a_5\alpha^2
            + \frac{1}{S}[1+32(E-1)\lambda_w^2] 2c_5\beta^2
    \notag\\
    &~~~ ~~~~~~
        +  [E + 200(E-1) \lambda_w^2]4b_4d_5 N^2B_{\max}^2\beta^2
       +  128(E-1) [1+ E(E-1)]\lambda_w^2 2d_5  N^2B_{\max}^2\beta^2
        \bigg]
        \notag\\
        &~~~ ~~~~~~+   \frac{(\alpha + N^2B_{\max}^2 \beta) L^2}{5N^3(1-\lambda_{w}^2)^4} \bigg[  N^2 B_{\max}^2  2(b_5+b_4d_5)\beta^2
    +    2a_5\alpha^2
    \bigg]
        \bigg\}T\frac{\sigma^2}{|\mathcal{I}|}  \notag\\
    &~~~
        + \bigg\{\bigg( \frac{2\alpha L^2}{ SN^3} +\frac{2\beta B_{\max}^2 L^2}{SN } \bigg) [1+ 32(E-1)\lambda_w^2]  \frac{2}{(1-\lambda_{w}^2)^4}  \notag\\
         &~~~~ ~~
        +  \frac{  L^2}{5N^3 } (\alpha + N^2B_{\max}^2 \beta)  \frac{1}{(1-\lambda_{w}^2)^4}   2a_3\alpha^2
        \bigg\}\phi_z^0 \notag\\
    &~~~
        +\bigg\{ \bigg[ \bigg( \frac{2\alpha L^2}{ N^3} \!+\!\frac{2\beta B_{\max}^2 L^2}{N } \bigg)   \frac{  N^2B_{\max}^2 \beta^2  }{(1\!-\!\lambda_{w}^2)^4}   \bigg\{2b_4 [E \!+\! 200(E\!-\!1) \lambda_w^2] \!+\!1
        \! \notag\\
         &~~~~ ~~+\! 128(E\!-\!1) [1\!+\! E(E\!-\!1) ]\lambda_w^2
         \bigg]
          +
            \frac{   L^2    B_{\max}^2}{5N}( \alpha+N^2B_{\max}^2 \beta) \frac{1}{(1-\lambda_{w}^2)^4} 2b_4 \beta^2
         \bigg\}
         \phi_u^0 \label{eq: ave gr sum 2}\\
& \leq \tilde{F}(\bar \thetab^0, \bar\xb^{0} ) - \underline{F} \notag\\
    &~~~ 
        + [4EL\alpha^2  +10(E-1)L\beta^2]T\frac{\sigma^2}{|\mathcal{I}|} \notag\\
    &~~~ 
        +  \frac{(\alpha + N^2B_{\max}^2 \beta) L^2}{5N^3(1-\lambda_{w}^2)^4} \bigg\{
             [1+ 7(E-1)\lambda_w^2] 4 a_5\alpha^2
            + \frac{1}{S}[1+32(E-1)\lambda_w^2] 4c_5\beta^2
    \notag\\
    &~~~ ~~~
        +  [E + 200(E-1) \lambda_w^2]8b_4d_5 N^2B_{\max}^2\beta^2
       +  128(E-1) [1+ E(E-1)]\lambda_w^2 4d_5  N^2B_{\max}^2\beta^2
        \notag\\
        &~~~~~~
        +     N^2 B_{\max}^2  2(b_5+b_4d_5)\beta^2
         +    2a_5\alpha^2
        \bigg\}T\frac{\sigma^2}{|\mathcal{I}|}  \notag\\
    &~~~
        + \frac{6(\alpha + N^2B_{\max}^2 \beta) L^2}{ SN^3(1-\lambda_{w}^2)^4}
            [1+ 32(E-1)\lambda_w^2]  \phi_z^0 \notag\\
    &~~~
        + \frac{(\alpha + N^2B_{\max}^2 \beta) L^2}{5N^3(1-\lambda_{w}^2)^4}
        \bigg\{      5N^2B_{\max}^2 \beta^2    \bigg[ 10b_4 [E + 200(E-1) \lambda_w^2] +2
           \notag\\
         &~~~~~~   +640(E-1) [1+ E(E-1) ]\lambda_w^2\bigg]
          +
             2b_4B_{\max}^2 \beta^2
         \bigg\}
         \phi_u^0 \notag\\
& \leq  \tilde{F}(\bar \thetab^0, \bar\xb^{0} ) - \underline{F}
        + 2EL[2\alpha^2  +5\beta^2]T\frac{\sigma^2}{|\mathcal{I}|} \notag\\
    &~~~ +  (\alpha + N^2B_{\max}^2 \beta)   (327\alpha^2 + 221326N^2B_{\max}^2 \beta^2)
        \frac{  E^3 L^2(1+7\lambda_w^2)^2  T \sigma^2}
            { N^2 (1-\lambda_{w}^2)^6 |\mathcal{I}|}  \notag\\
    &~~~ +  \frac{(\alpha + N^2B_{\max}^2 \beta)198 EL^2 }
            {S N^2(1-\lambda_{w}^2)^4}\phi_z^0
            +    \frac{(\alpha + N^2B_{\max}^2 \beta) E^3 L^2(1+7\lambda_w^2)   58681 B_{\max}^2 \beta^2}
            { N^2(1-\lambda_{w}^2)^5  } \phi_u^0
, \label{eq: gr sum theta x HBL 3}
\end{align}
where the last two inequality is by letting
\begin{align}
&\alpha \leq \frac{\sqrt{2 [1+32(E-1)\lambda_w^2]}}{\sqrt{5S a_3}},
\end{align}
so that in \eqref{eq: ave gr sum 2}, the second coefficient is smaller than the first coefficient before $\phi_z^0$, the last inequality is by inserting $a_5, b_4,b_5,c_5,d_5$ referred in Lemma \ref{lem: consensus HBL}
and Lemma \ref{lem: tracking HBL}, in addtion to using the properties of $\lambda_w^2< 1$ and $E\geq 1$ .

Let $\eta \triangleq \min\{\alpha, \beta\}$, then by dividing $\textstyle \frac{\eta E T}{  10}$ on both sides of \eqref{eq: gr sum theta x HBL 3}, we have
\begin{align}\label{eq: gr sum theta x HBL 4}
&\frac{1}{T}\sum_{r=0}^{T-1} \frac{1}{N} \sum_{n=1}^N \bigg\{ \mathbb{E}\bigg[\bigg\| \frac{1}{NS}\sum_{i=1}^S\nabla_{\thetab} f\bigg(\sum_{t=1}^N \Bb_{t,i}\xb_n^r, \thetab_n^r\bigg)\bigg\|^2\bigg] \notag\\
&~~~ +  \mathbb{E}\bigg[\bigg\|   \frac{1}{NS}\sum_{i=1}^S\sum_{t=1}^N \Bb_{t,i}^\top  \nabla_{\Bb_i\xb} f\bigg(\sum_{t=1}^N \Bb_{t,i}\xb_n^r, \thetab_n^r)\bigg\|^2\bigg]\bigg\} \notag\\
& \leq \frac{10(\bar{\underline{F}}^0 -\underline{F})}{\eta ET}
        +\frac{200L (2 \alpha^2  +5\beta^2 ) \sigma^2}{\eta E  N |\mathcal{I}| }   \notag\\
    &~~~ + 10 (\alpha + N^2B_{\max}^2 \beta)   (327\alpha^2 + 221326N^2B_{\max}^2 \beta^2)
        \frac{  E^2 L^2(1+7\lambda_w^2)^2    \sigma^2}
            {\eta  N^2 (1-\lambda_{w}^2)^6 |\mathcal{I}|}  \notag\\
    &~~~ +  \frac{(\alpha + N^2B_{\max}^2 \beta)1980  L^2 }
            {\eta TS N^2(1-\lambda_{w}^2)^4}\phi_z^0
    +    \frac{(\alpha + N^2B_{\max}^2 \beta) E^2 L^2(1+7\lambda_w^2)   586810 B_{\max}^2 \beta^2}
            {\eta T N^2(1-\lambda_{w}^2)^5  } \phi_u^0
, .
\end{align}

Without loss of generality, we assume that $N^2B_{\max}^2 \geq 1$, which can be realized by a properly big $N$ and normalizing data to enlarge $B_{\max}$.
Then, by inserting  $\alpha = \beta = \eta = \textstyle \frac{\sqrt{N}}
{100\sqrt{ET}}$ and $E\leq \frac{T^{\frac{1}{3}}}{N^3}$ into \eqref{eq: gr sum theta x HBL 4}, we obtain
\begin{align}
&\frac{1}{T}\sum_{r=0}^{T-1} \frac{1}{N} \sum_{n=1}^N \bigg\{ \mathbb{E}\bigg[\bigg\| \frac{1}{NS}\sum_{i=1}^S\nabla_{\thetab} f\bigg(\sum_{t=1}^N \Bb_{t,i}\xb_n^r, \thetab_n^r\bigg)\bigg\|^2\bigg] \notag\\
&~~~ +  \mathbb{E}\bigg[\bigg\|   \frac{1}{NS}\sum_{i=1}^S\sum_{t=1}^N \Bb_{t,i}^\top  \nabla_{\Bb_i\xb} f\bigg(\sum_{t=1}^N \Bb_{t,i}\xb_n^r, \thetab_n^r)\bigg\|^2\bigg]\bigg\} \notag\\
& \leq \frac{1000[\tilde{F}(\bar \thetab^0, \!\bar\xb^0) \! -\!\underline{F}]}{\sqrt{N ET}}
        \!+\!\frac{14L   \sigma^2}{  \sqrt{N ET} |\mathcal{I}| }
       \!  + \!
        \frac{   222 L^2(1\!+\!7\lambda_w^2)^2 B_{\max}^2 \sigma^2}
            {  \sqrt{N ET}  N (1\!-\!\lambda_{w}^2)^6 |\mathcal{I}|}  \notag\\
    &~~~ +  \frac{2L^2 B_{\max}^2  }
            {\sqrt{NET} N  (1-\lambda_w^2)^4}
            \bigg[\frac{990}{SE N }\phi_z^0
                + \frac{59 (1+7\lambda_w^2)}{(1-\lambda_w^2) T }
            \phi_u^0
            \bigg]   .
\end{align}
Thus, we have completed the proof Theorem \ref{thm: HBL}.
. \hfill $\blacksquare$
\end{document}